\documentclass[10pt]{article}
\usepackage[all]{xy}
\usepackage{amsmath}
\usepackage{amsfonts}
\usepackage{amssymb}
\usepackage{amscd}
\usepackage{amsthm}
\usepackage{latexsym}
\usepackage{amsbsy}

\begin{document}

\title{The Gauss-Bonnet Theorem for Noncommutative Two Tori With a General Conformal Structure}
\author {Farzad Fathi Zadeh and Masoud Khalkhali
\\ Department of Mathematics,  University of Western Ontario
\\ London, Ontario, Canada}
\date{}
\maketitle
\begin{abstract}
In this paper we give a proof of  the Gauss-Bonnet theorem of Connes and Tretkoff
for  noncommutative two tori $\mathbb{T}_{\theta}^2$
equipped with an arbitrary translation invariant complex structure. More
precisely, we show that for any
complex number $\tau$ in the upper  half plane, representing the conformal class of a metric on
$\mathbb{T}_{\theta}^2$, and a Weyl factor given by a positive invertible element $k \in C^{\infty}(\mathbb{T}_{\theta}^2)$,
the value at the origin,  $\zeta (0)$,  of
the spectral zeta function of the Laplacian $\triangle'$ attached to
$(\mathbb{T}_{\theta}^2, \tau, k)$ is independent of $\tau$ and
$k$.
\end{abstract}

\section{Introduction}
In this paper we extend the Gauss-Bonnet theorem of Connes and Tretkoff
\cite{contre} (cf. also \cite{cohcon} for a preliminary version)
for the noncommutative two torus $\mathbb{T}_{\theta}^2$,  to noncommutative two tori
equipped with an arbitrary translation invariant complex structure.  In fact this more general result was   already
 stated  in \cite{contre}.
 More
precisely, we show that for any
complex number $\tau$ in the upper  half plane, representing the conformal class of a metric on
$\mathbb{T}_{\theta}^2$, and a Weyl factor given by a positive invertible element $k \in C^{\infty}(\mathbb{T}_{\theta}^2)$,
the value at the origin,  $\zeta (0)$,  of
the spectral zeta function of the Laplacian $\triangle'$ attached to
$(\mathbb{T}_{\theta}^2, \tau, k)$ is independent of $\tau$ and
$k$. The main ideas in the proof are those employed and invented in
\cite{contre} to treat the case $\tau=\sqrt{-1}$, namely a pseudodifferential calculus for the canonical dynamical system associated to the noncommutative torus, and
  the use of the asymptotic expansion of the heat kernel in computing
the zeta values. This, however,  by itself is not enough and, similar to
\cite{contre},
one needs an extra and intricate argument to express $\zeta (0)$ in terms of the modular operator defined by the Weyl factor.
As a first step, the calculation of the asymptotic expansion of the heat operator for arbitrary values of the conformal class $\tau$
is quite involved and must be performed
by   a computer. We found it impossible to carry this step without the
use of symbolic calculations. Finally we should mention that, as is explained in \cite{contre},  there is a
close relationship between the subject of this paper and scale
invariance in spectral action \cite{chacon, chacon1}  on the one hand,
and non-unimodular (or twisted) spectral triples \cite{conmos} on the other hand.

It is a great pleasure to thank Alain Connes for motivating and enlightening discussions on the topic of this paper. We would also
like to   thank
German Combariza Gonzales, Farideh Fazayeli, and Eric Schost for their generous help with symbolic calculations with GAP.

\section{Preliminaries} \label{preliminaries}

Let $\Sigma$ be a closed, oriented,  2-dimensional smooth manifold equipped with a Riemannian metric $g$. Recall that the spectral  zeta function $\zeta (s)$
 associated to the Laplacian $\triangle_{g} = d^*d$ of  $(\Sigma, g)$
  is given by the following summation over the non-zero eigenvalues $\lambda_{j}$ of
$\triangle_{g}$:
\begin{equation} \label{zetadef}
\zeta (s) = \sum \lambda_{j}^{-s}, \,\,\, \quad \textnormal{Re}(s) > 1.
\end{equation}
The zeta function has a meromorphic continuation to $\mathbb{C}$ with a
 unique (simple) pole at $s=1$. In particular $\zeta (0)$ is
defined and  we have
\[ \zeta(0) + \textnormal{Card} \{j | \lambda_j = 0\} = \frac{1}{12 \pi} \int_{\Sigma} R = \frac{1}{6}
\chi (\Sigma), \nonumber \]
where $R$ is the scalar curvature and $\chi (\Sigma)$ is the Euler-Poincar\'e characteristic. Thus $\zeta (0)$ is a
topological invariant, and,   in particular, it remains invariant under  the conformal perturbation $g \to e^fg$ of the metric \cite{gil, ros}.

Next, we  recall the definition of the spectral zeta function for the Laplacian on the noncommutative two torus from \cite{contre}. Let  $\theta$ be an irrational number. Recall that the  irrational rotation $C^*$-algebra
 $A_{\theta}$ is, by definition, the universal unital $C^*$-algebra  generated by two unitaries $U, V$ satisfying
 $$VU=e^{2 \pi i \theta} UV.$$
One usually thinks of $A_{\theta}$ as the algebra of continuous functions on the noncommutative 2-torus
$\mathbb{T}_{\theta}$.
There is a continuous action of $\mathbb{T}^2$, $\mathbb{T}= \mathbb{R}/2\pi \mathbb{Z}$, on $A_{\theta}$ by $C^*$-algebra
automorphisms  $\{ \alpha_s\}$, $s\in \mathbb{R}^2$, defined by
\[\alpha_s(U^mV^n)=e^{is.(m,n)}U^mV^n.\]
The space of smooth elements for this action, that is those elements $a \in A_{\theta}$ for which the map $s \mapsto
\alpha_s (a)$  is $C^{\infty}$ will be denoted by
$A_{\theta}^{\infty}$. It is a dense subalgebra of $A_{\theta}$  which can be alternatively
described as the algebra of elements in $A_{\theta}$
whose (noncommutative) Fourier expansion has rapidly decreasing coefficients:
\[A_{\theta}^{\infty}=\big \{\sum_{m,n\in \mathbb{Z}}a_{m,n}U^mV^n: (|m|^k|n|^q|a_{m,n}|)_{m,n}
\text{\, is bounded for any positive\,} k,q  \big \}.\]
There is a unique normalized trace   $\tau_0$  on $A_{\theta}$ whose restriction on smooth elements is given by
\[\tau_0(\sum_{m,n\in \mathbb{Z}}a_{m,n}U^mV^n)=a_{0,0}.\]

 The infinitesimal generators of the above action of $\mathbb{T}^2$ on $A_{\theta}$ are the
 derivations  $\delta_1, \, \delta_2: A_{\theta}^{\infty} \to A_{\theta}^{\infty}$ defined by
\[\delta_1(U)=U, \,\, \delta_1(V)=0, \quad  \delta_2(U)=0, \,\, \delta_2(V)=V.\]
In fact, $\delta_1, \delta_2$ are analogues of the differential operators
$\frac{1}{i}\partial/\partial x, \frac{1}{i}\partial/\partial y$ acting on the smooth functions on
the ordinary two torus. We have $\delta_j(a^*)= -\delta_j(a)^* $ for $j=1, 2$ and all $a\in A_{\theta}^{\infty}$.
 Moreover, since $\tau_0 \circ \delta_j =0$, for $j=1, 2$, we have the
 analogue of the
integration by parts formula:
\[ \tau_0(a\delta_j(b)) = -\tau_0(\delta_j(a)b), \,\,\, \forall a,b \in A_{\theta}^{\infty}. \]

We define an inner product on  $A_{\theta}$ by
\[ \langle a, b \rangle = \tau_0(b^*a), \,\,\, a,b \in A_{\theta}, \nonumber \]
and complete $A_{\theta}$ with respect to this inner product to obtain a Hilbert space denoted by
$\mathcal{H}_0$. The derivations $\delta_1, \delta_2$, as unbounded
operators on $\mathcal{H}_0$,  are formally selfadjoint and have unique extensions to selfadjoint operators.

We introduce a complex structure associated with a complex number $\tau = \tau_1+i\tau_2, \, \tau_2 >0,$
by defining
\[ \partial = \delta_1 + \tau \delta_2, \,\,\, \partial^*=  \delta_1 + \overline{\tau} \delta_2. \]
To the conformal structure defined by $\tau$,  corresponds a  positive
Hochschild  two cocycle on  $A_{\theta}^{\infty}$ given by (cf.
\cite{con1})
\[ \psi (a, b, c)=- \tau_0 (a\partial b \partial^* c).\]

We note that   $\partial$ is an unbounded operator on  $\mathcal{H}_0$ and $\partial^*$ is
 its formal adjoint. The analogue of the space of $(1,0)-$forms on the ordinary two torus
is defined to be the Hilbert space completion of the space of finite sums $\sum a \partial b$, $a,b
\in A_{\theta}^{\infty}$, with respect to the inner product defined above, and it is denoted by
$\mathcal{H}^{(1,0)}$. Now we view $\partial$ as an unbounded operator from $\mathcal{H}_0$ to
$\mathcal{H}^{(1,0)}$ and define the Laplacian $\triangle$ on $A_{\theta}^{\infty}$ by:
\[\triangle := \partial^* \partial = \delta_1^2 + 2 \tau_1 \delta_1\delta_2 +
|\tau|^2\delta_2^2.  \]

In order to investigate the analogue of the Gauss-Bonnet theorem for the noncommutative
torus \cite{contre}, we vary the conformal class of the metric by choosing a smooth
selfadjoint element $h=h^* \in A_{\theta}^{\infty}$, and define a linear functional
$\varphi$ on $A_{\theta}$ by
\[\varphi(a)= \tau_0(ae^{-h}), \,\,\, a \in A_{\theta}.\]
In fact, $\varphi$ is a positive linear functional which is not a trace, however,
it is a twisted trace, and satisfies the KMS condition at $\beta = 1$ for the
1-parameter group $\{\sigma_t \}$, $t \in \mathbb{R}$ of inner automorphisms
$\sigma_t= \Delta^{-it}$ where the modular operator for $\varphi$ is given by
($cf.$ \cite{contre})
\[\Delta(x)=e^{-h}xe^{h};\]
moreover, the 1-parameter group of automorphisms $\sigma_t$ is generated by the derivation
$- \log \Delta$ where
\[\log \Delta (x) = [-h,x], \,\,\, x \in A_{\theta}^{\infty}. \]

We  define an inner product $\langle \, , \, \rangle_{\varphi}$ on $A_{\theta}$ by
\[ \langle a,b  \rangle_{\varphi} = \varphi(b^*a), \,\,\, a,b \in A_{\theta}. \]
The Hilbert space obtained from completing $A_{\theta}$ with respect to this
inner product will be denoted by $\mathcal{H}_{\varphi}$.
Now let $\partial_{\varphi}$ be the same operator as  $\partial$, but viewed as an unbounded
operator from $\mathcal{H}_{\varphi}$ to $\mathcal{H}^{(1,0)}$, and define the modified
Laplacian $\triangle'$ by
\[\triangle':=\partial_{\varphi}^* \partial_{\varphi}. \]
Obviously, $\triangle'$ is a positive unbounded operator acting in $\mathcal{H}_{\varphi}.$

\newtheorem{main1}{Lemma}[section]
\begin{main1} \label{main1}
Let $k=e^{h/2} \in A_{\theta}$ act in the Hilbert space $\mathcal{H}_0$ by left
multiplication. Then the operator $\triangle'$ acting in the Hilbert space
$\mathcal{H}_{\varphi}$ is anti-unitarily equivalent to the operator $k\triangle k$
acting in $\mathcal{H}_0$.
\begin{proof}
The proof is very similar to that of Connes and Tretkoff given in Lemma 2.1 in
\cite{contre}. The right multiplication by $k$ extends to a unitary map 
$W:\mathcal{H}_0 \to \mathcal{H}_{\varphi}$ since
\[ \langle W(a), W(b) \rangle_{\varphi} = \tau_0 ((bk)^*(ak)k^{-2})=\tau_0(kb^*akk^{-2})=
\tau_0(b^*a)=\langle a, b\rangle. \]
Obviously $\partial_{\varphi} \circ W= \partial_{\varphi} \circ R_k,$ both acting in
$\mathcal{H}_0$. Therefore $(\partial_{\varphi} \circ W)^*(\partial_{\varphi} \circ W)=
(\partial_{\varphi} \circ R_k)^*(\partial_{\varphi} \circ R_k)$
which, since $R_k^*=R_k$, yields:
\[  W^* \triangle' W = R_k \triangle R_k  . \]
Therefore the operators  $\triangle'$ and $R_k \triangle R_k$, acting in $\mathcal{H}_{\varphi}$
and $\mathcal{H}_0$ respectively, are unitarily equivalent. Now let $J$ be the involution
on $\mathcal{H}_0$ given by $J(a)=a^*$ for all $a\in A_{\theta}$ which is anti-unitary.
Since the operator $J$ commutes with $\triangle= \delta_1^2 + 2 \tau_1 \delta_1\delta_2 +
|\tau|^2\delta_2^2$ and satisfies $JR_k J= k$, we have
\[JR_k\triangle R_kJ= JR_kJ J\triangle J J R_kJ= k \triangle k.\]
Therefore  the operators $\triangle'$ and $k \triangle k$ are anti-unitarily equivalent.
\end{proof}
\end{main1}

\newtheorem{adjoint}[main1]{Remark}
\begin{adjoint}
To have an explicit formula for the adjoint of $\partial_{\varphi}$, one can easily
check that
\[\partial_{\varphi}^*=R_{k^2} \partial^*.\]
\end{adjoint}

\section{Pseudodifferential Operators and $\zeta(0)$} \label{zeta(0)}
As we saw in Section \ref{preliminaries}, the operator
$\triangle' \sim k \triangle k $ is a positive unbounded operator.
Similar to the classical case, using the formula \eqref{zetadef},
one can associate a spectral zeta function to this operator. The main result
of this paper states that the value at the origin of the zeta function of the
operator $\triangle' \sim k \triangle k$ is independent of the choice of the
invertible positive element  $k \in A_{\theta}^{\infty}$ and the conformal class of the metric defined by $\tau \in \mathbb{C} \setminus
\mathbb{R}$. Accordingly, this
section is intended to briefly discuss the notion of a pseudodifferential
operator associated to the dynamical system $(A_{\theta}^{\infty}, \alpha_s)$ \cite{baj, con},
and to  explain how this will help to compute the value at the origin
of the spectral zeta function of an elliptic operator. For a more detailed
discussion, we refer the reader to \cite{contre, gil}.\\

For a non-negative integer $n$, the space of differential operators
on $A_{\theta}^{\infty}$ of order at most $n$ is defined to be the
vector space of operators of the form
\[ \sum_{j_1+j_2 \leq n } a_{j_1,j_2}
\delta_1^{j_1} \delta_2^{j_2}, \,\,\, j_1, \, j_2 \geq 0, \,\,\, a_{j_1,j_2}
\in A_{\theta}^{\infty}.\]

The notion of a differential operator on $A_{\theta}^{\infty}$ can be generalized
to the notion of a pseudodifferential operator using operator valued symbols. In the
sequel, we shall use the notation
$\partial_1=\frac{\partial}{\partial \xi_1}$, $\partial_2=\frac{\partial}{\partial \xi_2}$.

\newtheorem{symboldef}[main1]{Definition}
\begin{symboldef}
For an integer $n$, a smooth map $\rho: \mathbb{R}^2 \to A_{\theta}^{\infty}$ is  said
to be a symbol of order $n$, if for all non-negative integers $i_1, i_2, j_1,
j_2,$
\[ ||\delta_1^{i_1} \delta_2^{i_2} \partial_1^{j_1} \partial_2^{j_2} \rho(\xi) ||
\leq c (1+|\xi|)^{n-j_1-j_2},\]
where $c$ is a constant, and if there exists a smooth map $k: \mathbb{R}^2 \to
A_{\theta}^{\infty}$ such that
\[\lim_{\lambda \to \infty} \lambda^{-n} \rho(\lambda\xi_1, \lambda\xi_2) = k (\xi_1, \xi_2).\]
The space of symbols of order $n$ is denoted by $S_n$.
\end{symboldef}

To a symbol $\rho$ of order $n$, one can associate an operator on $A_{\theta}^{\infty}$,
denoted by $P_{\rho}$, given by
\[ P_{\rho}(a) = (2 \pi)^{-2} \int \int e^{-is \cdot \xi} \rho(\xi) \alpha_s(a) \,ds \,
d\xi. \]
The operator $P_{\rho}$ is said to be a pseudodifferential operator of order $n$. For
example, the differential operator $\sum_{j_1+j_2 \leq n } a_{j_1,j_2}
\delta_1^{j_1} \delta_2^{j_2}$ is associated with the symbol $\sum_{j_1+j_2 \leq n } a_{j_1,j_2}
\xi_1^{j_1} \xi_2^{j_2}$ via the above formula.

\newtheorem{symbolequivalence}[main1]{Definition}
\begin{symbolequivalence}
Two symbols $\rho$, $\rho'\in S_k$  are said to be equivalent if and only if $\rho-\rho'$ is in
$S_n$ for all integers $n$. The equivalence of the symbols will be denoted by  $\rho \sim \rho'$.
\end{symbolequivalence}

The following lemma shows that the space of pseudodifferential operators
is an algebra and one can find the symbol of the product of
pseudodifferential operators up to the above equivalence. Also,
the adjoint of a pseudodifferential operator, with respect to the inner
product defined on $\mathcal{H}_0$ in Section
\ref{preliminaries}, is a pseudodifferential operator with the symbol given in the following proposition up to the equivalence ($cf.$ \cite{contre}).

\newtheorem{symbolcalculus}[main1]{Proposition}
\begin{symbolcalculus} \label{symbolcalculus}
Let $P$ and $Q$ be pseudodifferential operators with the symbols
$\rho$ and $\rho'$ respectively. Then the adjoint $P^*$ and
the product $PQ$ are pseudodifferential operators with the following
symbols
\[
\sigma(P^*) \sim \sum_{\ell_1, \ell_2 \geq 0} \frac{1}{\ell_1! \ell_2!}
\partial_1^{\ell_1} \partial_2^{\ell_2}\delta_1^{\ell_1}\delta_2^{\ell_2}
(\rho(\xi))^*,
\]

\[
\sigma (P Q) \sim \sum_{\ell_1, \ell_2 \geq 0} \frac{1}{\ell_1! \ell_2!}
\partial_1^{\ell_1} \partial_2^{\ell_2}(\rho (\xi))
\delta_1^{\ell_1}\delta_2^{\ell_2} (\rho'(\xi)).
\]

\end{symbolcalculus}

\newtheorem{elliptic}[main1]{Definition}
\begin{elliptic} \label{elliptic}
Let $\rho$ be a symbol of order $n$. It is said to be elliptic if $\rho(\xi)$ is invertible for $\xi \neq 0$, and if there exists a constant $c$ such that
\[ || \rho(\xi)^{-1} || \leq c (1+|\xi|)^{-n} \]
for sufficiently large $|\xi|.$
\end{elliptic}

The operator $\triangle= \delta_1^2 + 2 \tau_1 \delta_1\delta_2 +
|\tau|^2\delta_2^2$ defined in Section \ref{preliminaries} is an example of an elliptic operator.\\

In order to study the value at the origin of the $\zeta$ function of the
operator $\triangle' \sim k \triangle k$, one can use the Mellin transform to
write:

\[ \zeta(s) = \text{Tr} (\triangle'^{-s}) = \frac{1}{\Gamma(s)}
\int_{0}^{\infty} \text{Trace}^+ (e^{-t\triangle'})t^{s-1} \, dt, \]
where $\text{Trace}^+ (e^{-t\triangle'})= \text{Trace} (e^{-t\triangle'}) -
\text{Dim Ker} (\triangle')$.

Using the Cauchy integral formula, one has
\[ e^{-t\triangle'} = \frac{1}{2\pi i} \int_C e^{-t \lambda} (\triangle' - \lambda)^{-1} \, d \lambda \]
where $C$ is a curve in the complex plane that goes around the non-negative
real axis in the anti-clockwise direction without touching it. Appealing to this formula and  using similar arguments to those of \cite{gil}, one can derive the following asymptotic expansion:
\[ \text{Trace}(e^{-t \triangle'}) \sim t^{-1} \sum_{n=0}^{\infty} B_{2n} (\triangle')
t^n, \,\,\, t \to 0^+.  \]
Considering the fact that $\Gamma$ has a simple pole at $s=0$ with residue 1, and using the above asymptotic expansion, one can see that
\[ \zeta(0)= B_2(\triangle').\]
In order to find the latter, similar to the formula in \cite{gil}, one can approximate the inverse of the operator
$(\triangle'-\lambda)$ by a pseudodifferential
operator $B_{\lambda}$ whose symbol has an expansion of the form
\[ b_0(\xi, \lambda) +  b_1(\xi, \lambda) +  b_2(\xi, \lambda) + \cdots \]
where $b_j(\xi, \lambda)$ is a symbol of order $-2-j$, and
\[ \sigma(B_{\lambda} (\triangle'-\lambda) ) \sim 1 .\]
Then, one can see that
\[ B_2(\triangle') = \frac{1}{2 \pi i} \int \int_C e^{-\lambda} \tau_0(b_2(\xi, \lambda)) d \lambda d \xi.  \]
To compute the latter, by a homogeneity argument ($cf.$ \cite{contre}), one can set $\lambda=-1$, and multiply the answer by $-1$. Therefore
\[ \zeta(0)= - \int \tau_0(b_2(\xi,-1)) \, d \xi.\]

\section{The Computation of $\zeta(0)$} \label{computationofzeta(0)}
In order to compute the value at the origin of the zeta function of the operator
$\triangle' \sim k\triangle k$
by the method explained in Section \ref{zeta(0)}, first we find the
symbol of this operator.

\newtheorem{lapsymbol}[main1]{Lemma}
\begin{lapsymbol}
The operator $k\triangle k$ has symbol $a_2(\xi)+a_1(\xi)+
a_0(\xi)$
where
\[  a_2(\xi)= \xi_1^2k^2+|\tau|^2\xi_2^2k^2+2\tau_1\xi_1\xi_2k^2,  \]
\[a_1(\xi)=2\xi_1k\delta_1(k) + 2|\tau|^2\xi_2k\delta_2(k) +
2\tau_1\xi_1k\delta_2(k)+2\tau_1\xi_2k\delta_1(k),\]
\[ a_0(\xi)= k\delta_1^2(k)+ |\tau|^2k\delta_2^2(k) + 2\tau_1k\delta_1\delta_2(k).\]
\begin{proof}
The above expressions are derived by applying the formula for
the symbol of product of
pseudodifferential operators explained in Proposition \ref{symbolcalculus}
to the symbol of
left multiplication by $k$ which is $k$, and the symbol of $\triangle$
which is
$\xi_1^2+2\tau_1\xi_1\xi_2+|\tau|^2\xi_2^2$.
\end{proof}
\end{lapsymbol}

Similar to the calculation of the inverse of the symbol of
$\triangle'-\lambda$ in \cite{gil}, now we can find the symbols
$b_j(\xi)$ for $j=0,1,2,\dots,$ such that each $b_j$ is
of order $-2-j$ and
\[ (b_0+b_1+b_2+\cdots)\,\sigma(\triangle'+1) = (b_0+b_1+b_2+\cdots)
((a_2+1)+a_1+a_0) \sim 1 .\]
Note that, without any loss in the generality for computing $\zeta(0)$,
we have set
$\lambda=-1$. Let $a_2'=a_2+1, a_1'=a_1, a_0'=a_0$. In fact, 1 is treated
as a symbol
of order 2. Then, the above equation yields
\[ \sum_{\substack{j, \ell_1, \ell_2 \geq 0,\\ k=0, 1, 2}} \frac{1}{\ell_1! \ell_2!}
\partial_1^{\ell_1} \partial_2^{\ell_2}(b_j)
\delta_1^{\ell_1}\delta_2^{\ell_2} (a_k') \sim 1.\]
Since $\partial_1^{\ell_1} \partial_2^{\ell_2}(b_j)
\delta_1^{\ell_1}\delta_2^{\ell_2} (a_k')$ is a symbol of order
$-2-j-(\ell_1+\ell_2)+k$, one can decompose the above series into
homogeneous terms of order $-n$ for $n \geq 0$, and write the above equation as
\[\sum_{n\geq 0} \,\,\, \sum_{\substack{j, \ell_1, \ell_2 \geq 0,\, k=0, 1, 2,
\\2+j+\ell_1+\ell_2-k=n}} \frac{1}{\ell_1! \ell_2!}
\partial_1^{\ell_1} \partial_2^{\ell_2}(b_j)
\delta_1^{\ell_1}\delta_2^{\ell_2} (a_k') \sim 1 .\]
Now by comparing the symbols of order $n=0$ on both sides, we have
$b_0a_2'=1$, therefore
\[b_0=a_2'^{-1}=(a_2+1)^{-1}=(\xi_1^2k^2+|\tau|^2\xi_2^2k^2+
2\tau_1\xi_1\xi_2k^2+1)^{-1}.\]
Inductively, one gets
\[b_n a_2' = - \sum_{\substack{2+j+\ell_1+\ell_2-k=n, \\ 0 \leq j <n, \, 0 \leq k \leq 2}}
\frac{1}{\ell_1! \ell_2!}
\partial_1^{\ell_1} \partial_2^{\ell_2}(b_j)
\delta_1^{\ell_1}\delta_2^{\ell_2} (a_k') .\]
If k=2 then $\ell_1+\ell_2>0$, therefore
$\delta_1^{\ell_1}\delta_2^{\ell_2} (a_k')=\delta_1^{\ell_1}\delta_2^{\ell_2}
(a_k)$
in the latter. Hence, for $n > 0$:

\[b_n  = - \sum_{\substack{2+j+\ell_1+\ell_2-k=n, \\ 0 \leq j <n, \, 0 \leq k \leq 2}} \frac{1}
{\ell_1! \ell_2!}
\partial_1^{\ell_1} \partial_2^{\ell_2}(b_j)
\delta_1^{\ell_1}\delta_2^{\ell_2} (a_k) b_0.\]
In particular, we have:

\begin{eqnarray}
b_1&=&
-(b_0a_1b_0  + \partial_1(b_0)\delta_1(a_2)b_0 +
\partial_2(b_0)\delta_2(a_2)b_0), \nonumber
\end{eqnarray}
\begin{eqnarray}
b_2&=&
-(b_0a_0b_0 + b_1a_1b_0 + \partial_1(b_0)\delta_1(a_1)b_0 +
\partial_2(b_0)\delta_2(a_1)b_0 + \partial_1(b_1)\delta_1(a_2)b_0 + \nonumber \\
&&\partial_2(b_1)\delta_2(a_2)b_0 + (1/2)\partial_{11}(b_0)\delta_1^2(a_2)b_0 +
(1/2)\partial_{22}(b_0)\delta_2^2(a_2)b_0  +
\partial_{12}(b_0)\delta_{12}(a_2)b_0). \nonumber
\end{eqnarray}

After computing these terms directly, one finds that $b_2$ is equal to\\
$ $\\
$-b_0k\delta_1^2(k)b_0-
2\tau_1b_0k\delta_1\delta_2(k)b_0
-|\tau|^2b_0k\delta_2^2
(k)b_0+\\
6\xi_1^2b_0^2k^2\delta_1(k)^2b_0+\xi_1^2b_0^2k^
2\delta_1^2(k)b_0k+
5\xi_1^2b_0^2k^3\delta_1^2(k)b_0+\\
2\xi_1^2b_0k\delta_1(k)b_0\delta_1(k)b_0k+
6\xi_1^2b_0k\delta_1(k)b_0k\delta_1(k)b_0+\\
6\tau_1\xi_1^2b_0^2k^2\delta_1(k)\delta_2(k)b_0+
6\tau_1\xi_1^2b_0^2k^2\delta_2(k)\delta_1(k)b_0+\\
2\tau_1\xi_1^2b_0^2k^2\delta_1\delta_2(k)b_0k+
10\tau_1\xi_1^2b_0^2k^3\delta_1\delta_2(k)b_0+\\
2\tau_1\xi_1^2b_0k\delta_1(k)b_0\delta_2(k)b_0k+
6\tau_1\xi_1^2b_0k\delta_1(k)b_0k\delta_2(k)b_0+\\
2\tau_1\xi_1^2b_0k\delta_2(k)b_0\delta_1(k)b_0k+
6\tau_1\xi_1^2b_0k\delta_2(k)b_0k\delta_1(k)b_0+\\
12\tau_1\xi_1\xi_2b_0^2k^2\delta_1(k)^2b_0+
2\tau_1\xi_1\xi_2b_0^2k^2\delta_1^2(k)b_0k+\\
10\tau_1\xi_1\xi_2b_0^2k^3\delta_1^2(k)b_0+
4\tau_1\xi_1\xi_2b_0k\delta_1(k)b_0\delta_1(k)b_0k+\\
12\tau_1\xi_1\xi_2b_0k\delta_1(k)b_0k\delta_1(k)b_0+
4\tau_1^2\xi_1^2b_0^2k^2\delta_2(k)^2b_0+\\
4\tau_1^2\xi_1^2b_0^2k^3\delta_2^2(k)b_0+
4\tau_1^2\xi_1^2b_0k\delta_2(k)b_0k\delta_2(k)b_0+\\
8\tau_1^2\xi_1\xi_2b_0^2k^2\delta_1(k)\delta_2(k)b_0+
8\tau_1^2\xi_1\xi_2b_0^2k^2\delta_2(k)\delta_1(k)b_0+\\
4\tau_1^2\xi_1\xi_2b_0^2k^2\delta_1\delta_2(k)b_0k+
12\tau_1^2\xi_1\xi_2b_0^2k^3\delta_1\delta_2(k)b_0+\\
4\tau_1^2\xi_1\xi_2b_0k\delta_1(k)b_0\delta_2(k)b_0k+
8\tau_1^2\xi_1\xi_2b_0k\delta_1(k)b_0k\delta_2(k)b_0+\\
4\tau_1^2\xi_1\xi_2b_0k\delta_2(k)b_0\delta_1(k)b_0k+
8\tau_1^2\xi_1\xi_2b_0k\delta_2(k)b_0k\delta_1(k)b_0+\\
4\tau_1^2\xi_2^2b_0^2k^2\delta_1(k)^2b_0+
4\tau_1^2\xi_2^2b_0^2k^3\delta_1^2(k)b_0+\\
4\tau_1^2\xi_2^2b_0k\delta_1(k)b_0k\delta_1(k)b_0+
2|\tau|^2\xi_1^2b_0^2k^2\delta_2(k)^2b_0+\\
|\tau|^2\xi_1^2b_0^2k^2\delta_2^2(k)b_0k+
|\tau|^2\xi_1^2b_0^2k^3\delta_2^2(k)b_0+\\
2|\tau|^2\xi_1^2b_0k\delta_2(k)b_0\delta_2(k)b_0k+
2|\tau|^2\xi_1^2b_0k\delta_2(k)b_0k\delta_2(k)b_0+\\
4|\tau|^2\xi_1\xi_2b_0^2k^2\delta_1(k)\delta_2(k)b_0+
4|\tau|^2\xi_1\xi_2b_0^2k^2\delta_2(k)\delta_1(k)b_0+\\
8|\tau|^2\xi_1\xi_2b_0^2k^3\delta_1\delta_2(k)b_0+
4|\tau|^2\xi_1\xi_2b_0k\delta_1(k)b_0k\delta_2(k)b_0+\\
4|\tau|^2\xi_1\xi_2b_0k\delta_2(k)b_0k\delta_1(k)b_0+
2|\tau|^2\xi_2^2b_0^2k^2\delta_1(k)^2b_0+\\
|\tau|^2\xi_2^2b_0^2k^2\delta_1^2(k)b_0k+
|\tau|^2\xi_2^2b_0^2k^3\delta_1^2(k)b_0+\\
2|\tau|^2\xi_2^2b_0k\delta_1(k)b_0\delta_1(k)b_0k+
2|\tau|^2\xi_2^2b_0k\delta_1(k)b_0k\delta_1(k)b_0+\\
12\tau_1|\tau|^2\xi_1\xi_2b_0^2k^2\delta_2(k)^2b_0+
2\tau_1|\tau|^2\xi_1\xi_2b_0^2k^2\delta_2^2(k)b_0k+\\
10\tau_1|\tau|^2\xi_1\xi_2b_0^2k^3\delta_2^2(k)b_0+
4\tau_1|\tau|^2\xi_1\xi_2b_0k\delta_2(k)b_0\delta_2(k)b_0k+\\
12\tau_1|\tau|^2\xi_1\xi_2b_0k\delta_2(k)b_0k\delta_2(k)b_0+
6\tau_1|\tau|^2\xi_2^2b_0^2k^2\delta_1(k)\delta_2(k)b_0+\\
6\tau_1|\tau|^2\xi_2^2b_0^2k^2\delta_2(k)\delta_1(k)b_0+
2\tau_1|\tau|^2\xi_2^2b_0^2k^2\delta_1\delta_2(k)b_0k+\\
10\tau_1|\tau|^2\xi_2^2b_0^2k^3\delta_1\delta_2(k)b_0+
2\tau_1|\tau|^2\xi_2^2b_0k\delta_1(k)b_0\delta_2(k)b_0k+\\
6\tau_1|\tau|^2\xi_2^2b_0k\delta_1(k)b_0k\delta_2(k)b_0+
2\tau_1|\tau|^2\xi_2^2b_0k\delta_2(k)b_0\delta_1(k)b_0k+\\
6\tau_1|\tau|^2\xi_2^2b_0k\delta_2(k)b_0k\delta_1(k)b_0+
6|\tau|^4\xi_2^2b_0^2k^2\delta_2(k)^2b_0+\\
|\tau|^4\xi_2^2b_0^2k^2\delta_2^2(k)b_0k+
5|\tau|^4\xi_2^2b_0^2k^3\delta_2^2(k)b_0+\\
2|\tau|^4\xi_2^2b_0k\delta_2(k)b_0\delta_2(k)b_0k+
6|\tau|^4\xi_2^2b_0k\delta_2(k)b_0k\delta_2(k)b_0
-\\8\xi_1^4b_0^3k^4\delta_1(k)^2b_0
-4\xi_1^4b_0^3k^4\delta_1^2(k)b_0k
-\\4\xi_1^4b_0^3k^5\delta_1^2(k)b_0-
6\xi_1^4b_0^2k^2\delta_1(k)b_0k\delta_1(k)b_0k
-\\10\xi_1^4b_0^2k^2\delta_1(k)b_0k^2\delta_1(k)b_0
-10\xi_1^4b_0^2k^3\delta_1(k)b_0\delta_1(k)b_0k-\\
14\xi_1^4b_0^2k^3\delta_1(k)b_0k\delta_1(k)b_0
-4\xi_1^4b_0k\delta_1(k)b_0^2k^2\delta_1(k)b_0k
-\\4\xi_1^4b_0k\delta_1(k)b_0^2k^3\delta_1(k)b_0-
8\tau_1\xi_1^4b_0^3k^4\delta_1(k)\delta_2(k)b_0
-\\8\tau_1\xi_1^4b_0^3k^4\delta_2(k)\delta_1(k)b_0
-8\tau_1\xi_1^4b_0^3k^4\delta_1\delta_2(k)b_0k-\\
8\tau_1\xi_1^4b_0^3k^5\delta_1\delta_2(k)b_0-
6\tau_1\xi_1^4b_0^2k^2\delta_1(k)b_0k\delta_2(k)b_0k
-\\10\tau_1\xi_1^4b_0^2k^2\delta_1(k)b_0k^2\delta_2(k)b_0-
6\tau_1\xi_1^4b_0^2k^2\delta_2(k)b_0k\delta_1(k)b_0k
-\\10\tau_1\xi_1^4b_0^2k^2\delta_2(k)b_0k^2\delta_1(k)b_0
-10\tau_1\xi_1^4b_0^2k^3\delta_1(k)b_0\delta_2(k)b_0k-\\
14\tau_1\xi_1^4b_0^2k^3\delta_1(k)b_0k\delta_2(k)b_0
-10\tau_1\xi_1^4b_0^2k^3\delta_2(k)b_0\delta_1(k)b_0k
-\\14\tau_1\xi_1^4b_0^2k^3\delta_2(k)b_0k\delta_1(k)b_0-
4\tau_1\xi_1^4b_0k\delta_1(k)b_0^2k^2\delta_2(k)b_0k
-\\4\tau_1\xi_1^4b_0k\delta_1(k)b_0^2k^3\delta_2(k)b_0
-4\tau_1\xi_1^4b_0k\delta_2(k)b_0^2k^2\delta_1(k)b_0k-\\
4\tau_1\xi_1^4b_0k\delta_2(k)b_0^2k^3\delta_1(k)b_0-
32\tau_1\xi_1^3\xi_2b_0^3k^4\delta_1(k)^2b_0-\\
16\tau_1\xi_1^3\xi_2b_0^3k^4\delta_1^2(k)b_0k-
16\tau_1\xi_1^3\xi_2b_0^3k^5\delta_1^2(k)b_0
-\\24\tau_1\xi_1^3\xi_2b_0^2k^2\delta_1(k)b_0k\delta_1(k)b_0k-
40\tau_1\xi_1^3\xi_2b_0^2k^2\delta_1(k)b_0k^2\delta_1(k)b_0-\\
40\tau_1\xi_1^3\xi_2b_0^2k^3\delta_1(k)b_0\delta_1(k)b_0k-
56\tau_1\xi_1^3\xi_2b_0^2k^3\delta_1(k)b_0k\delta_1(k)b_0-\\
16\tau_1\xi_1^3\xi_2b_0k\delta_1(k)b_0^2k^2\delta_1(k)b_0k-
16\tau_1\xi_1^3\xi_2b_0k\delta_1(k)b_0^2k^3\delta_1(k)b_0
-\\8\tau_1^2\xi_1^4b_0^3k^4\delta_2(k)^2b_0-
4\tau_1^2\xi_1^4b_0^3k^4\delta_2^2(k)b_0k-\\
4\tau_1^2\xi_1^4b_0^3k^5\delta_2^2(k)b_0-
4\tau_1^2\xi_1^4b_0^2k^2\delta_2(k)b_0k\delta_2(k)b_0k-\\
8\tau_1^2\xi_1^4b_0^2k^2\delta_2(k)b_0k^2\delta_2(k)b_0-
8\tau_1^2\xi_1^4b_0^2k^3\delta_2(k)b_0\delta_2(k)b_0k-\\
12\tau_1^2\xi_1^4b_0^2k^3\delta_2(k)b_0k\delta_2(k)b_0
-4\tau_1^2\xi_1^4b_0k\delta_2(k)b_0^2k^2\delta_2(k)b_0k-\\
4\tau_1^2\xi_1^4b_0k\delta_2(k)b_0^2k^3\delta_2(k)b_0
-24\tau_1^2\xi_1^3\xi_2b_0^3k^4\delta_1(k)\delta_2(k)b_0
-\\24\tau_1^2\xi_1^3\xi_2b_0^3k^4\delta_2(k)\delta_1(k)b_0-
24\tau_1^2\xi_1^3\xi_2b_0^3k^4\delta_1\delta_2(k)b_0k-\\
24\tau_1^2\xi_1^3\xi_2b_0^3k^5\delta_1\delta_2(k)b_0
-20\tau_1^2\xi_1^3\xi_2b_0^2k^2\delta_1(k)b_0k\delta_2(k)b_0k-\\
32\tau_1^2\xi_1^3\xi_2b_0^2k^2\delta_1(k)b_0k^2\delta_2(k)b_0
-20\tau_1^2\xi_1^3\xi_2b_0^2k^2\delta_2(k)b_0k\delta_1(k)b_0k
-\\32\tau_1^2\xi_1^3\xi_2b_0^2k^2\delta_2(k)b_0k^2\delta_1(k)b_0-
32\tau_1^2\xi_1^3\xi_2b_0^2k^3\delta_1(k)b_0\delta_2(k)b_0k-\\
44\tau_1^2\xi_1^3\xi_2b_0^2k^3\delta_1(k)b_0k\delta_2(k)b_0
-32\tau_1^2\xi_1^3\xi_2b_0^2k^3\delta_2(k)b_0\delta_1(k)b_0k-\\
44\tau_1^2\xi_1^3\xi_2b_0^2k^3\delta_2(k)b_0k\delta_1(k)b_0
-12\tau_1^2\xi_1^3\xi_2b_0k\delta_1(k)b_0^2k^2\delta_2(k)b_0k
-\\12\tau_1^2\xi_1^3\xi_2b_0k\delta_1(k)b_0^2k^3\delta_2(k)b_0-
12\tau_1^2\xi_1^3\xi_2b_0k\delta_2(k)b_0^2k^2\delta_1(k)b_0k
-\\12\tau_1^2\xi_1^3\xi_2b_0k\delta_2(k)b_0^2k^3\delta_1(k)b_0-
40\tau_1^2\xi_1^2\xi_2^2b_0^3k^4\delta_1(k)^2b_0-\\
20\tau_1^2\xi_1^2\xi_2^2b_0^3k^4\delta_1^2(k)b_0k-
20\tau_1^2\xi_1^2\xi_2^2b_0^3k^5\delta_1^2(k)b_0-\\
28\tau_1^2\xi_1^2\xi_2^2b_0^2k^2\delta_1(k)b_0k\delta_1(k)b_0k-
48\tau_1^2\xi_1^2\xi_2^2b_0^2k^2\delta_1(k)b_0k^2\delta_1(k)b_0-\\
48\tau_1^2\xi_1^2\xi_2^2b_0^2k^3\delta_1(k)b_0\delta_1(k)b_0k-
68\tau_1^2\xi_1^2\xi_2^2b_0^2k^3\delta_1(k)b_0k\delta_1(k)b_0-\\
20\tau_1^2\xi_1^2\xi_2^2b_0k\delta_1(k)b_0^2k^2\delta_1(k)b_0k
-20\tau_1^2\xi_1^2\xi_2^2b_0k\delta_1(k)b_0^2k^3\delta_1(k)b_0
-\\2|\tau|^2\xi_1^4b_0^2k^2\delta_2(k)b_0k\delta_2(k)b_0k-
2|\tau|^2\xi_1^4b_0^2k^2\delta_2(k)b_0k^2\delta_2(k)b_0
-\\2|\tau|^2\xi_1^4b_0^2k^3\delta_2(k)b_0\delta_2(k)b_0k
-2|\tau|^2\xi_1^4b_0^2k^3\delta_2(k)b_0k\delta_2(k)b_0-\\
8|\tau|^2\xi_1^3\xi_2b_0^3k^4\delta_1(k)\delta_2(k)b_0-8
|\tau|^2\xi_1^3\xi_2b_0^3k^4\delta_2(k)\delta_1(k)b_0-\\
8|\tau|^2\xi_1^3\xi_2b_0^3k^4\delta_1\delta_2(k)b_0k-
8|\tau|^2\xi_1^3\xi_2b_0^3k^5\delta_1\delta_2(k)b_0-\\4
|\tau|^2\xi_1^3\xi_2b_0^2k^2\delta_1(k)b_0k\delta_2(k)b_0k-8
|\tau|^2\xi_1^3\xi_2b_0^2k^2\delta_1(k)b_0k^2\delta_2(k)b_0-\\
4|\tau|^2\xi_1^3\xi_2b_0^2k^2\delta_2(k)b_0k\delta_1(k)b_0k-
8|\tau|^2\xi_1^3\xi_2b_0^2k^2\delta_2(k)b_0k^2\delta_1(k)b_0-\\8
|\tau|^2\xi_1^3\xi_2b_0^2k^3\delta_1(k)b_0\delta_2(k)b_0k-
12|\tau|^2\xi_1^3\xi_2b_0^2k^3\delta_1(k)b_0k\delta_2(k)b_0
-\\8|\tau|^2\xi_1^3\xi_2b_0^2k^3\delta_2(k)b_0\delta_1(k)b_0k-
12|\tau|^2\xi_1^3\xi_2b_0^2k^3\delta_2(k)b_0k\delta_1(k)b_0-\\
4|\tau|^2\xi_1^3\xi_2b_0k\delta_1(k)b_0^2k^2\delta_2(k)b_0k-4
|\tau|^2\xi_1^3\xi_2b_0k\delta_1(k)b_0^2k^3\delta_2(k)b_0-\\
4|\tau|^2\xi_1^3\xi_2b_0k\delta_2(k)b_0^2k^2\delta_1(k)b_0k-
4|\tau|^2\xi_1^3\xi_2b_0k\delta_2(k)b_0^2k^3\delta_1(k)b_0-\\8
|\tau|^2\xi_1^2\xi_2^2b_0^3k^4\delta_1(k)^2b_0-
4|\tau|^2\xi_1^2\xi_2^2b_0^3k^4\delta_1^2(k)b_0k-\\
4|\tau|^2\xi_1^2\xi_2^2b_0^3k^5\delta_1^2(k)b_0
-8|\tau|^2\xi_1^2\xi_2^2b_0^2k^2\delta_1(k)b_0k\delta_1(k)b_0k-\\
12|\tau|^2\xi_1^2\xi_2^2b_0^2k^2\delta_1(k)b_0k^2\delta_1(k)b_0-12
|\tau|^2\xi_1^2\xi_2^2b_0^2k^3\delta_1(k)b_0\delta_1(k)b_0k-\\
16|\tau|^2\xi_1^2\xi_2^2b_0^2k^3\delta_1(k)b_0k\delta_1(k)b_0-
4|\tau|^2\xi_1^2\xi_2^2b_0k\delta_1(k)b_0^2k^2\delta_1(k)b_0k
-\\4|\tau|^2\xi_1^2\xi_2^2b_0k\delta_1(k)b_0^2k^3\delta_1(k)b_0-
16\tau_1^3\xi_1^3\xi_2b_0^3k^4\delta_2(k)^2b_0-\\
8\tau_1^3\xi_1^3\xi_2b_0^3k^4\delta_2^2(k)b_0k
-8\tau_1^3\xi_1^3\xi_2b_0^3k^5\delta_2^2(k)b_0-\\
8\tau_1^3\xi_1^3\xi_2b_0^2k^2\delta_2(k)b_0k\delta_2(k)b_0k-
16\tau_1^3\xi_1^3\xi_2b_0^2k^2\delta_2(k)b_0k^2\delta_2(k)b_0
-\\16\tau_1^3\xi_1^3\xi_2b_0^2k^3\delta_2(k)b_0\delta_2(k)b_0k-
24\tau_1^3\xi_1^3\xi_2b_0^2k^3\delta_2(k)b_0k\delta_2(k)b_0-\\
8\tau_1^3\xi_1^3\xi_2b_0k\delta_2(k)b_0^2k^2\delta_2(k)b_0k-8
\tau_1^3\xi_1^3\xi_2b_0k\delta_2(k)b_0^2k^3\delta_2(k)b_0-\\16
\tau_1^3\xi_1^2\xi_2^2b_0^3k^4\delta_1(k)\delta_2(k)b_0-
16\tau_1^3\xi_1^2\xi_2^2b_0^3k^4\delta_2(k)\delta_1(k)b_0
-\\16\tau_1^3\xi_1^2\xi_2^2b_0^3k^4\delta_1\delta_2(k)b_0k
-16\tau_1^3\xi_1^2\xi_2^2b_0^3k^5\delta_1\delta_2(k)b_0-\\
16\tau_1^3\xi_1^2\xi_2^2b_0^2k^2\delta_1(k)b_0k\delta_2(k)b_0k-
24\tau_1^3\xi_1^2\xi_2^2b_0^2k^2\delta_1(k)b_0k^2\delta_2(k)b_0-\\
16\tau_1^3\xi_1^2\xi_2^2b_0^2k^2\delta_2(k)b_0k\delta_1(k)b_0k-
24\tau_1^3\xi_1^2\xi_2^2b_0^2k^2\delta_2(k)b_0k^2\delta_1(k)b_0-\\
24\tau_1^3\xi_1^2\xi_2^2b_0^2k^3\delta_1(k)b_0\delta_2(k)b_0k-
32\tau_1^3\xi_1^2\xi_2^2b_0^2k^3\delta_1(k)b_0k\delta_2(k)b_0-\\
24\tau_1^3\xi_1^2\xi_2^2b_0^2k^3\delta_2(k)b_0\delta_1(k)b_0k-32
\tau_1^3\xi_1^2\xi_2^2b_0^2k^3\delta_2(k)b_0k\delta_1(k)b_0-\\8
\tau_1^3\xi_1^2\xi_2^2b_0k\delta_1(k)b_0^2k^2\delta_2(k)b_0k-
8\tau_1^3\xi_1^2\xi_2^2b_0k\delta_1(k)b_0^2k^3\delta_2(k)b_0
-\\8\tau_1^3\xi_1^2\xi_2^2b_0k\delta_2(k)b_0^2k^2\delta_1(k)b_0k
-8\tau_1^3\xi_1^2\xi_2^2b_0k\delta_2(k)b_0^2k^3\delta_1(k)b_0
-\\
16\tau_1^3\xi_1\xi_2^3b_0^3k^4\delta_1(k)^2b_0-
8\tau_1^3\xi_1\xi_2^3b_0^3k^4\delta_1^2(k)b_0k
-\\8\tau_1^3\xi_1\xi_2^3b_0^3k^5\delta_1^2(k)b_0-
8\tau_1^3\xi_1\xi_2^3b_0^2k^2\delta_1(k)b_0k\delta_1(k)b_0k
-\\16\tau_1^3\xi_1\xi_2^3b_0^2k^2\delta_1(k)b_0k^2\delta_1(k)b_0
-16\tau_1^3\xi_1\xi_2^3b_0^2k^3\delta_1(k)b_0\delta_1(k)b_0k-\\
24\tau_1^3\xi_1\xi_2^3b_0^2k^3\delta_1(k)b_0k\delta_1(k)b_0
-8\tau_1^3\xi_1\xi_2^3b_0k\delta_1(k)b_0^2k^2\delta_1(k)b_0k-\\
8\tau_1^3\xi_1\xi_2^3b_0k\delta_1(k)b_0^2k^3\delta_1(k)b_0-
16\tau_1|\tau|^2\xi_1^3\xi_2b_0^3k^4\delta_2(k)^2b_0-\\
8\tau_1|\tau|^2\xi_1^3\xi_2b_0^3k^4\delta_2^2(k)b_0k-8
\tau_1|\tau|^2\xi_1^3\xi_2b_0^3k^5\delta_2^2(k)b_0-\\
16\tau_1|\tau|^2\xi_1^3\xi_2b_0^2k^2\delta_2(k)b_0k\delta_2(k)b_0k
-24\tau_1|\tau|^2\xi_1^3\xi_2b_0^2k^2\delta_2(k)b_0k^2\delta_2(k)b_0
-\\24\tau_1|\tau|^2\xi_1^3\xi_2b_0^2k^3\delta_2(k)b_0\delta_2(k)b_0k-
32\tau_1|\tau|^2\xi_1^3\xi_2b_0^2k^3\delta_2(k)b_0k\delta_2(k)b_0-\\
8\tau_1|\tau|^2\xi_1^3\xi_2b_0k\delta_2(k)b_0^2k^2\delta_2(k)b_0k-
8\tau_1|\tau|^2\xi_1^3\xi_2b_0k\delta_2(k)b_0^2k^3\delta_2(k)b_0-\\
32\tau_1|\tau|^2\xi_1^2\xi_2^2b_0^3k^4\delta_1(k)\delta_2(k)b_0
-32\tau_1|\tau|^2\xi_1^2\xi_2^2b_0^3k^4\delta_2(k)\delta_1(k)b_0-\\
32\tau_1|\tau|^2\xi_1^2\xi_2^2b_0^3k^4\delta_1\delta_2(k)b_0k-
32\tau_1|\tau|^2\xi_1^2\xi_2^2b_0^3k^5\delta_1\delta_2(k)b_0-\\
20\tau_1|\tau|^2\xi_1^2\xi_2^2b_0^2k^2\delta_1(k)b_0k\delta_2(k)b_0k-
36\tau_1|\tau|^2\xi_1^2\xi_2^2b_0^2k^2\delta_1(k)b_0k^2\delta_2(k)b_0-\\
20\tau_1|\tau|^2\xi_1^2\xi_2^2b_0^2k^2\delta_2(k)b_0k\delta_1(k)b_0k-
36\tau_1|\tau|^2\xi_1^2\xi_2^2b_0^2k^2\delta_2(k)b_0k^2\delta_1(k)b_0-\\
36\tau_1|\tau|^2\xi_1^2\xi_2^2b_0^2k^3\delta_1(k)b_0\delta_2(k)b_0k-
52\tau_1|\tau|^2\xi_1^2\xi_2^2b_0^2k^3\delta_1(k)b_0k\delta_2(k)b_0-\\
36\tau_1|\tau|^2\xi_1^2\xi_2^2b_0^2k^3\delta_2(k)b_0\delta_1(k)b_0k-
52\tau_1|\tau|^2\xi_1^2\xi_2^2b_0^2k^3\delta_2(k)b_0k\delta_1(k)b_0-\\
16\tau_1|\tau|^2\xi_1^2\xi_2^2b_0k\delta_1(k)b_0^2k^2\delta_2(k)b_0k-
16\tau_1|\tau|^2\xi_1^2\xi_2^2b_0k\delta_1(k)b_0^2k^3\delta_2(k)b_0-\\
16\tau_1|\tau|^2\xi_1^2\xi_2^2b_0k\delta_2(k)b_0^2k^2\delta_1(k)b_0k-
16\tau_1|\tau|^2\xi_1^2\xi_2^2b_0k\delta_2(k)b_0^2k^3\delta_1(k)b_0
-\\16\tau_1|\tau|^2\xi_1\xi_2^3b_0^3k^4\delta_1(k)^2b_0-
8\tau_1|\tau|^2\xi_1\xi_2^3b_0^3k^4\delta_1^2(k)b_0k-\\
8\tau_1|\tau|^2\xi_1\xi_2^3b_0^3k^5\delta_1^2(k)b_0-
16\tau_1|\tau|^2\xi_1\xi_2^3b_0^2k^2\delta_1(k)b_0k\delta_1(k)b_0k
-\\24\tau_1|\tau|^2\xi_1\xi_2^3b_0^2k^2\delta_1(k)b_0k^2\delta_1(k)b_0-
24\tau_1|\tau|^2\xi_1\xi_2^3b_0^2k^3\delta_1(k)b_0\delta_1(k)b_0k-\\32
\tau_1|\tau|^2\xi_1\xi_2^3b_0^2k^3\delta_1(k)b_0k\delta_1(k)b_0-
8\tau_1|\tau|^2\xi_1\xi_2^3b_0k\delta_1(k)b_0^2k^2\delta_1(k)b_0k-\\
8\tau_1|\tau|^2\xi_1\xi_2^3b_0k\delta_1(k)b_0^2k^3\delta_1(k)b_0-
40\tau_1^2|\tau|^2\xi_1^2\xi_2^2b_0^3k^4\delta_2(k)^2b_0-\\
20\tau_1^2|\tau|^2\xi_1^2\xi_2^2b_0^3k^4\delta_2^2(k)b_0k-
20\tau_1^2|\tau|^2\xi_1^2\xi_2^2b_0^3k^5\delta_2^2(k)b_0-\\
28\tau_1^2|\tau|^2\xi_1^2\xi_2^2b_0^2k^2\delta_2(k)b_0k\delta_2(k)b_0k-
48\tau_1^2|\tau|^2\xi_1^2\xi_2^2b_0^2k^2\delta_2(k)b_0k^2\delta_2(k)b_0-\\
48\tau_1^2|\tau|^2\xi_1^2\xi_2^2b_0^2k^3\delta_2(k)b_0\delta_2(k)b_0k-
68\tau_1^2|\tau|^2\xi_1^2\xi_2^2b_0^2k^3\delta_2(k)b_0k\delta_2(k)b_0-\\
20\tau_1^2|\tau|^2\xi_1^2\xi_2^2b_0k\delta_2(k)b_0^2k^2\delta_2(k)b_0k-
20\tau_1^2|\tau|^2\xi_1^2\xi_2^2b_0k\delta_2(k)b_0^2k^3\delta_2(k)b_0-\\
24\tau_1^2|\tau|^2\xi_1\xi_2^3b_0^3k^4\delta_1(k)\delta_2(k)b_0-
24\tau_1^2|\tau|^2\xi_1\xi_2^3b_0^3k^4\delta_2(k)\delta_1(k)b_0
-\\24\tau_1^2|\tau|^2\xi_1\xi_2^3b_0^3k^4\delta_1\delta_2(k)b_0k-
24\tau_1^2|\tau|^2\xi_1\xi_2^3b_0^3k^5\delta_1\delta_2(k)b_0
-\\20\tau_1^2|\tau|^2\xi_1\xi_2^3b_0^2k^2\delta_1(k)b_0k\delta_2(k)b_0k-
32\tau_1^2|\tau|^2\xi_1\xi_2^3b_0^2k^2\delta_1(k)b_0k^2\delta_2(k)b_0
-\\20\tau_1^2|\tau|^2\xi_1\xi_2^3b_0^2k^2\delta_2(k)b_0k\delta_1(k)b_0k-
32\tau_1^2|\tau|^2\xi_1\xi_2^3b_0^2k^2\delta_2(k)b_0k^2\delta_1(k)b_0
-\\32\tau_1^2|\tau|^2\xi_1\xi_2^3b_0^2k^3\delta_1(k)b_0\delta_2(k)b_0k-
44\tau_1^2|\tau|^2\xi_1\xi_2^3b_0^2k^3\delta_1(k)b_0k\delta_2(k)b_0
-\\32\tau_1^2|\tau|^2\xi_1\xi_2^3b_0^2k^3\delta_2(k)b_0\delta_1(k)b_0k-
44\tau_1^2|\tau|^2\xi_1\xi_2^3b_0^2k^3\delta_2(k)b_0k\delta_1(k)b_0
-\\12\tau_1^2|\tau|^2\xi_1\xi_2^3b_0k\delta_1(k)b_0^2k^2\delta_2(k)b_0k-
12\tau_1^2|\tau|^2\xi_1\xi_2^3b_0k\delta_1(k)b_0^2k^3\delta_2(k)b_0-\\
12\tau_1^2|\tau|^2\xi_1\xi_2^3b_0k\delta_2(k)b_0^2k^2\delta_1(k)b_0k-
12\tau_1^2|\tau|^2\xi_1\xi_2^3b_0k\delta_2(k)b_0^2k^3\delta_1(k)b_0-\\
8\tau_1^2|\tau|^2\xi_2^4b_0^3k^4\delta_1(k)^2b_0
-4\tau_1^2|\tau|^2\xi_2^4b_0^3k^4\delta_1^2(k)b_0k-\\
4\tau_1^2|\tau|^2\xi_2^4b_0^3k^5\delta_1^2(k)b_0-
4\tau_1^2|\tau|^2\xi_2^4b_0^2k^2\delta_1(k)b_0k\delta_1(k)b_0k-\\
8\tau_1^2|\tau|^2\xi_2^4b_0^2k^2\delta_1(k)b_0k^2\delta_1(k)b_0
-8\tau_1^2|\tau|^2\xi_2^4b_0^2k^3\delta_1(k)b_0\delta_1(k)b_0k-\\
12\tau_1^2|\tau|^2\xi_2^4b_0^2k^3\delta_1(k)b_0k\delta_1(k)b_0-
4\tau_1^2|\tau|^2\xi_2^4b_0k\delta_1(k)b_0^2k^2\delta_1(k)b_0k-\\
4\tau_1^2|\tau|^2\xi_2^4b_0k\delta_1(k)b_0^2k^3\delta_1(k)b_0-
8|\tau|^4\xi_1^2\xi_2^2b_0^3k^4\delta_2(k)^2b_0-\\
4|\tau|^4\xi_1^2\xi_2^2b_0^3k^4\delta_2^2(k)b_0k-
4|\tau|^4\xi_1^2\xi_2^2b_0^3k^5\delta_2^2(k)b_0-\\
8|\tau|^4\xi_1^2\xi_2^2b_0^2k^2\delta_2(k)b_0k\delta_2(k)b_0k
-12|\tau|^4\xi_1^2\xi_2^2b_0^2k^2\delta_2(k)b_0k^2\delta_2(k)b_0-\\
12|\tau|^4\xi_1^2\xi_2^2b_0^2k^3\delta_2(k)b_0\delta_2(k)b_0k-
16|\tau|^4\xi_1^2\xi_2^2b_0^2k^3\delta_2(k)b_0k\delta_2(k)b_0-\\
4|\tau|^4\xi_1^2\xi_2^2b_0k\delta_2(k)b_0^2k^2\delta_2(k)b_0k-
4|\tau|^4\xi_1^2\xi_2^2b_0k\delta_2(k)b_0^2k^3\delta_2(k)b_0-\\
8|\tau|^4\xi_1\xi_2^3b_0^3k^4\delta_1(k)\delta_2(k)b_0-
8|\tau|^4\xi_1\xi_2^3b_0^3k^4\delta_2(k)\delta_1(k)b_0-\\
8|\tau|^4\xi_1\xi_2^3b_0^3k^4\delta_1\delta_2(k)b_0k-
8|\tau|^4\xi_1\xi_2^3b_0^3k^5\delta_1\delta_2(k)b_0-\\
4|\tau|^4\xi_1\xi_2^3b_0^2k^2\delta_1(k)b_0k\delta_2(k)b_0k-
8|\tau|^4\xi_1\xi_2^3b_0^2k^2\delta_1(k)b_0k^2\delta_2(k)b_0-\\
4|\tau|^4\xi_1\xi_2^3b_0^2k^2\delta_2(k)b_0k\delta_1(k)b_0k-
8|\tau|^4\xi_1\xi_2^3b_0^2k^2\delta_2(k)b_0k^2\delta_1(k)b_0-\\
8|\tau|^4\xi_1\xi_2^3b_0^2k^3\delta_1(k)b_0\delta_2(k)b_0k
-12|\tau|^4\xi_1\xi_2^3b_0^2k^3\delta_1(k)b_0k\delta_2(k)b_0-\\
8|\tau|^4\xi_1\xi_2^3b_0^2k^3\delta_2(k)b_0\delta_1(k)b_0k-
12|\tau|^4\xi_1\xi_2^3b_0^2k^3\delta_2(k)b_0k\delta_1(k)b_0-\\
4|\tau|^4\xi_1\xi_2^3b_0k\delta_1(k)b_0^2k^2\delta_2(k)b_0k
-4|\tau|^4\xi_1\xi_2^3b_0k\delta_1(k)b_0^2k^3\delta_2(k)b_0-\\
4|\tau|^4\xi_1\xi_2^3b_0k\delta_2(k)b_0^2k^2\delta_1(k)b_0k-
4|\tau|^4\xi_1\xi_2^3b_0k\delta_2(k)b_0^2k^3\delta_1(k)b_0-\\
2|\tau|^4\xi_2^4b_0^2k^2\delta_1(k)b_0k\delta_1(k)b_0k
-2|\tau|^4\xi_2^4b_0^2k^2\delta_1(k)b_0k^2\delta_1(k)b_0-\\
2|\tau|^4\xi_2^4b_0^2k^3\delta_1(k)b_0\delta_1(k)b_0k-
2|\tau|^4\xi_2^4b_0^2k^3\delta_1(k)b_0k\delta_1(k)b_0+\\
8\xi_1^6b_0^3k^4\delta_1(k)b_0k\delta_1(k)b_0k+
8\xi_1^6b_0^3k^4\delta_1(k)b_0k^2\delta_1(k)b_0+\\
8\xi_1^6b_0^3k^5\delta_1(k)b_0\delta_1(k)b_0k+
8\xi_1^6b_0^3k^5\delta_1(k)b_0k\delta_1(k)b_0+\\
4\xi_1^6b_0^2k^2\delta_1(k)b_0^2k^3\delta_1(k)b_0k+
4\xi_1^6b_0^2k^2\delta_1(k)b_0^2k^4\delta_1(k)b_0+\\
4\xi_1^6b_0^2k^3\delta_1(k)b_0^2k^2\delta_1(k)b_0k+
4\xi_1^6b_0^2k^3\delta_1(k)b_0^2k^3\delta_1(k)b_0-\\
32\tau_1|\tau|^4\xi_1\xi_2^3b_0^3k^4\delta_2(k)^2b_0
-16\tau_1|\tau|^4\xi_1\xi_2^3b_0^3k^4\delta_2^2(k)b_0k-\\
16\tau_1|\tau|^4\xi_1\xi_2^3b_0^3k^5\delta_2^2(k)b_0-24
\tau_1|\tau|^4\xi_1\xi_2^3b_0^2k^2\delta_2(k)b_0k\delta_2(k)b_0k-\\
40\tau_1|\tau|^4\xi_1\xi_2^3b_0^2k^2\delta_2(k)b_0k^2\delta_2(k)b_0-
40\tau_1|\tau|^4\xi_1\xi_2^3b_0^2k^3\delta_2(k)b_0\delta_2(k)b_0k-\\
56\tau_1|\tau|^4\xi_1\xi_2^3b_0^2k^3\delta_2(k)b_0k\delta_2(k)b_0
-16\tau_1|\tau|^4\xi_1\xi_2^3b_0k\delta_2(k)b_0^2k^2\delta_2(k)b_0k-\\
16\tau_1|\tau|^4\xi_1\xi_2^3b_0k\delta_2(k)b_0^2k^3\delta_2(k)b_0-
8\tau_1|\tau|^4\xi_2^4b_0^3k^4\delta_1(k)\delta_2(k)b_0-\\
8\tau_1|\tau|^4\xi_2^4b_0^3k^4\delta_2(k)\delta_1(k)b_0-
8\tau_1|\tau|^4\xi_2^4b_0^3k^4\delta_1\delta_2(k)b_0k-\\
8\tau_1|\tau|^4\xi_2^4b_0^3k^5\delta_1\delta_2(k)b_0-
6\tau_1|\tau|^4\xi_2^4b_0^2k^2\delta_1(k)b_0k\delta_2(k)b_0k-\\
10\tau_1|\tau|^4\xi_2^4b_0^2k^2\delta_1(k)b_0k^2\delta_2(k)b_0
-6\tau_1|\tau|^4\xi_2^4b_0^2k^2\delta_2(k)b_0k\delta_1(k)b_0k-\\
10\tau_1|\tau|^4\xi_2^4b_0^2k^2\delta_2(k)b_0k^2\delta_1(k)b_0-
10\tau_1|\tau|^4\xi_2^4b_0^2k^3\delta_1(k)b_0\delta_2(k)b_0k-\\
14\tau_1|\tau|^4\xi_2^4b_0^2k^3\delta_1(k)b_0k\delta_2(k)b_0-
10\tau_1|\tau|^4\xi_2^4b_0^2k^3\delta_2(k)b_0\delta_1(k)b_0k-\\
14\tau_1|\tau|^4\xi_2^4b_0^2k^3\delta_2(k)b_0k\delta_1(k)b_0-
4\tau_1|\tau|^4\xi_2^4b_0k\delta_1(k)b_0^2k^2\delta_2(k)b_0k-\\
4\tau_1|\tau|^4\xi_2^4b_0k\delta_1(k)b_0^2k^3\delta_2(k)b_0-
4\tau_1|\tau|^4\xi_2^4b_0k\delta_2(k)b_0^2k^2\delta_1(k)b_0k-\\
4\tau_1|\tau|^4\xi_2^4b_0k\delta_2(k)b_0^2k^3\delta_1(k)b_0+
8\tau_1\xi_1^6b_0^3k^4\delta_1(k)b_0k\delta_2(k)b_0k+\\
8\tau_1\xi_1^6b_0^3k^4\delta_1(k)b_0k^2\delta_2(k)b_0+
8\tau_1\xi_1^6b_0^3k^4\delta_2(k)b_0k\delta_1(k)b_0k+\\
8\tau_1\xi_1^6b_0^3k^4\delta_2(k)b_0k^2\delta_1(k)b_0+
8\tau_1\xi_1^6b_0^3k^5\delta_1(k)b_0\delta_2(k)b_0k+\\
8\tau_1\xi_1^6b_0^3k^5\delta_1(k)b_0k\delta_2(k)b_0+
8\tau_1\xi_1^6b_0^3k^5\delta_2(k)b_0\delta_1(k)b_0k+\\
8\tau_1\xi_1^6b_0^3k^5\delta_2(k)b_0k\delta_1(k)b_0+
4\tau_1\xi_1^6b_0^2k^2\delta_1(k)b_0^2k^3\delta_2(k)b_0k+\\
4\tau_1\xi_1^6b_0^2k^2\delta_1(k)b_0^2k^4\delta_2(k)b_0+
4\tau_1\xi_1^6b_0^2k^2\delta_2(k)b_0^2k^3\delta_1(k)b_0k+\\
4\tau_1\xi_1^6b_0^2k^2\delta_2(k)b_0^2k^4\delta_1(k)b_0+
4\tau_1\xi_1^6b_0^2k^3\delta_1(k)b_0^2k^2\delta_2(k)b_0k+\\
4\tau_1\xi_1^6b_0^2k^3\delta_1(k)b_0^2k^3\delta_2(k)b_0+
4\tau_1\xi_1^6b_0^2k^3\delta_2(k)b_0^2k^2\delta_1(k)b_0k+\\
4\tau_1\xi_1^6b_0^2k^3\delta_2(k)b_0^2k^3\delta_1(k)b_0+
48\tau_1\xi_1^5\xi_2b_0^3k^4\delta_1(k)b_0k\delta_1(k)b_0k+\\
48\tau_1\xi_1^5\xi_2b_0^3k^4\delta_1(k)b_0k^2\delta_1(k)b_0+
48\tau_1\xi_1^5\xi_2b_0^3k^5\delta_1(k)b_0\delta_1(k)b_0k+\\
48\tau_1\xi_1^5\xi_2b_0^3k^5\delta_1(k)b_0k\delta_1(k)b_0+
24\tau_1\xi_1^5\xi_2b_0^2k^2\delta_1(k)b_0^2k^3\delta_1(k)b_0k+\\
24\tau_1\xi_1^5\xi_2b_0^2k^2\delta_1(k)b_0^2k^4\delta_1(k)b_0+
24\tau_1\xi_1^5\xi_2b_0^2k^3\delta_1(k)b_0^2k^2\delta_1(k)b_0k+\\
24\tau_1\xi_1^5\xi_2b_0^2k^3\delta_1(k)b_0^2k^3\delta_1(k)b_0-
8|\tau|^6\xi_2^4b_0^3k^4\delta_2(k)^2b_0-
4|\tau|^6\xi_2^4b_0^3k^4\delta_2^2(k)b_0k-\\
4|\tau|^6\xi_2^4b_0^3k^5\delta_2^2(k)b_0-
6|\tau|^6\xi_2^4b_0^2k^2\delta_2(k)b_0k\delta_2(k)b_0k-\\
10|\tau|^6\xi_2^4b_0^2k^2\delta_2(k)b_0k^2\delta_2(k)b_0-
10|\tau|^6\xi_2^4b_0^2k^3\delta_2(k)b_0\delta_2(k)b_0k-\\14|\tau|^
6\xi_2^4b_0^2k^3\delta_2(k)b_0k\delta_2(k)b_0-4|\tau|^6\xi_2^4b_0k
\delta_2(k)b_0^2k^2\delta_2(k)b_0k-\\4|\tau|^6\xi_2^4b_0k\delta_2(k)b_0^2k^
3\delta_2(k)b_0+8\tau_1^2\xi_1^6b_0^3k^4\delta_2(k)b_0k\delta_2(k)b_0k+
\\8\tau_1^2\xi_1^6b_0^3k^4\delta_2(k)b_0k^2\delta_2(k)b_0+
8\tau_1^2\xi_1^
6b_0^3k^5\delta_2(k)b_0\delta_2(k)b_0k+\\8\tau_1^2\xi_1^6b_0^3k^5\delta_2(k)b_0k\delta_2(k)b_0+4\tau_1^2\xi_1^6b_0^2k^2\delta_2(k)b_0^2k^
3\delta_2(k)b_0k+\\
4\tau_1^2\xi_1^6b_0^2k^2\delta_2(k)b_0^2k^4\delta_2(k)b_0+4\tau_1^2\xi_1^6b_0^2k^3\delta_2(k)b_0^2k^2\delta_2(k)b_0k+\\4\tau_1^
2\xi_1^6b_0^2k^3\delta_2(k)b_0^2k^3\delta_2(k)b_0+
40\tau_1^2\xi_1^5\xi_2b_0^3k^4\delta_1(k)b_0k\delta_2(k)b_0k+\\40\tau_1^2\xi_1^5\xi_2b_0^3k^
4\delta_1(k)b_0k^2\delta_2(k)b_0+40\tau_1^2\xi_1^5\xi_2b_0^3k^4\delta_2(k)b_0k\delta_1(k)b_0k+\\40\tau_1^2\xi_1^5\xi_2b_0^3k^4\delta_2(k)b_0k^
2\delta_1(k)b_0+40\tau_1^2\xi_1^5\xi_2b_0^3k^5\delta_1(k)b_0\delta_2(k)b_0k+\\40\tau_1^2\xi_1^5\xi_2b_0^3k^5\delta_1(k)b_0k\delta_2(k)b_0+40\tau_1^
2\xi_1^5\xi_2b_0^3k^5\delta_2(k)b_0\delta_1(k)b_0k+\\40\tau_1^2\xi_1^5\xi_2b_0^3k^5\delta_2(k)b_0k\delta_1(k)b_0+20\tau_1^2\xi_1^5\xi_2b_0^2k^
2\delta_1(k)b_0^2k^3\delta_2(k)b_0k+\\20\tau_1^2\xi_1^5\xi_2b_0^2k^2\delta_1(k)b_0^2k^4\delta_2(k)b_0+20\tau_1^2\xi_1^5\xi_2b_0^2k^2\delta_2(k)b_0^
2k^3\delta_1(k)b_0k+\\20\tau_1^2\xi_1^5\xi_2b_0^2k^2\delta_2(k)b_0^2k^4\delta_1(k)b_0+20\tau_1^2\xi_1^5\xi_2b_0^2k^3\delta_1(k)b_0^2k^
2\delta_2(k)b_0k+\\20\tau_1^2\xi_1^5\xi_2b_0^2k^3\delta_1(k)b_0^2k^3\delta_2(k)b_0+20\tau_1^2\xi_1^5\xi_2b_0^2k^3\delta_2(k)b_0^2k^2\delta_1(k)b_0k+\\
20\tau_1^2\xi_1^5\xi_2b_0^2k^3\delta_2(k)b_0^2k^3\delta_1(k)b_0+104\tau_1^2\xi_1^4\xi_2^2b_0^3k^4\delta_1(k)b_0k\delta_1(k)b_0k+\\104\tau_1^2\xi_1^
4\xi_2^2b_0^3k^4\delta_1(k)b_0k^2\delta_1(k)b_0+104\tau_1^2\xi_1^4\xi_2^2b_0^3k^5\delta_1(k)b_0\delta_1(k)b_0k+\\104\tau_1^2\xi_1^4\xi_2^2b_0^3k^
5\delta_1(k)b_0k\delta_1(k)b_0+52\tau_1^2\xi_1^4\xi_2^2b_0^2k^2\delta_1(k)b_0^2k^3\delta_1(k)b_0k+\\52\tau_1^2\xi_1^4\xi_2^2b_0^2k^2\delta_1(k)b_0^
2k^4\delta_1(k)b_0+52\tau_1^2\xi_1^4\xi_2^2b_0^2k^3\delta_1(k)b_0^2k^2\delta_1(k)b_0k+\\52\tau_1^2\xi_1^4\xi_2^2b_0^2k^3\delta_1(k)b_0^2k^
3\delta_1(k)b_0+8|\tau|^2\xi_1^5\xi_2b_0^3k^4\delta_1(k)b_0k\delta_2(k)b_0k+\\8|\tau|^2\xi_1^5\xi_2b_0^3k^4\delta_1(k)b_0k^2\delta_2(k)b_0+8|\tau|^2\xi_1^
5\xi_2b_0^3k^4\delta_2(k)b_0k\delta_1(k)b_0k+\\8|\tau|^2\xi_1^5\xi_2b_0^3k^4\delta_2(k)b_0k^2\delta_1(k)b_0+8|\tau|^2\xi_1^5\xi_2b_0^3k^
5\delta_1(k)b_0\delta_2(k)b_0k+\\8|\tau|^2\xi_1^5\xi_2b_0^3k^5\delta_1(k)b_0k\delta_2(k)b_0+8|\tau|^2\xi_1^5\xi_2b_0^3k^5\delta_2(k)b_0\delta_1(k)b_0k+\\8|\tau|^
2\xi_1^5\xi_2b_0^3k^5\delta_2(k)b_0k\delta_1(k)b_0+4|\tau|^2\xi_1^5\xi_2b_0^2k^2\delta_1(k)b_0^2k^3\delta_2(k)b_0k+\\4|\tau|^2\xi_1^5\xi_2b_0^2k^
2\delta_1(k)b_0^2k^4\delta_2(k)b_0+4|\tau|^2\xi_1^5\xi_2b_0^2k^2\delta_2(k)b_0^2k^3\delta_1(k)b_0k+\\4|\tau|^2\xi_1^5\xi_2b_0^2k^2\delta_2(k)b_0^2k^
4\delta_1(k)b_0+4|\tau|^2\xi_1^5\xi_2b_0^2k^3\delta_1(k)b_0^2k^2\delta_2(k)b_0k+\\4|\tau|^2\xi_1^5\xi_2b_0^2k^3\delta_1(k)b_0^2k^3\delta_2(k)b_0+4|\tau|^2\xi_1^
5\xi_2b_0^2k^3\delta_2(k)b_0^2k^2\delta_1(k)b_0k+\\4|\tau|^2\xi_1^5\xi_2b_0^2k^3\delta_2(k)b_0^2k^3\delta_1(k)b_0+16|\tau|^2\xi_1^4\xi_2^2b_0^3k^
4\delta_1(k)b_0k\delta_1(k)b_0k+\\16|\tau|^2\xi_1^4\xi_2^2b_0^3k^4\delta_1(k)b_0k^2\delta_1(k)b_0+16|\tau|^2\xi_1^4\xi_2^2b_0^3k^
5\delta_1(k)b_0\delta_1(k)b_0k+\\16|\tau|^2\xi_1^4\xi_2^2b_0^3k^5\delta_1(k)b_0k\delta_1(k)b_0+8|\tau|^2\xi_1^4\xi_2^2b_0^2k^2\delta_1(k)b_0^2k^
3\delta_1(k)b_0k+\\8|\tau|^2\xi_1^4\xi_2^2b_0^2k^2\delta_1(k)b_0^2k^4\delta_1(k)b_0+8|\tau|^2\xi_1^4\xi_2^2b_0^2k^3\delta_1(k)b_0^2k^2\delta_1(k)b_0k+\\8|\tau|^
2\xi_1^4\xi_2^2b_0^2k^3\delta_1(k)b_0^2k^3\delta_1(k)b_0+32\tau_1^3\xi_1^5\xi_2b_0^3k^4\delta_2(k)b_0k\delta_2(k)b_0k+\\32\tau_1^3\xi_1^5\xi_2b_0^
3k^4\delta_2(k)b_0k^2\delta_2(k)b_0+32\tau_1^3\xi_1^5\xi_2b_0^3k^5\delta_2(k)b_0\delta_2(k)b_0k+\\32\tau_1^3\xi_1^5\xi_2b_0^3k^
5\delta_2(k)b_0k\delta_2(k)b_0+16\tau_1^3\xi_1^5\xi_2b_0^2k^2\delta_2(k)b_0^2k^3\delta_2(k)b_0k+\\16\tau_1^3\xi_1^5\xi_2b_0^2k^2\delta_2(k)b_0^2k^
4\delta_2(k)b_0+16\tau_1^3\xi_1^5\xi_2b_0^2k^3\delta_2(k)b_0^2k^2\delta_2(k)b_0k+\\16\tau_1^3\xi_1^5\xi_2b_0^2k^3\delta_2(k)b_0^2k^3\delta_2(k)b_0+
64\tau_1^3\xi_1^4\xi_2^2b_0^3k^4\delta_1(k)b_0k\delta_2(k)b_0k+\\64\tau_1^3\xi_1^4\xi_2^2b_0^3k^4\delta_1(k)b_0k^2\delta_2(k)b_0+64\tau_1^3\xi_1^
4\xi_2^2b_0^3k^4\delta_2(k)b_0k\delta_1(k)b_0k+\\64\tau_1^3\xi_1^4\xi_2^2b_0^3k^4\delta_2(k)b_0k^2\delta_1(k)b_0+64\tau_1^3\xi_1^4\xi_2^2b_0^3k^
5\delta_1(k)b_0\delta_2(k)b_0k+\\64\tau_1^3\xi_1^4\xi_2^2b_0^3k^5\delta_1(k)b_0k\delta_2(k)b_0+64\tau_1^3\xi_1^4\xi_2^2b_0^3k^
5\delta_2(k)b_0\delta_1(k)b_0k+\\64\tau_1^3\xi_1^4\xi_2^2b_0^3k^5\delta_2(k)b_0k\delta_1(k)b_0+32\tau_1^3\xi_1^4\xi_2^2b_0^2k^2\delta_1(k)b_0^2k^
3\delta_2(k)b_0k+\\32\tau_1^3\xi_1^4\xi_2^2b_0^2k^2\delta_1(k)b_0^2k^4\delta_2(k)b_0+32\tau_1^3\xi_1^4\xi_2^2b_0^2k^2\delta_2(k)b_0^2k^
3\delta_1(k)b_0k+\\32\tau_1^3\xi_1^4\xi_2^2b_0^2k^2\delta_2(k)b_0^2k^4\delta_1(k)b_0+32\tau_1^3\xi_1^4\xi_2^2b_0^2k^3\delta_1(k)b_0^2k^
2\delta_2(k)b_0k+\\32\tau_1^3\xi_1^4\xi_2^2b_0^2k^3\delta_1(k)b_0^2k^3\delta_2(k)b_0+32\tau_1^3\xi_1^4\xi_2^2b_0^2k^3\delta_2(k)b_0^2k^
2\delta_1(k)b_0k+\\32\tau_1^3\xi_1^4\xi_2^2b_0^2k^3\delta_2(k)b_0^2k^3\delta_1(k)b_0+96\tau_1^3\xi_1^3\xi_2^3b_0^3k^4\delta_1(k)b_0k\delta_1(k)b_0k+\\
96\tau_1^3\xi_1^3\xi_2^3b_0^3k^4\delta_1(k)b_0k^2\delta_1(k)b_0+96\tau_1^3\xi_1^3\xi_2^3b_0^3k^5\delta_1(k)b_0\delta_1(k)b_0k+\\96\tau_1^3\xi_1^
3\xi_2^3b_0^3k^5\delta_1(k)b_0k\delta_1(k)b_0+48\tau_1^3\xi_1^3\xi_2^3b_0^2k^2\delta_1(k)b_0^2k^3\delta_1(k)b_0k+\\48\tau_1^3\xi_1^3\xi_2^3b_0^2k^
2\delta_1(k)b_0^2k^4\delta_1(k)b_0+48\tau_1^3\xi_1^3\xi_2^3b_0^2k^3\delta_1(k)b_0^2k^2\delta_1(k)b_0k+\\48\tau_1^3\xi_1^3\xi_2^3b_0^2k^
3\delta_1(k)b_0^2k^3\delta_1(k)b_0+16\tau_1|\tau|^2\xi_1^5\xi_2b_0^3k^4\delta_2(k)b_0k\delta_2(k)b_0k+\\16\tau_1|\tau|^2\xi_1^5\xi_2b_0^3k^4\delta_2(k)b_0k^
2\delta_2(k)b_0+16\tau_1|\tau|^2\xi_1^5\xi_2b_0^3k^5\delta_2(k)b_0\delta_2(k)b_0k+\\16\tau_1|\tau|^2\xi_1^5\xi_2b_0^3k^5\delta_2(k)b_0k\delta_2(k)b_0+
8\tau_1|\tau|^2\xi_1^5\xi_2b_0^2k^2\delta_2(k)b_0^2k^3\delta_2(k)b_0k+\\8\tau_1|\tau|^2\xi_1^5\xi_2b_0^2k^2\delta_2(k)b_0^2k^4\delta_2(k)b_0+8\tau_1|\tau|^
2\xi_1^5\xi_2b_0^2k^3\delta_2(k)b_0^2k^2\delta_2(k)b_0k+\\8\tau_1|\tau|^2\xi_1^5\xi_2b_0^2k^3\delta_2(k)b_0^2k^3\delta_2(k)b_0+56\tau_1|\tau|^2\xi_1^4\xi_2^
2b_0^3k^4\delta_1(k)b_0k\delta_2(k)b_0k+\\56\tau_1|\tau|^2\xi_1^4\xi_2^2b_0^3k^4\delta_1(k)b_0k^2\delta_2(k)b_0+56\tau_1|\tau|^2\xi_1^4\xi_2^2b_0^3k^
4\delta_2(k)b_0k\delta_1(k)b_0k+\\56\tau_1|\tau|^2\xi_1^4\xi_2^2b_0^3k^4\delta_2(k)b_0k^2\delta_1(k)b_0+56\tau_1|\tau|^2\xi_1^4\xi_2^2b_0^3k^
5\delta_1(k)b_0\delta_2(k)b_0k+\\56\tau_1|\tau|^2\xi_1^4\xi_2^2b_0^3k^5\delta_1(k)b_0k\delta_2(k)b_0+56\tau_1|\tau|^2\xi_1^4\xi_2^2b_0^3k^
5\delta_2(k)b_0\delta_1(k)b_0k+\\56\tau_1|\tau|^2\xi_1^4\xi_2^2b_0^3k^5\delta_2(k)b_0k\delta_1(k)b_0+28\tau_1|\tau|^2\xi_1^4\xi_2^2b_0^2k^2\delta_1(k)b_0^2k^
3\delta_2(k)b_0k+\\28\tau_1|\tau|^2\xi_1^4\xi_2^2b_0^2k^2\delta_1(k)b_0^2k^4\delta_2(k)b_0+28\tau_1|\tau|^2\xi_1^4\xi_2^2b_0^2k^2\delta_2(k)b_0^2k^
3\delta_1(k)b_0k+\\28\tau_1|\tau|^2\xi_1^4\xi_2^2b_0^2k^2\delta_2(k)b_0^2k^4\delta_1(k)b_0+28\tau_1|\tau|^2\xi_1^4\xi_2^2b_0^2k^3\delta_1(k)b_0^2k^
2\delta_2(k)b_0k+\\28\tau_1|\tau|^2\xi_1^4\xi_2^2b_0^2k^3\delta_1(k)b_0^2k^3\delta_2(k)b_0+28\tau_1|\tau|^2\xi_1^4\xi_2^2b_0^2k^3\delta_2(k)b_0^2k^
2\delta_1(k)b_0k+\\28\tau_1|\tau|^2\xi_1^4\xi_2^2b_0^2k^3\delta_2(k)b_0^2k^3\delta_1(k)b_0+64\tau_1|\tau|^2\xi_1^3\xi_2^3b_0^3k^
4\delta_1(k)b_0k\delta_1(k)b_0k+\\64\tau_1|\tau|^2\xi_1^3\xi_2^3b_0^3k^4\delta_1(k)b_0k^2\delta_1(k)b_0+64\tau_1|\tau|^2\xi_1^3\xi_2^3b_0^3k^
5\delta_1(k)b_0\delta_1(k)b_0k+\\64\tau_1|\tau|^2\xi_1^3\xi_2^3b_0^3k^5\delta_1(k)b_0k\delta_1(k)b_0+32\tau_1|\tau|^2\xi_1^3\xi_2^3b_0^2k^2\delta_1(k)b_0^2k^
3\delta_1(k)b_0k+\\32\tau_1|\tau|^2\xi_1^3\xi_2^3b_0^2k^2\delta_1(k)b_0^2k^4\delta_1(k)b_0+32\tau_1|\tau|^2\xi_1^3\xi_2^3b_0^2k^3\delta_1(k)b_0^2k^
2\delta_1(k)b_0k+\\32\tau_1|\tau|^2\xi_1^3\xi_2^3b_0^2k^3\delta_1(k)b_0^2k^3\delta_1(k)b_0+32\tau_1^4\xi_1^4\xi_2^2b_0^3k^
4\delta_2(k)b_0k\delta_2(k)b_0k+\\32\tau_1^4\xi_1^4\xi_2^2b_0^3k^4\delta_2(k)b_0k^2\delta_2(k)b_0+32\tau_1^4\xi_1^4\xi_2^2b_0^3k^
5\delta_2(k)b_0\delta_2(k)b_0k+\\32\tau_1^4\xi_1^4\xi_2^2b_0^3k^5\delta_2(k)b_0k\delta_2(k)b_0+16\tau_1^4\xi_1^4\xi_2^2b_0^2k^2\delta_2(k)b_0^2k^
3\delta_2(k)b_0k+\\16\tau_1^4\xi_1^4\xi_2^2b_0^2k^2\delta_2(k)b_0^2k^4\delta_2(k)b_0+16\tau_1^4\xi_1^4\xi_2^2b_0^2k^3\delta_2(k)b_0^2k^
2\delta_2(k)b_0k+\\16\tau_1^4\xi_1^4\xi_2^2b_0^2k^3\delta_2(k)b_0^2k^3\delta_2(k)b_0+32\tau_1^4\xi_1^3\xi_2^3b_0^3k^4\delta_1(k)b_0k\delta_2(k)b_0k+\\
32\tau_1^4\xi_1^3\xi_2^3b_0^3k^4\delta_1(k)b_0k^2\delta_2(k)b_0+32\tau_1^4\xi_1^3\xi_2^3b_0^3k^4\delta_2(k)b_0k\delta_1(k)b_0k+\\32\tau_1^4\xi_1^
3\xi_2^3b_0^3k^4\delta_2(k)b_0k^2\delta_1(k)b_0+32\tau_1^4\xi_1^3\xi_2^3b_0^3k^5\delta_1(k)b_0\delta_2(k)b_0k+\\32\tau_1^4\xi_1^3\xi_2^3b_0^3k^
5\delta_1(k)b_0k\delta_2(k)b_0+32\tau_1^4\xi_1^3\xi_2^3b_0^3k^5\delta_2(k)b_0\delta_1(k)b_0k+\\32\tau_1^4\xi_1^3\xi_2^3b_0^3k^
5\delta_2(k)b_0k\delta_1(k)b_0+16\tau_1^4\xi_1^3\xi_2^3b_0^2k^2\delta_1(k)b_0^2k^3\delta_2(k)b_0k+\\16\tau_1^4\xi_1^3\xi_2^3b_0^2k^2\delta_1(k)b_0^
2k^4\delta_2(k)b_0+16\tau_1^4\xi_1^3\xi_2^3b_0^2k^2\delta_2(k)b_0^2k^3\delta_1(k)b_0k+\\16\tau_1^4\xi_1^3\xi_2^3b_0^2k^2\delta_2(k)b_0^2k^
4\delta_1(k)b_0+16\tau_1^4\xi_1^3\xi_2^3b_0^2k^3\delta_1(k)b_0^2k^2\delta_2(k)b_0k+\\16\tau_1^4\xi_1^3\xi_2^3b_0^2k^3\delta_1(k)b_0^2k^
3\delta_2(k)b_0+16\tau_1^4\xi_1^3\xi_2^3b_0^2k^3\delta_2(k)b_0^2k^2\delta_1(k)b_0k+\\16\tau_1^4\xi_1^3\xi_2^3b_0^2k^3\delta_2(k)b_0^2k^
3\delta_1(k)b_0+32\tau_1^4\xi_1^2\xi_2^4b_0^3k^4\delta_1(k)b_0k\delta_1(k)b_0k+\\32\tau_1^4\xi_1^2\xi_2^4b_0^3k^4\delta_1(k)b_0k^2\delta_1(k)b_0+
32\tau_1^4\xi_1^2\xi_2^4b_0^3k^5\delta_1(k)b_0\delta_1(k)b_0k+\\32\tau_1^4\xi_1^2\xi_2^4b_0^3k^5\delta_1(k)b_0k\delta_1(k)b_0+16\tau_1^4\xi_1^2\xi_2^
4b_0^2k^2\delta_1(k)b_0^2k^3\delta_1(k)b_0k+\\16\tau_1^4\xi_1^2\xi_2^4b_0^2k^2\delta_1(k)b_0^2k^4\delta_1(k)b_0+16\tau_1^4\xi_1^2\xi_2^4b_0^2k^
3\delta_1(k)b_0^2k^2\delta_1(k)b_0k+\\16\tau_1^4\xi_1^2\xi_2^4b_0^2k^3\delta_1(k)b_0^2k^3\delta_1(k)b_0+80\tau_1^2|\tau|^2\xi_1^4\xi_2^2b_0^3k^
4\delta_2(k)b_0k\delta_2(k)b_0k+\\80\tau_1^2|\tau|^2\xi_1^4\xi_2^2b_0^3k^4\delta_2(k)b_0k^2\delta_2(k)b_0+80\tau_1^2|\tau|^2\xi_1^4\xi_2^2b_0^3k^
5\delta_2(k)b_0\delta_2(k)b_0k+\\80\tau_1^2|\tau|^2\xi_1^4\xi_2^2b_0^3k^5\delta_2(k)b_0k\delta_2(k)b_0+40\tau_1^2|\tau|^2\xi_1^4\xi_2^2b_0^2k^2\delta_2(k)b_0^
2k^3\delta_2(k)b_0k+\\40\tau_1^2|\tau|^2\xi_1^4\xi_2^2b_0^2k^2\delta_2(k)b_0^2k^4\delta_2(k)b_0+40\tau_1^2|\tau|^2\xi_1^4\xi_2^2b_0^2k^3\delta_2(k)b_0^2k^
2\delta_2(k)b_0k+\\40\tau_1^2|\tau|^2\xi_1^4\xi_2^2b_0^2k^3\delta_2(k)b_0^2k^3\delta_2(k)b_0+112\tau_1^2|\tau|^2\xi_1^3\xi_2^3b_0^3k^
4\delta_1(k)b_0k\delta_2(k)b_0k+\\112\tau_1^2|\tau|^2\xi_1^3\xi_2^3b_0^3k^4\delta_1(k)b_0k^2\delta_2(k)b_0+112\tau_1^2|\tau|^2\xi_1^3\xi_2^3b_0^3k^
4\delta_2(k)b_0k\delta_1(k)b_0k+\\112\tau_1^2|\tau|^2\xi_1^3\xi_2^3b_0^3k^4\delta_2(k)b_0k^2\delta_1(k)b_0+112\tau_1^2|\tau|^2\xi_1^3\xi_2^3b_0^3k^
5\delta_1(k)b_0\delta_2(k)b_0k+\\112\tau_1^2|\tau|^2\xi_1^3\xi_2^3b_0^3k^5\delta_1(k)b_0k\delta_2(k)b_0+112\tau_1^2|\tau|^2\xi_1^3\xi_2^3b_0^3k^
5\delta_2(k)b_0\delta_1(k)b_0k+\\112\tau_1^2|\tau|^2\xi_1^3\xi_2^3b_0^3k^5\delta_2(k)b_0k\delta_1(k)b_0+56\tau_1^2|\tau|^2\xi_1^3\xi_2^3b_0^2k^2\delta_1(k)b_0^
2k^3\delta_2(k)b_0k+\\56\tau_1^2|\tau|^2\xi_1^3\xi_2^3b_0^2k^2\delta_1(k)b_0^2k^4\delta_2(k)b_0+56\tau_1^2|\tau|^2\xi_1^3\xi_2^3b_0^2k^2\delta_2(k)b_0^2k^
3\delta_1(k)b_0k+\\56\tau_1^2|\tau|^2\xi_1^3\xi_2^3b_0^2k^2\delta_2(k)b_0^2k^4\delta_1(k)b_0+56\tau_1^2|\tau|^2\xi_1^3\xi_2^3b_0^2k^3\delta_1(k)b_0^2k^
2\delta_2(k)b_0k+\\56\tau_1^2|\tau|^2\xi_1^3\xi_2^3b_0^2k^3\delta_1(k)b_0^2k^3\delta_2(k)b_0+56\tau_1^2|\tau|^2\xi_1^3\xi_2^3b_0^2k^3\delta_2(k)b_0^2k^
2\delta_1(k)b_0k+\\56\tau_1^2|\tau|^2\xi_1^3\xi_2^3b_0^2k^3\delta_2(k)b_0^2k^3\delta_1(k)b_0+80\tau_1^2|\tau|^2\xi_1^2\xi_2^4b_0^3k^
4\delta_1(k)b_0k\delta_1(k)b_0k+\\80\tau_1^2|\tau|^2\xi_1^2\xi_2^4b_0^3k^4\delta_1(k)b_0k^2\delta_1(k)b_0+80\tau_1^2|\tau|^2\xi_1^2\xi_2^4b_0^3k^
5\delta_1(k)b_0\delta_1(k)b_0k+\\80\tau_1^2|\tau|^2\xi_1^2\xi_2^4b_0^3k^5\delta_1(k)b_0k\delta_1(k)b_0+40\tau_1^2|\tau|^2\xi_1^2\xi_2^4b_0^2k^2\delta_1(k)b_0^
2k^3\delta_1(k)b_0k+\\40\tau_1^2|\tau|^2\xi_1^2\xi_2^4b_0^2k^2\delta_1(k)b_0^2k^4\delta_1(k)b_0+40\tau_1^2|\tau|^2\xi_1^2\xi_2^4b_0^2k^3\delta_1(k)b_0^2k^
2\delta_1(k)b_0k+\\40\tau_1^2|\tau|^2\xi_1^2\xi_2^4b_0^2k^3\delta_1(k)b_0^2k^3\delta_1(k)b_0+8|\tau|^4\xi_1^4\xi_2^2b_0^3k^4\delta_2(k)b_0k\delta_2(k)b_0k+\\
8|\tau|^4\xi_1^4\xi_2^2b_0^3k^4\delta_2(k)b_0k^2\delta_2(k)b_0+8|\tau|^4\xi_1^4\xi_2^2b_0^3k^5\delta_2(k)b_0\delta_2(k)b_0k+\\8|\tau|^4\xi_1^4\xi_2^2b_0^3k^
5\delta_2(k)b_0k\delta_2(k)b_0+4|\tau|^4\xi_1^4\xi_2^2b_0^2k^2\delta_2(k)b_0^2k^3\delta_2(k)b_0k+\\4|\tau|^4\xi_1^4\xi_2^2b_0^2k^2\delta_2(k)b_0^2k^
4\delta_2(k)b_0+4|\tau|^4\xi_1^4\xi_2^2b_0^2k^3\delta_2(k)b_0^2k^2\delta_2(k)b_0k+\\4|\tau|^4\xi_1^4\xi_2^2b_0^2k^3\delta_2(k)b_0^2k^3\delta_2(k)b_0+16|\tau|^
4\xi_1^3\xi_2^3b_0^3k^4\delta_1(k)b_0k\delta_2(k)b_0k+\\16|\tau|^4\xi_1^3\xi_2^3b_0^3k^4\delta_1(k)b_0k^2\delta_2(k)b_0+16|\tau|^4\xi_1^3\xi_2^3b_0^3k^
4\delta_2(k)b_0k\delta_1(k)b_0k+\\16|\tau|^4\xi_1^3\xi_2^3b_0^3k^4\delta_2(k)b_0k^2\delta_1(k)b_0+16|\tau|^4\xi_1^3\xi_2^3b_0^3k^
5\delta_1(k)b_0\delta_2(k)b_0k+\\16|\tau|^4\xi_1^3\xi_2^3b_0^3k^5\delta_1(k)b_0k\delta_2(k)b_0+16|\tau|^4\xi_1^3\xi_2^3b_0^3k^5\delta_2(k)b_0\delta_1(k)b_0k+\\
16|\tau|^4\xi_1^3\xi_2^3b_0^3k^5\delta_2(k)b_0k\delta_1(k)b_0+8|\tau|^4\xi_1^3\xi_2^3b_0^2k^2\delta_1(k)b_0^2k^3\delta_2(k)b_0k+\\8|\tau|^4\xi_1^3\xi_2^3b_0^
2k^2\delta_1(k)b_0^2k^4\delta_2(k)b_0+8|\tau|^4\xi_1^3\xi_2^3b_0^2k^2\delta_2(k)b_0^2k^3\delta_1(k)b_0k+\\8|\tau|^4\xi_1^3\xi_2^3b_0^2k^2\delta_2(k)b_0^2k^
4\delta_1(k)b_0+8|\tau|^4\xi_1^3\xi_2^3b_0^2k^3\delta_1(k)b_0^2k^2\delta_2(k)b_0k+\\8|\tau|^4\xi_1^3\xi_2^3b_0^2k^3\delta_1(k)b_0^2k^3\delta_2(k)b_0+8|\tau|^
4\xi_1^3\xi_2^3b_0^2k^3\delta_2(k)b_0^2k^2\delta_1(k)b_0k+\\8|\tau|^4\xi_1^3\xi_2^3b_0^2k^3\delta_2(k)b_0^2k^3\delta_1(k)b_0+8|\tau|^4\xi_1^2\xi_2^4b_0^3k^
4\delta_1(k)b_0k\delta_1(k)b_0k+\\8|\tau|^4\xi_1^2\xi_2^4b_0^3k^4\delta_1(k)b_0k^2\delta_1(k)b_0+8|\tau|^4\xi_1^2\xi_2^4b_0^3k^
5\delta_1(k)b_0\delta_1(k)b_0k+\\8|\tau|^4\xi_1^2\xi_2^4b_0^3k^5\delta_1(k)b_0k\delta_1(k)b_0+4|\tau|^4\xi_1^2\xi_2^4b_0^2k^2\delta_1(k)b_0^2k^
3\delta_1(k)b_0k+\\4|\tau|^4\xi_1^2\xi_2^4b_0^2k^2\delta_1(k)b_0^2k^4\delta_1(k)b_0+4|\tau|^4\xi_1^2\xi_2^4b_0^2k^3\delta_1(k)b_0^2k^2\delta_1(k)b_0k+\\4|\tau|^
4\xi_1^2\xi_2^4b_0^2k^3\delta_1(k)b_0^2k^3\delta_1(k)b_0+96\tau_1^3|\tau|^2\xi_1^3\xi_2^3b_0^3k^4\delta_2(k)b_0k\delta_2(k)b_0k+\\96\tau_1^3|\tau|^2\xi_1^
3\xi_2^3b_0^3k^4\delta_2(k)b_0k^2\delta_2(k)b_0+96\tau_1^3|\tau|^2\xi_1^3\xi_2^3b_0^3k^5\delta_2(k)b_0\delta_2(k)b_0k+\\96\tau_1^3|\tau|^2\xi_1^3\xi_2^3b_0^
3k^5\delta_2(k)b_0k\delta_2(k)b_0+48\tau_1^3|\tau|^2\xi_1^3\xi_2^3b_0^2k^2\delta_2(k)b_0^2k^3\delta_2(k)b_0k+\\48\tau_1^3|\tau|^2\xi_1^3\xi_2^3b_0^2k^
2\delta_2(k)b_0^2k^4\delta_2(k)b_0+48\tau_1^3|\tau|^2\xi_1^3\xi_2^3b_0^2k^3\delta_2(k)b_0^2k^2\delta_2(k)b_0k+\\48\tau_1^3|\tau|^2\xi_1^3\xi_2^3b_0^2k^
3\delta_2(k)b_0^2k^3\delta_2(k)b_0+64\tau_1^3|\tau|^2\xi_1^2\xi_2^4b_0^3k^4\delta_1(k)b_0k\delta_2(k)b_0k+\\64\tau_1^3|\tau|^2\xi_1^2\xi_2^4b_0^3k^
4\delta_1(k)b_0k^2\delta_2(k)b_0+64\tau_1^3|\tau|^2\xi_1^2\xi_2^4b_0^3k^4\delta_2(k)b_0k\delta_1(k)b_0k+\\64\tau_1^3|\tau|^2\xi_1^2\xi_2^4b_0^3k^
4\delta_2(k)b_0k^2\delta_1(k)b_0+64\tau_1^3|\tau|^2\xi_1^2\xi_2^4b_0^3k^5\delta_1(k)b_0\delta_2(k)b_0k+\\64\tau_1^3|\tau|^2\xi_1^2\xi_2^4b_0^3k^
5\delta_1(k)b_0k\delta_2(k)b_0+64\tau_1^3|\tau|^2\xi_1^2\xi_2^4b_0^3k^5\delta_2(k)b_0\delta_1(k)b_0k+\\64\tau_1^3|\tau|^2\xi_1^2\xi_2^4b_0^3k^
5\delta_2(k)b_0k\delta_1(k)b_0+32\tau_1^3|\tau|^2\xi_1^2\xi_2^4b_0^2k^2\delta_1(k)b_0^2k^3\delta_2(k)b_0k+\\32\tau_1^3|\tau|^2\xi_1^2\xi_2^4b_0^2k^
2\delta_1(k)b_0^2k^4\delta_2(k)b_0+32\tau_1^3|\tau|^2\xi_1^2\xi_2^4b_0^2k^2\delta_2(k)b_0^2k^3\delta_1(k)b_0k+\\32\tau_1^3|\tau|^2\xi_1^2\xi_2^4b_0^2k^
2\delta_2(k)b_0^2k^4\delta_1(k)b_0+32\tau_1^3|\tau|^2\xi_1^2\xi_2^4b_0^2k^3\delta_1(k)b_0^2k^2\delta_2(k)b_0k+\\32\tau_1^3|\tau|^2\xi_1^2\xi_2^4b_0^2k^
3\delta_1(k)b_0^2k^3\delta_2(k)b_0+32\tau_1^3|\tau|^2\xi_1^2\xi_2^4b_0^2k^3\delta_2(k)b_0^2k^2\delta_1(k)b_0k+\\32\tau_1^3|\tau|^2\xi_1^2\xi_2^4b_0^2k^
3\delta_2(k)b_0^2k^3\delta_1(k)b_0+32\tau_1^3|\tau|^2\xi_1\xi_2^5b_0^3k^4\delta_1(k)b_0k\delta_1(k)b_0k+\\32\tau_1^3|\tau|^2\xi_1\xi_2^5b_0^3k^
4\delta_1(k)b_0k^2\delta_1(k)b_0+32\tau_1^3|\tau|^2\xi_1\xi_2^5b_0^3k^5\delta_1(k)b_0\delta_1(k)b_0k+\\32\tau_1^3|\tau|^2\xi_1\xi_2^5b_0^3k^
5\delta_1(k)b_0k\delta_1(k)b_0+16\tau_1^3|\tau|^2\xi_1\xi_2^5b_0^2k^2\delta_1(k)b_0^2k^3\delta_1(k)b_0k+\\16\tau_1^3|\tau|^2\xi_1\xi_2^5b_0^2k^
2\delta_1(k)b_0^2k^4\delta_1(k)b_0+16\tau_1^3|\tau|^2\xi_1\xi_2^5b_0^2k^3\delta_1(k)b_0^2k^2\delta_1(k)b_0k+\\16\tau_1^3|\tau|^2\xi_1\xi_2^5b_0^2k^
3\delta_1(k)b_0^2k^3\delta_1(k)b_0+64\tau_1|\tau|^4\xi_1^3\xi_2^3b_0^3k^4\delta_2(k)b_0k\delta_2(k)b_0k+\\64\tau_1|\tau|^4\xi_1^3\xi_2^3b_0^3k^
4\delta_2(k)b_0k^2\delta_2(k)b_0+64\tau_1|\tau|^4\xi_1^3\xi_2^3b_0^3k^5\delta_2(k)b_0\delta_2(k)b_0k+\\64\tau_1|\tau|^4\xi_1^3\xi_2^3b_0^3k^
5\delta_2(k)b_0k\delta_2(k)b_0+32\tau_1|\tau|^4\xi_1^3\xi_2^3b_0^2k^2\delta_2(k)b_0^2k^3\delta_2(k)b_0k+\\32\tau_1|\tau|^4\xi_1^3\xi_2^3b_0^2k^
2\delta_2(k)b_0^2k^4\delta_2(k)b_0+32\tau_1|\tau|^4\xi_1^3\xi_2^3b_0^2k^3\delta_2(k)b_0^2k^2\delta_2(k)b_0k+\\32\tau_1|\tau|^4\xi_1^3\xi_2^3b_0^2k^
3\delta_2(k)b_0^2k^3\delta_2(k)b_0+56\tau_1|\tau|^4\xi_1^2\xi_2^4b_0^3k^4\delta_1(k)b_0k\delta_2(k)b_0k+\\56\tau_1|\tau|^4\xi_1^2\xi_2^4b_0^3k^
4\delta_1(k)b_0k^2\delta_2(k)b_0+56\tau_1|\tau|^4\xi_1^2\xi_2^4b_0^3k^4\delta_2(k)b_0k\delta_1(k)b_0k+\\56\tau_1|\tau|^4\xi_1^2\xi_2^4b_0^3k^
4\delta_2(k)b_0k^2\delta_1(k)b_0+56\tau_1|\tau|^4\xi_1^2\xi_2^4b_0^3k^5\delta_1(k)b_0\delta_2(k)b_0k+\\56\tau_1|\tau|^4\xi_1^2\xi_2^4b_0^3k^
5\delta_1(k)b_0k\delta_2(k)b_0+56\tau_1|\tau|^4\xi_1^2\xi_2^4b_0^3k^5\delta_2(k)b_0\delta_1(k)b_0k+\\56\tau_1|\tau|^4\xi_1^2\xi_2^4b_0^3k^
5\delta_2(k)b_0k\delta_1(k)b_0+28\tau_1|\tau|^4\xi_1^2\xi_2^4b_0^2k^2\delta_1(k)b_0^2k^3\delta_2(k)b_0k+\\28\tau_1|\tau|^4\xi_1^2\xi_2^4b_0^2k^
2\delta_1(k)b_0^2k^4\delta_2(k)b_0+28\tau_1|\tau|^4\xi_1^2\xi_2^4b_0^2k^2\delta_2(k)b_0^2k^3\delta_1(k)b_0k+\\28\tau_1|\tau|^4\xi_1^2\xi_2^4b_0^2k^
2\delta_2(k)b_0^2k^4\delta_1(k)b_0+28\tau_1|\tau|^4\xi_1^2\xi_2^4b_0^2k^3\delta_1(k)b_0^2k^2\delta_2(k)b_0k+\\28\tau_1|\tau|^4\xi_1^2\xi_2^4b_0^2k^
3\delta_1(k)b_0^2k^3\delta_2(k)b_0+28\tau_1|\tau|^4\xi_1^2\xi_2^4b_0^2k^3\delta_2(k)b_0^2k^2\delta_1(k)b_0k+\\28\tau_1|\tau|^4\xi_1^2\xi_2^4b_0^2k^
3\delta_2(k)b_0^2k^3\delta_1(k)b_0+16\tau_1|\tau|^4\xi_1\xi_2^5b_0^3k^4\delta_1(k)b_0k\delta_1(k)b_0k+\\16\tau_1|\tau|^4\xi_1\xi_2^5b_0^3k^4\delta_1(k)b_0k^
2\delta_1(k)b_0+16\tau_1|\tau|^4\xi_1\xi_2^5b_0^3k^5\delta_1(k)b_0\delta_1(k)b_0k+\\16\tau_1|\tau|^4\xi_1\xi_2^5b_0^3k^5\delta_1(k)b_0k\delta_1(k)b_0+
8\tau_1|\tau|^4\xi_1\xi_2^5b_0^2k^2\delta_1(k)b_0^2k^3\delta_1(k)b_0k+\\8\tau_1|\tau|^4\xi_1\xi_2^5b_0^2k^2\delta_1(k)b_0^2k^4\delta_1(k)b_0+8\tau_1|\tau|^
4\xi_1\xi_2^5b_0^2k^3\delta_1(k)b_0^2k^2\delta_1(k)b_0k+\\8\tau_1|\tau|^4\xi_1\xi_2^5b_0^2k^3\delta_1(k)b_0^2k^3\delta_1(k)b_0+104\tau_1^2|\tau|^4\xi_1^
2\xi_2^4b_0^3k^4\delta_2(k)b_0k\delta_2(k)b_0k+\\104\tau_1^2|\tau|^4\xi_1^2\xi_2^4b_0^3k^4\delta_2(k)b_0k^2\delta_2(k)b_0+104\tau_1^2|\tau|^4\xi_1^2\xi_2^
4b_0^3k^5\delta_2(k)b_0\delta_2(k)b_0k+\\104\tau_1^2|\tau|^4\xi_1^2\xi_2^4b_0^3k^5\delta_2(k)b_0k\delta_2(k)b_0+52\tau_1^2|\tau|^4\xi_1^2\xi_2^4b_0^2k^
2\delta_2(k)b_0^2k^3\delta_2(k)b_0k+\\52\tau_1^2|\tau|^4\xi_1^2\xi_2^4b_0^2k^2\delta_2(k)b_0^2k^4\delta_2(k)b_0+52\tau_1^2|\tau|^4\xi_1^2\xi_2^4b_0^2k^
3\delta_2(k)b_0^2k^2\delta_2(k)b_0k+\\52\tau_1^2|\tau|^4\xi_1^2\xi_2^4b_0^2k^3\delta_2(k)b_0^2k^3\delta_2(k)b_0+40\tau_1^2|\tau|^4\xi_1\xi_2^5b_0^3k^
4\delta_1(k)b_0k\delta_2(k)b_0k+\\40\tau_1^2|\tau|^4\xi_1\xi_2^5b_0^3k^4\delta_1(k)b_0k^2\delta_2(k)b_0+40\tau_1^2|\tau|^4\xi_1\xi_2^5b_0^3k^
4\delta_2(k)b_0k\delta_1(k)b_0k+\\40\tau_1^2|\tau|^4\xi_1\xi_2^5b_0^3k^4\delta_2(k)b_0k^2\delta_1(k)b_0+40\tau_1^2|\tau|^4\xi_1\xi_2^5b_0^3k^
5\delta_1(k)b_0\delta_2(k)b_0k+\\40\tau_1^2|\tau|^4\xi_1\xi_2^5b_0^3k^5\delta_1(k)b_0k\delta_2(k)b_0+40\tau_1^2|\tau|^4\xi_1\xi_2^5b_0^3k^
5\delta_2(k)b_0\delta_1(k)b_0k+\\40\tau_1^2|\tau|^4\xi_1\xi_2^5b_0^3k^5\delta_2(k)b_0k\delta_1(k)b_0+20\tau_1^2|\tau|^4\xi_1\xi_2^5b_0^2k^2\delta_1(k)b_0^2k^
3\delta_2(k)b_0k+\\20\tau_1^2|\tau|^4\xi_1\xi_2^5b_0^2k^2\delta_1(k)b_0^2k^4\delta_2(k)b_0+20\tau_1^2|\tau|^4\xi_1\xi_2^5b_0^2k^2\delta_2(k)b_0^2k^
3\delta_1(k)b_0k+\\20\tau_1^2|\tau|^4\xi_1\xi_2^5b_0^2k^2\delta_2(k)b_0^2k^4\delta_1(k)b_0+20\tau_1^2|\tau|^4\xi_1\xi_2^5b_0^2k^3\delta_1(k)b_0^2k^
2\delta_2(k)b_0k+\\20\tau_1^2|\tau|^4\xi_1\xi_2^5b_0^2k^3\delta_1(k)b_0^2k^3\delta_2(k)b_0+20\tau_1^2|\tau|^4\xi_1\xi_2^5b_0^2k^3\delta_2(k)b_0^2k^
2\delta_1(k)b_0k+\\20\tau_1^2|\tau|^4\xi_1\xi_2^5b_0^2k^3\delta_2(k)b_0^2k^3\delta_1(k)b_0+8\tau_1^2|\tau|^4\xi_2^6b_0^3k^4\delta_1(k)b_0k\delta_1(k)b_0k+\\
8\tau_1^2|\tau|^4\xi_2^6b_0^3k^4\delta_1(k)b_0k^2\delta_1(k)b_0+8\tau_1^2|\tau|^4\xi_2^6b_0^3k^5\delta_1(k)b_0\delta_1(k)b_0k+\\8\tau_1^2|\tau|^4\xi_2^6b_0^3k^
5\delta_1(k)b_0k\delta_1(k)b_0+4\tau_1^2|\tau|^4\xi_2^6b_0^2k^2\delta_1(k)b_0^2k^3\delta_1(k)b_0k+\\4\tau_1^2|\tau|^4\xi_2^6b_0^2k^2\delta_1(k)b_0^2k^
4\delta_1(k)b_0+4\tau_1^2|\tau|^4\xi_2^6b_0^2k^3\delta_1(k)b_0^2k^2\delta_1(k)b_0k+\\4\tau_1^2|\tau|^4\xi_2^6b_0^2k^3\delta_1(k)b_0^2k^3\delta_1(k)b_0+16|\tau|^
6\xi_1^2\xi_2^4b_0^3k^4\delta_2(k)b_0k\delta_2(k)b_0k+\\16|\tau|^6\xi_1^2\xi_2^4b_0^3k^4\delta_2(k)b_0k^2\delta_2(k)b_0+16|\tau|^6\xi_1^2\xi_2^4b_0^3k^
5\delta_2(k)b_0\delta_2(k)b_0k+\\16|\tau|^6\xi_1^2\xi_2^4b_0^3k^5\delta_2(k)b_0k\delta_2(k)b_0+8|\tau|^6\xi_1^2\xi_2^4b_0^2k^2\delta_2(k)b_0^2k^
3\delta_2(k)b_0k+\\8|\tau|^6\xi_1^2\xi_2^4b_0^2k^2\delta_2(k)b_0^2k^4\delta_2(k)b_0+8|\tau|^6\xi_1^2\xi_2^4b_0^2k^3\delta_2(k)b_0^2k^2\delta_2(k)b_0k+\\8|\tau|^
6\xi_1^2\xi_2^4b_0^2k^3\delta_2(k)b_0^2k^3\delta_2(k)b_0+8|\tau|^6\xi_1\xi_2^5b_0^3k^4\delta_1(k)b_0k\delta_2(k)b_0k+\\8|\tau|^6\xi_1\xi_2^5b_0^3k^
4\delta_1(k)b_0k^2\delta_2(k)b_0+8|\tau|^6\xi_1\xi_2^5b_0^3k^4\delta_2(k)b_0k\delta_1(k)b_0k+\\8|\tau|^6\xi_1\xi_2^5b_0^3k^4\delta_2(k)b_0k^2\delta_1(k)b_0+
8|\tau|^6\xi_1\xi_2^5b_0^3k^5\delta_1(k)b_0\delta_2(k)b_0k+\\8|\tau|^6\xi_1\xi_2^5b_0^3k^5\delta_1(k)b_0k\delta_2(k)b_0+8|\tau|^6\xi_1\xi_2^5b_0^3k^
5\delta_2(k)b_0\delta_1(k)b_0k+\\8|\tau|^6\xi_1\xi_2^5b_0^3k^5\delta_2(k)b_0k\delta_1(k)b_0+4|\tau|^6\xi_1\xi_2^5b_0^2k^2\delta_1(k)b_0^2k^3\delta_2(k)b_0k+\\
4|\tau|^6\xi_1\xi_2^5b_0^2k^2\delta_1(k)b_0^2k^4\delta_2(k)b_0+4|\tau|^6\xi_1\xi_2^5b_0^2k^2\delta_2(k)b_0^2k^3\delta_1(k)b_0k+\\4|\tau|^6\xi_1\xi_2^5b_0^2k^
2\delta_2(k)b_0^2k^4\delta_1(k)b_0+4|\tau|^6\xi_1\xi_2^5b_0^2k^3\delta_1(k)b_0^2k^2\delta_2(k)b_0k+\\4|\tau|^6\xi_1\xi_2^5b_0^2k^3\delta_1(k)b_0^2k^
3\delta_2(k)b_0+4|\tau|^6\xi_1\xi_2^5b_0^2k^3\delta_2(k)b_0^2k^2\delta_1(k)b_0k+\\4|\tau|^6\xi_1\xi_2^5b_0^2k^3\delta_2(k)b_0^2k^3\delta_1(k)b_0+48\tau_1|\tau|^
6\xi_1\xi_2^5b_0^3k^4\delta_2(k)b_0k\delta_2(k)b_0k+\\48\tau_1|\tau|^6\xi_1\xi_2^5b_0^3k^4\delta_2(k)b_0k^2\delta_2(k)b_0+48\tau_1|\tau|^6\xi_1\xi_2^5b_0^
3k^5\delta_2(k)b_0\delta_2(k)b_0k+\\48\tau_1|\tau|^6\xi_1\xi_2^5b_0^3k^5\delta_2(k)b_0k\delta_2(k)b_0+24\tau_1|\tau|^6\xi_1\xi_2^5b_0^2k^2\delta_2(k)b_0^2k^
3\delta_2(k)b_0k+\\24\tau_1|\tau|^6\xi_1\xi_2^5b_0^2k^2\delta_2(k)b_0^2k^4\delta_2(k)b_0+24\tau_1|\tau|^6\xi_1\xi_2^5b_0^2k^3\delta_2(k)b_0^2k^
2\delta_2(k)b_0k+\\24\tau_1|\tau|^6\xi_1\xi_2^5b_0^2k^3\delta_2(k)b_0^2k^3\delta_2(k)b_0+8\tau_1|\tau|^6\xi_2^6b_0^3k^4\delta_1(k)b_0k\delta_2(k)b_0k+\\
8\tau_1|\tau|^6\xi_2^6b_0^3k^4\delta_1(k)b_0k^2\delta_2(k)b_0+8\tau_1|\tau|^6\xi_2^6b_0^3k^4\delta_2(k)b_0k\delta_1(k)b_0k+\\8\tau_1|\tau|^6\xi_2^6b_0^3k^
4\delta_2(k)b_0k^2\delta_1(k)b_0+8\tau_1|\tau|^6\xi_2^6b_0^3k^5\delta_1(k)b_0\delta_2(k)b_0k+\\8\tau_1|\tau|^6\xi_2^6b_0^3k^5\delta_1(k)b_0k\delta_2(k)b_0+
8\tau_1|\tau|^6\xi_2^6b_0^3k^5\delta_2(k)b_0\delta_1(k)b_0k+\\8\tau_1|\tau|^6\xi_2^6b_0^3k^5\delta_2(k)b_0k\delta_1(k)b_0+4\tau_1|\tau|^6\xi_2^6b_0^2k^
2\delta_1(k)b_0^2k^3\delta_2(k)b_0k+\\4\tau_1|\tau|^6\xi_2^6b_0^2k^2\delta_1(k)b_0^2k^4\delta_2(k)b_0+4\tau_1|\tau|^6\xi_2^6b_0^2k^2\delta_2(k)b_0^2k^
3\delta_1(k)b_0k+\\4\tau_1|\tau|^6\xi_2^6b_0^2k^2\delta_2(k)b_0^2k^4\delta_1(k)b_0+4\tau_1|\tau|^6\xi_2^6b_0^2k^3\delta_1(k)b_0^2k^2\delta_2(k)b_0k+\\
4\tau_1|\tau|^6\xi_2^6b_0^2k^3\delta_1(k)b_0^2k^3\delta_2(k)b_0+4\tau_1|\tau|^6\xi_2^6b_0^2k^3\delta_2(k)b_0^2k^2\delta_1(k)b_0k+\\4\tau_1|\tau|^6\xi_2^6b_0^
2k^3\delta_2(k)b_0^2k^3\delta_1(k)b_0+8|\tau|^8\xi_2^6b_0^3k^4\delta_2(k)b_0k\delta_2(k)b_0k+\\8|\tau|^8\xi_2^6b_0^3k^4\delta_2(k)b_0k^2\delta_2(k)b_0+8|\tau|^
8\xi_2^6b_0^3k^5\delta_2(k)b_0\delta_2(k)b_0k+\\8|\tau|^8\xi_2^6b_0^3k^5\delta_2(k)b_0k\delta_2(k)b_0+4|\tau|^8\xi_2^6b_0^2k^2\delta_2(k)b_0^2k^
3\delta_2(k)b_0k+\\
4|\tau|^8\xi_2^6b_0^2k^2\delta_2(k)b_0^2k^4\delta_2(k)b_0+4|\tau|^8\xi_2^6b_0^2k^3\delta_2(k)b_0^2k^2\delta_2(k)b_0k+\\4|\tau|^8\xi_2^6b_0^
2k^3\delta_2(k)b_0^2k^3\delta_2(k)b_0.$\\

Invoking the trace property of $\tau_0$, we can permute the factors of each term cyclically
and obtain the following terms for computing $\zeta(0)$\\
$ $\\
$-b_0^2k\delta_1^2(k)-2\tau_1b_0^2k\delta_1\delta_2(k)-|\tau|^2b_0^2k\delta_2^2(k)+
(6\xi_1^2+12\tau_1\xi_1\xi_2+4\tau_1^2\xi_2^2+2|\tau|^2\xi_2^2)b_0^3k^3\delta_1^2(k)+\\
(6\xi_1^2+12\tau_1\xi_1\xi_2+4\tau_1^2\xi_2^2+2|\tau|^2\xi_2^2)b_0^3k^2\delta_1(k)^2+
(2\xi_1^2+4\tau_1\xi_1\xi_2+2|\tau|^2\xi_2^2)b_0^2k^2\delta_1(k)b_0\delta_1(k)+\\
(6\xi_1^2+12\tau_1\xi_1\xi_2+4\tau_1^2\xi_2^2+2|\tau|^2\xi_2^2)b_0^2k\delta_1(k)b_0k\delta_1(k)+\\
(12\tau_1\xi_1^2+16\tau_1^2\xi_1\xi_2+8|\tau|^2\xi_1\xi_2+12\tau_1|\tau|^2\xi_2^2)b_0^3k^3\delta_1\delta_2(k)+\\
(6\tau_1\xi_1^2+8\tau_1^2\xi_1\xi_2+4|\tau|^2\xi_1\xi_2+6\tau_1|\tau|^2\xi_2^2)b_0^3k^2\delta_1(k)\delta_2(k)+\\
(6\tau_1\xi_1^2+8\tau_1^2\xi_1\xi_2+4|\tau|^2\xi_1\xi_2+6\tau_1|\tau|^2\xi_2^2)b_0^3k^2\delta_2(k)\delta_1(k)+\\
(2\tau_1\xi_1^2+4\tau_1^2\xi_1\xi_2+2\tau_1|\tau|^2\xi_2^2)b_0^2k^2\delta_1(k)b_0\delta_2(k)+\\
(2\tau_1\xi_1^2+4\tau_1^2\xi_1\xi_2+2\tau_1|\tau|^2\xi_2^2)b_0^2k^2\delta_2(k)b_0\delta_1(k)+\\
(6\tau_1\xi_1^2+8\tau_1^2\xi_1\xi_2+4|\tau|^2\xi_1\xi_2+6\tau_1|\tau|^2\xi_2^2)b_0^2k\delta_1(k)b_0k\delta_2(k)+\\
(6\tau_1\xi_1^2+8\tau_1^2\xi_1\xi_2+4|\tau|^2\xi_1\xi_2+6\tau_1|\tau|^2\xi_2^2)b_0^2k\delta_2(k)b_0k\delta_1(k)+\\
(4\tau_1^2\xi_1^2+2|\tau|^2\xi_1^2+12\tau_1|\tau|^2\xi_1\xi_2+6|\tau|^4\xi_2^2)b_0^3k^3\delta_2^2(k)+\\
(4\tau_1^2\xi_1^2+2|\tau|^2\xi_1^2+12\tau_1|\tau|^2\xi_1\xi_2+6|\tau|^4\xi_2^2)b_0^3k^2\delta_2(k)^2+\\
(4\tau_1^2\xi_1^2+2|\tau|^2\xi_1^2+12\tau_1|\tau|^2\xi_1\xi_2+6|\tau|^4\xi_2^2)b_0^2k\delta_2(k)b_0k\delta_2(k)+\\
(2|\tau|^2\xi_1^2+4\tau_1|\tau|^2\xi_1\xi_2+2|\tau|^4\xi_2^2)b_0^2k^2\delta_2(k)b_0\delta_2(k)+\\
(-8\xi_1^4-32\tau_1\xi_1^3\xi_2-40\tau_1^2\xi_1^2\xi_2^2-8|\tau|^2\xi_1^2\xi_2^2-16\tau_1^3\xi_1\xi_2^3-16\tau_1|\tau|^2\xi_1\xi_2^3-8\tau_1^2|\tau|^2\xi_2^4)
b_0^4k^5\delta_1^2(k)+\\
(-8\xi_1^4-32\tau_1\xi_1^3\xi_2-40\tau_1^2\xi_1^2\xi_2^2-8|\tau|^2\xi_1^2\xi_2^2-16\tau_1^3\xi_1\xi_2^3-16\tau_1|\tau|^2\xi_1\xi_2^3-8\tau_1^2|\tau|^2\xi_2^4)
b_0^4k^4\delta_1(k)^2+\\
(-10\xi_1^4-40\tau_1\xi_1^3\xi_2-48\tau_1^2\xi_1^2\xi_2^2-12|\tau|^2\xi_1^2\xi_2^2-16\tau_1^3\xi_1\xi_2^3-24\tau_1|\tau|^2\xi_1\xi_2^3-8\tau_1^2|\tau|^2\xi_2^4
-\\
2|\tau|^4\xi_2^4)b_0^3k^4\delta_1(k)b_0\delta_1(k)+\\
(-20\xi_1^4-80\tau_1\xi_1^3\xi_2-96\tau_1^2\xi_1^2\xi_2^2-24|\tau|^2\xi_1^2\xi_2^2-32\tau_1^3\xi_1\xi_2^3-48\tau_1|\tau|^2\xi_1\xi_2^3-16\tau_1^2|\tau|^2\xi_2^4-\\
4|\tau|^4\xi_2^4)b_0^3k^3\delta_1(k)b_0k\delta_1(k)+\\
(-10\xi_1^4-40\tau_1\xi_1^3\xi_2-48\tau_1^2\xi_1^2\xi_2^2-12|\tau|^2\xi_1^2\xi_2^2-16\tau_1^3\xi_1\xi_2^3-24\tau_1|\tau|^2\xi_1\xi_2^3-8\tau_1^2|\tau|^2\xi_2^4
-\\
2|\tau|^4\xi_2^4)b_0^3k^2\delta_1(k)b_0k^2\delta_1(k)+\\
(-4\xi_1^4-16\tau_1\xi_1^3\xi_2-20\tau_1^2\xi_1^2\xi_2^2-4|\tau|^2\xi_1^2\xi_2^2-8\tau_1^3\xi_1\xi_2^3-8\tau_1|\tau|^2\xi_1\xi_2^3-\\
4\tau_1^2|\tau|^2\xi_2^4)
b_0^2k^2\delta_1(k)b_0^2k^2\delta_1(k)+\\
(-4\xi_1^4-16\tau_1\xi_1^3\xi_2-20\tau_1^2\xi_1^2\xi_2^2-4|\tau|^2\xi_1^2\xi_2^2-8\tau_1^3\xi_1\xi_2^3-8\tau_1|\tau|^2\xi_1\xi_2^3-\\
4\tau_1^2|\tau|^2\xi_2^4)b_0^2k\delta_1(k)b_0^2k^3\delta_1(k)+\\
(-16\tau_1\xi_1^4-48\tau_1^2\xi_1^3\xi_2-16|\tau|^2\xi_1^3\xi_2-32\tau_1^3\xi_1^2\xi_2^2-64\tau_1|\tau|^2\xi_1^2\xi_2^2-48\tau_1^2|\tau|^2\xi_1\xi_2^3-16|\tau|^4\xi_1\xi_2^3-\\
16\tau_1|\tau|^4\xi_2^4)b_0^4k^5\delta_1\delta_2(k)+\\
(-8\tau_1\xi_1^4-24\tau_1^2\xi_1^3\xi_2-8|\tau|^2\xi_1^3\xi_2-16\tau_1^3\xi_1^2\xi_2^2-32\tau_1|\tau|^2\xi_1^2\xi_2^2-24\tau_1^2|\tau|^2\xi_1\xi_2^3-8|\tau|^4\xi_1\xi_2^3-\\
8\tau_1|\tau|^4\xi_2^4)b_0^4k^4\delta_1(k)\delta_2(k)+\\
(-8\tau_1\xi_1^4-24\tau_1^2\xi_1^3\xi_2-8|\tau|^2\xi_1^3\xi_2-16\tau_1^3\xi_1^2\xi_2^2-32\tau_1|\tau|^2\xi_1^2\xi_2^2-24\tau_1^2|\tau|^2\xi_1\xi_2^3-8|\tau|^4\xi_1\xi_2^3-\\
8\tau_1|\tau|^4\xi_2^4)b_0^4k^4\delta_2(k)\delta_1(k)+\\
(-10\tau_1\xi_1^4-32\tau_1^2\xi_1^3\xi_2-8|\tau|^2\xi_1^3\xi_2-24\tau_1^3\xi_1^2\xi_2^2-36\tau_1|\tau|^2\xi_1^2\xi_2^2-32\tau_1^2|\tau|^2\xi_1\xi_2^3-8|\tau|^4\xi_1\xi_2^3-\\
10\tau_1|\tau|^4\xi_2^4)b_0^3k^4\delta_1(k)b_0\delta_2(k)+\\
(-10\tau_1\xi_1^4-32\tau_1^2\xi_1^3\xi_2-8|\tau|^2\xi_1^3\xi_2-24\tau_1^3\xi_1^2\xi_2^2-36\tau_1|\tau|^2\xi_1^2\xi_2^2-32\tau_1^2|\tau|^2\xi_1\xi_2^3-8|\tau|^4\xi_1\xi_2^3-\\
10\tau_1|\tau|^4\xi_2^4)b_0^3k^4\delta_2(k)b_0\delta_1(k)+\\
(-20\tau_1\xi_1^4-64\tau_1^2\xi_1^3\xi_2-16|\tau|^2\xi_1^3\xi_2-48\tau_1^3\xi_1^2\xi_2^2-72\tau_1|\tau|^2\xi_1^2\xi_2^2-64\tau_1^2|\tau|^2\xi_1\xi_2^3-16|\tau|^4\xi_1\xi_2^3-\\
20\tau_1|\tau|^4\xi_2^4)b_0^3k^3\delta_1(k)b_0k\delta_2(k)+\\
(-20\tau_1\xi_1^4-64\tau_1^2\xi_1^3\xi_2-16|\tau|^2\xi_1^3\xi_2-48\tau_1^3\xi_1^2\xi_2^2-72\tau_1|\tau|^2\xi_1^2\xi_2^2-64\tau_1^2|\tau|^2\xi_1\xi_2^3-16|\tau|^4\xi_1\xi_2^3-\\
20\tau_1|\tau|^4\xi_2^4)b_0^3k^3\delta_2(k)b_0k\delta_1(k)+\\
(-10\tau_1\xi_1^4-32\tau_1^2\xi_1^3\xi_2-8|\tau|^2\xi_1^3\xi_2-24\tau_1^3\xi_1^2\xi_2^2-36\tau_1|\tau|^2\xi_1^2\xi_2^2-32\tau_1^2|\tau|^2\xi_1\xi_2^3-8|\tau|^4\xi_1\xi_2^3-\\
10\tau_1|\tau|^4\xi_2^4)b_0^3k^2\delta_1(k)b_0k^2\delta_2(k)+\\
(-10\tau_1\xi_1^4-32\tau_1^2\xi_1^3\xi_2-8|\tau|^2\xi_1^3\xi_2-24\tau_1^3\xi_1^2\xi_2^2-36\tau_1|\tau|^2\xi_1^2\xi_2^2-32\tau_1^2|\tau|^2\xi_1\xi_2^3-8|\tau|^4\xi_1\xi_2^3-\\
10\tau_1|\tau|^4\xi_2^4)b_0^3k^2\delta_2(k)b_0k^2\delta_1(k)+\\
(-8\tau_1\xi_1^4-24\tau_1^2\xi_1^3\xi_2-8|\tau|^2\xi_1^3\xi_2-16\tau_1^3\xi_1^2\xi_2^2-32\tau_1|\tau|^2\xi_1^2\xi_2^2
-24\tau_1^2|\tau|^2\xi_1\xi_2^3-8|\tau|^4\xi_1\xi_2^3-\\
8\tau_1|\tau|^4\xi_2^4)b_0^2k^2\delta_1(k)b_0^2k^2\delta_2(k)+\\
(-4\tau_1\xi_1^4-12\tau_1^2\xi_1^3\xi_2-4|\tau|^2\xi_1^3\xi_2-8\tau_1^3\xi_1^2\xi_2^2-16\tau_1|\tau|^2\xi_1^2\xi_2^2-12\tau_1^2|\tau|^2\xi_1\xi_2^3+
-4|\tau|^4\xi_1\xi_2^3-\\
4\tau_1|\tau|^4\xi_2^4)b_0^2k\delta_1(k)b_0^2k^3\delta_2(k)+\\
(-4\tau_1\xi_1^4-12\tau_1^2\xi_1^3\xi_2-4|\tau|^2\xi_1^3\xi_2-8\tau_1^3\xi_1^2\xi_2^2-16\tau_1|\tau|^2\xi_1^2\xi_2^2-12\tau_1^2|\tau|^2\xi_1\xi_2^3
-4|\tau|^4\xi_1\xi_2^3-\\
4\tau_1|\tau|^4\xi_2^4)b_0^2k\delta_2(k)b_0^2k^3\delta_1(k)+\\
(-8\tau_1^2\xi_1^4-16\tau_1^3\xi_1^3\xi_2-16\tau_1|\tau|^2\xi_1^3\xi_2-40\tau_1^2|\tau|^2\xi_1^2\xi_2^2-8|\tau|^4\xi_1^2\xi_2^2-32\tau_1|\tau|^4\xi_1\xi_2^3-8|\tau|^6\xi_2^4)b_0^4k^5\delta_2^2(k)+
(-8\tau_1^2\xi_1^4-16\tau_1^3\xi_1^3\xi_2-16\tau_1|\tau|^2\xi_1^3\xi_2-40\tau_1^2|\tau|^2\xi_1^2\xi_2^2-8|\tau|^4\xi_1^2\xi_2^2-32\tau_1|\tau|^4\xi_1\xi_2^3-8|\tau|^6\xi_2^4)b_0^4k^4\delta_2(k)^2+
(-8\tau_1^2\xi_1^4-2|\tau|^2\xi_1^4-16\tau_1^3\xi_1^3\xi_2-24\tau_1|\tau|^2\xi_1^3\xi_2-48\tau_1^2|\tau|^2\xi_1^2\xi_2^2-12|\tau|^4\xi_1^2\xi_2^2-40\tau_1|\tau|^4\xi_1\xi_2^3-\\
10|\tau|^6\xi_2^4)b_0^3k^4\delta_2(k)b_0\delta_2(k)+ \\
(-16\tau_1^2\xi_1^4-4|\tau|^2\xi_1^4-32\tau_1^3\xi_1^3\xi_2-48\tau_1|\tau|^2\xi_1^3\xi_2-96\tau_1^2|\tau|^2\xi_1^2\xi_2^2-24|\tau|^4\xi_1^2\xi_2^2-80\tau_1|\tau|^4\xi_1\xi_2^3-\\
20|\tau|^6\xi_2^4)b_0^3k^3\delta_2(k)b_0k\delta_2(k)+\\
(-8\tau_1^2\xi_1^4-2|\tau|^2\xi_1^4-16\tau_1^3\xi_1^3\xi_2-24\tau_1|\tau|^2\xi_1^3\xi_2-48\tau_1^2|\tau|^2\xi_1^2\xi_2^2-12|\tau|^4\xi_1^2\xi_2^2-40\tau_1|\tau|^4\xi_1\xi_2^3-\\
10|\tau|^6\xi_2^4)b_0^3k^2\delta_2(k)b_0k^2\delta_2(k)+\\
(-4\tau_1^2\xi_1^4-8\tau_1^3\xi_1^3\xi_2-8\tau_1|\tau|^2\xi_1^3\xi_2-20\tau_1^2|\tau|^2\xi_1^2\xi_2^2-4|\tau|^4\xi_1^2\xi_2^2-16\tau_1|\tau|^4\xi_1\xi_2^3-4|\tau|^6\xi_2^4)b_0^2k^2\delta_2(k)b_0^2k^2\delta_2(k)+
(-4\tau_1^2\xi_1^4-8\tau_1^3\xi_1^3\xi_2-8\tau_1|\tau|^2\xi_1^3\xi_2-20\tau_1^2|\tau|^2\xi_1^2\xi_2^2-4|\tau|^4\xi_1^2\xi_2^2-16\tau_1|\tau|^4\xi_1\xi_2^3-4|\tau|^6\xi_2^4)b_0^2k\delta_2(k)b_0^2k^3\delta_2(k)+
(8\xi_1^6+48\tau_1\xi_1^5\xi_2+104\tau_1^2\xi_1^4\xi_2^2+16|\tau|^2\xi_1^4\xi_2^2+96\tau_1^3\xi_1^3\xi_2^3+64\tau_1|\tau|^2\xi_1^3\xi_2^3+32\tau_1^4\xi_1^2\xi_2^4+\\
80\tau_1^2|\tau|^2\xi_1^2\xi_2^4+8|\tau|^4\xi_1^2\xi_2^4+32\tau_1^3|\tau|^2\xi_1\xi_2^5+16\tau_1|\tau|^4\xi_1\xi_2^5+8\tau_1^2|\tau|^4\xi_2^6)b_0^4k^6\delta_1(k)b_0\delta_1(k)+\\
(16\xi_1^6+96\tau_1\xi_1^5\xi_2+(208)\tau_1^2\xi_1^4\xi_2^2+32|\tau|^2\xi_1^4\xi_2^2+192\tau_1^3\xi_1^3\xi_2^3+128\tau_1|\tau|^2\xi_1^3\xi_2^3+64\tau_1^4\xi_1^2\xi_2^4\\
+160\tau_1^2|\tau|^2\xi_1^2\xi_2^4+16|\tau|^4\xi_1^2\xi_2^4+64\tau_1^3|\tau|^2\xi_1\xi_2^5+32\tau_1|\tau|^4\xi_1\xi_2^5+16\tau_1^2|\tau|^4\xi_2^6)b_0^4k^5\delta_1(k)b_0k\delta_1(k)\\
+ (8\xi_1^6+48\tau_1\xi_1^5\xi_2+104\tau_1^2\xi_1^4\xi_2^2+16|\tau|^2\xi_1^4\xi_2^2+96\tau_1^3\xi_1^3\xi_2^3+64\tau_1|\tau|^2\xi_1^3\xi_2^3+32\tau_1^4\xi_1^2\xi_2^4+\\
80\tau_1^2|\tau|^2\xi_1^2\xi_2^4+8|\tau|^4\xi_1^2\xi_2^4+32\tau_1^3|\tau|^2\xi_1\xi_2^5+16\tau_1|\tau|^4\xi_1\xi_2^5+8\tau_1^2|\tau|^4\xi_2^6)b_0^4k^4\delta_1(k)b_0k^2\delta_1(k)+\\
(4\xi_1^6+(24)\tau_1\xi_1^5\xi_2+52\tau_1^2\xi_1^4\xi_2^2+8|\tau|^2\xi_1^4\xi_2^2+48\tau_1^3\xi_1^3\xi_2^3+32\tau_1|\tau|^2\xi_1^3\xi_2^3+16\tau_1^4\xi_1^2\xi_2^4+
\\
40\tau_1^2|\tau|^2\xi_1^2\xi_2^4+4|\tau|^4\xi_1^2\xi_2^4+16\tau_1^3|\tau|^2\xi_1\xi_2^5+8\tau_1|\tau|^4\xi_1\xi_2^5+4\tau_1^2|\tau|^4\xi_2^6)b_0^3k^4\delta_1(k)b_0^2k^2\delta_1(k)+\\
(8\xi_1^6+48\tau_1\xi_1^5\xi_2+104\tau_1^2\xi_1^4\xi_2^2+16|\tau|^2\xi_1^4\xi_2^2+96\tau_1^3\xi_1^3\xi_2^3+64\tau_1|\tau|^2\xi_1^3\xi_2^3+32\tau_1^4\xi_1^2\xi_2^4+\\
80\tau_1^2|\tau|^2\xi_1^2\xi_2^4+8|\tau|^4\xi_1^2\xi_2^4+32\tau_1^3|\tau|^2\xi_1\xi_2^5+16\tau_1|\tau|^4\xi_1\xi_2^5+8\tau_1^2|\tau|^4\xi_2^6)b_0^3k^3\delta_1(k)b_0^2k^3\delta_1(k)\\
+(4\xi_1^6+(24)\tau_1\xi_1^5\xi_2+52\tau_1^2\xi_1^4\xi_2^2+8|\tau|^2\xi_1^4\xi_2^2+48\tau_1^3\xi_1^3\xi_2^3+32\tau_1|\tau|^2\xi_1^3\xi_2^3+16\tau_1^4\xi_1^2\xi_2^4+
\\
40\tau_1^2|\tau|^2\xi_1^2\xi_2^4+4|\tau|^4\xi_1^2\xi_2^4+16\tau_1^3|\tau|^2\xi_1\xi_2^5+8\tau_1|\tau|^4\xi_1\xi_2^5+4\tau_1^2|\tau|^4\xi_2^6)b_0^3k^2\delta_1(k)b_0^2k^4\delta_1(k)+\\
(8\tau_1\xi_1^6+40\tau_1^2\xi_1^5\xi_2+8|\tau|^2\xi_1^5\xi_2+64\tau_1^3\xi_1^4\xi_2^2+56\tau_1|\tau|^2\xi_1^4\xi_2^2+32\tau_1^4\xi_1^3\xi_2^3+112\tau_1^2|\tau|^2\xi_1^3\xi_2^3+\\
16|\tau|^4\xi_1^3\xi_2^3+64\tau_1^3|\tau|^2\xi_1^2\xi_2^4+56\tau_1|\tau|^4\xi_1^2\xi_2^4+40\tau_1^2|\tau|^4\xi_1\xi_2^5+8|\tau|^6\xi_1\xi_2^5+\\
8\tau_1|\tau|^6\xi_2^6)b_0^4k^6\delta_1(k)b_0\delta_2(k)+\\
(8\tau_1\xi_1^6+40\tau_1^2\xi_1^5\xi_2+8|\tau|^2\xi_1^5\xi_2+64\tau_1^3\xi_1^4\xi_2^2+56\tau_1|\tau|^2\xi_1^4\xi_2^2+32\tau_1^4\xi_1^3\xi_2^3+112\tau_1^2|\tau|^2\xi_1^3\xi_2^3+\\
16|\tau|^4\xi_1^3\xi_2^3+64\tau_1^3|\tau|^2\xi_1^2\xi_2^4+56\tau_1|\tau|^4\xi_1^2\xi_2^4+40\tau_1^2|\tau|^4\xi_1\xi_2^5+8|\tau|^6\xi_1\xi_2^5+\\
8\tau_1|\tau|^6\xi_2^6)b_0^4k^6\delta_2(k)b_0\delta_1(k)+\\
(16\tau_1\xi_1^6+80\tau_1^2\xi_1^5\xi_2+16|\tau|^2\xi_1^5\xi_2+128\tau_1^3\xi_1^4\xi_2^2+112\tau_1|\tau|^2\xi_1^4\xi_2^2+64\tau_1^4\xi_1^3\xi_2^3+224\tau_1^2|\tau|^2\xi_1^3\xi_2^3+\\32|\tau|^4\xi_1^3\xi_2^3+128\tau_1^3|\tau|^2\xi_1^2\xi_2^4+112\tau_1|\tau|^4\xi_1^2\xi_2^4+80\tau_1^2|\tau|^4\xi_1\xi_2^5+16|\tau|^6\xi_1\xi_2^5+\\
16\tau_1|\tau|^6\xi_2^6)b_0^4k^5\delta_1(k)b_0k\delta_2(k)+\\
(16\tau_1\xi_1^6+80\tau_1^2\xi_1^5\xi_2+16|\tau|^2\xi_1^5\xi_2+128\tau_1^3\xi_1^4\xi_2^2+112\tau_1|\tau|^2\xi_1^4\xi_2^2+64\tau_1^4\xi_1^3\xi_2^3+224\tau_1^2|\tau|^2\xi_1^3\xi_2^3+\\
32|\tau|^4\xi_1^3\xi_2^3+128\tau_1^3|\tau|^2\xi_1^2\xi_2^4+112\tau_1|\tau|^4\xi_1^2\xi_2^4+80\tau_1^2|\tau|^4\xi_1\xi_2^5+16|\tau|^6\xi_1\xi_2^5+\\
16\tau_1|\tau|^6\xi_2^6)b_0^4k^5\delta_2(k)b_0k\delta_1(k)+ \\
(8\tau_1\xi_1^6+40\tau_1^2\xi_1^5\xi_2+8|\tau|^2\xi_1^5\xi_2+64\tau_1^3\xi_1^4\xi_2^2+56\tau_1|\tau|^2\xi_1^4\xi_2^2+32\tau_1^4\xi_1^3\xi_2^3+112\tau_1^2|\tau|^2\xi_1^3\xi_2^3+\\
16|\tau|^4\xi_1^3\xi_2^3+64\tau_1^3|\tau|^2\xi_1^2\xi_2^4+56\tau_1|\tau|^4\xi_1^2\xi_2^4+40\tau_1^2|\tau|^4\xi_1\xi_2^5+
8|\tau|^6\xi_1\xi_2^5+
\\8\tau_1|\tau|^6\xi_2^6)b_0^4k^4\delta_1(k)b_0k^2\delta_2(k)+ \\
(8\tau_1\xi_1^6+40\tau_1^2\xi_1^5\xi_2+8|\tau|^2\xi_1^5\xi_2+64\tau_1^3\xi_1^4\xi_2^2+56\tau_1|\tau|^2\xi_1^4\xi_2^2+32\tau_1^4\xi_1^3\xi_2^3+112\tau_1^2|\tau|^2\xi_1^3\xi_2^3+\\
16|\tau|^4\xi_1^3\xi_2^3+64\tau_1^3|\tau|^2\xi_1^2\xi_2^4+56\tau_1|\tau|^4\xi_1^2\xi_2^4+40\tau_1^2|\tau|^4\xi_1\xi_2^5+8|\tau|^6\xi_1\xi_2^5+\\
8\tau_1|\tau|^6\xi_2^6)b_0^4k^4\delta_2(k)b_0k^2\delta_1(k)+\\
(4\tau_1\xi_1^6+20\tau_1^2\xi_1^5\xi_2+4|\tau|^2\xi_1^5\xi_2+32\tau_1^3\xi_1^4\xi_2^2+28\tau_1|\tau|^2\xi_1^4\xi_2^2+16\tau_1^4\xi_1^3\xi_2^3+56\tau_1^2|\tau|^2\xi_1^3\xi_2^3+\\
8|\tau|^4\xi_1^3\xi_2^3+32\tau_1^3|\tau|^2\xi_1^2\xi_2^4+28\tau_1|\tau|^4\xi_1^2\xi_2^4+20\tau_1^2|\tau|^4\xi_1\xi_2^5+4|\tau|^6\xi_1\xi_2^5+\\4\tau_1|\tau|^6\xi_2^6)
b_0^3k^4\delta_1(k)b_0^2k^2\delta_2(k)+\\
(4\tau_1\xi_1^6+20\tau_1^2\xi_1^5\xi_2+4|\tau|^2\xi_1^5\xi_2+32\tau_1^3\xi_1^4\xi_2^2+28\tau_1|\tau|^2\xi_1^4\xi_2^2+16\tau_1^4\xi_1^3\xi_2^3+56\tau_1^2|\tau|^2\xi_1^3\xi_2^3+\\8|\tau|^4\xi_1^3\xi_2^3+32\tau_1^3|\tau|^2\xi_1^2\xi_2^4+28\tau_1|\tau|^4\xi_1^2\xi_2^4+20\tau_1^2|\tau|^4\xi_1\xi_2^5+4|\tau|^6\xi_1\xi_2^5+\\4\tau_1|\tau|^6\xi_2^6)
b_0^3k^4\delta_2(k)b_0^2k^2\delta_1(k)+\\
(8\tau_1\xi_1^6+40\tau_1^2\xi_1^5\xi_2+8|\tau|^2\xi_1^5\xi_2+64\tau_1^3\xi_1^4\xi_2^2+56\tau_1|\tau|^2\xi_1^4\xi_2^2+32\tau_1^4\xi_1^3\xi_2^3+112\tau_1^2|\tau|^2\xi_1^3\xi_2^3+\\16|\tau|^4\xi_1^3\xi_2^3+64\tau_1^3|\tau|^2\xi_1^2\xi_2^4+56\tau_1|\tau|^4\xi_1^2\xi_2^4+40\tau_1^2|\tau|^4\xi_1\xi_2^5+8|\tau|^6\xi_1\xi_2^5+\\
8\tau_1|\tau|^6\xi_2^6)b_0^3k^3\delta_1(k)b_0^2k^3\delta_2(k)+\\
(8\tau_1\xi_1^6+40\tau_1^2\xi_1^5\xi_2+8|\tau|^2\xi_1^5\xi_2+64\tau_1^3\xi_1^4\xi_2^2+56\tau_1|\tau|^2\xi_1^4\xi_2^2+32\tau_1^4\xi_1^3\xi_2^3+112\tau_1^2|\tau|^2\xi_1^3\xi_2^3+\\
16|\tau|^4\xi_1^3\xi_2^3+64\tau_1^3|\tau|^2\xi_1^2\xi_2^4+56\tau_1|\tau|^4\xi_1^2\xi_2^4+40\tau_1^2|\tau|^4\xi_1\xi_2^5+8|\tau|^6\xi_1\xi_2^5+\\8\tau_1|\tau|^6\xi_2^6)b_0^3k^3\delta_2(k)b_0^2k^3\delta_1(k)+\\
(4\tau_1\xi_1^6+20\tau_1^2\xi_1^5\xi_2+4|\tau|^2\xi_1^5\xi_2+32\tau_1^3\xi_1^4\xi_2^2+28\tau_1|\tau|^2\xi_1^4\xi_2^2+16\tau_1^4\xi_1^3\xi_2^3+56\tau_1^2|\tau|^2\xi_1^3\xi_2^3+\\
8|\tau|^4\xi_1^3\xi_2^3+32\tau_1^3|\tau|^2\xi_1^2\xi_2^4+28\tau_1|\tau|^4\xi_1^2\xi_2^4+20\tau_1^2|\tau|^4\xi_1\xi_2^5+
4|\tau|^6\xi_1\xi_2^5+\\4\tau_1|\tau|^6\xi_2^6)b_0^3k^2\delta_1(k)b_0^2k^4\delta_2(k)+\\
(4\tau_1\xi_1^6+20\tau_1^2\xi_1^5\xi_2+4|\tau|^2\xi_1^5\xi_2+32\tau_1^3\xi_1^4\xi_2^2+28\tau_1|\tau|^2\xi_1^4\xi_2^2+16\tau_1^4\xi_1^3\xi_2^3+56\tau_1^2|\tau|^2\xi_1^3\xi_2^3+\\
8|\tau|^4\xi_1^3\xi_2^3+
32\tau_1^3|\tau|^2\xi_1^2\xi_2^4+28\tau_1|\tau|^4\xi_1^2\xi_2^4+20\tau_1^2|\tau|^4\xi_1\xi_2^5+4|\tau|^6\xi_1\xi_2^5+\\
4\tau_1|\tau|^6\xi_2^6)b_0^3k^2\delta_2(k)b_0^2k^4\delta_1(k)+\\
(8\tau_1^2\xi_1^6+32\tau_1^3\xi_1^5\xi_2+16\tau_1|\tau|^2\xi_1^5\xi_2+32\tau_1^4\xi_1^4\xi_2^2+80\tau_1^2|\tau|^2\xi_1^4\xi_2^2+8|\tau|^4\xi_1^4\xi_2^2+96\tau_1^3|\tau|^2\xi_1^3\xi_2^3+\\
64\tau_1|\tau|^4\xi_1^3\xi_2^3+104\tau_1^2|\tau|^4\xi_1^2\xi_2^4+16|\tau|^6\xi_1^2\xi_2^4+48\tau_1|\tau|^6\xi_1\xi_2^5+8|\tau|^8\xi_2^6)b_0^4k^6\delta_2(k)b_0\delta_2(k)+\\
(16\tau_1^2\xi_1^6+64\tau_1^3\xi_1^5\xi_2+32\tau_1|\tau|^2\xi_1^5\xi_2+64\tau_1^4\xi_1^4\xi_2^2+160\tau_1^2|\tau|^2\xi_1^4\xi_2^2+16|\tau|^4\xi_1^4\xi_2^2+192\tau_1^3|\tau|^2\xi_1^3\xi_2^3
\\
+128\tau_1|\tau|^4\xi_1^3\xi_2^3+208\tau_1^2|\tau|^4\xi_1^2\xi_2^4+32|\tau|^6\xi_1^2\xi_2^4+96\tau_1|\tau|^6\xi_1\xi_2^5+16|\tau|^8\xi_2^6)b_0^4k^5\delta_2(k)b_0k\delta_2(k)+\\
(8\tau_1^2\xi_1^6+32\tau_1^3\xi_1^5\xi_2+16\tau_1|\tau|^2\xi_1^5\xi_2+32\tau_1^4\xi_1^4\xi_2^2+80\tau_1^2|\tau|^2\xi_1^4\xi_2^2+8|\tau|^4\xi_1^4\xi_2^2+96\tau_1^3|\tau|^2\xi_1^3\xi_2^3+64\tau_1|\tau|^4\xi_1^3\xi_2^3+104\tau_1^2|\tau|^4\xi_1^2\xi_2^4+16|\tau|^6\xi_1^2\xi_2^4+48\tau_1|\tau|^6\xi_1\xi_2^5+8|\tau|^8\xi_2^6)b_0^4k^4\delta_2(k)b_0k^2\delta_2(k)+\\
(4\tau_1^2\xi_1^6+16\tau_1^3\xi_1^5\xi_2+8\tau_1|\tau|^2\xi_1^5\xi_2+16\tau_1^4\xi_1^4\xi_2^2+40\tau_1^2|\tau|^2\xi_1^4\xi_2^2+4|\tau|^4\xi_1^4\xi_2^2+48\tau_1^3|\tau|^2\xi_1^3\xi_2^3+32\tau_1|\tau|^4\xi_1^3\xi_2^3+52\tau_1^2|\tau|^4\xi_1^2\xi_2^4+8|\tau|^6\xi_1^2\xi_2^4+(24)\tau_1|\tau|^6\xi_1\xi_2^5+4|\tau|^8\xi_2^6)b_0^3k^4\delta_2(k)b_0^2k^2\delta_2(k)+\\
(8\tau_1^2\xi_1^6+32\tau_1^3\xi_1^5\xi_2+16\tau_1|\tau|^2\xi_1^5\xi_2+32\tau_1^4\xi_1^4\xi_2^2+80\tau_1^2|\tau|^2\xi_1^4\xi_2^2+8|\tau|^4\xi_1^4\xi_2^2+96\tau_1^3|\tau|^2\xi_1^3\xi_2^3+64\tau_1|\tau|^4\xi_1^3\xi_2^3+104\tau_1^2|\tau|^4\xi_1^2\xi_2^4+16|\tau|^6\xi_1^2\xi_2^4+48\tau_1|\tau|^6\xi_1\xi_2^5+8|\tau|^8\xi_2^6)b_0^3k^3\delta_2(k)b_0^2k^3\delta_2(k)+\\
(4\tau_1^2\xi_1^6+16\tau_1^3\xi_1^5\xi_2+8\tau_1|\tau|^2\xi_1^5\xi_2+16\tau_1^4\xi_1^4\xi_2^2+40\tau_1^2|\tau|^2\xi_1^4\xi_2^2+4|\tau|^4\xi_1^4\xi_2^2+48\tau_1^3|\tau|^2\xi_1^3\xi_2^3+32\tau_1|\tau|^4\xi_1^3\xi_2^3+52\tau_1^2|\tau|^4\xi_1^2\xi_2^4+8|\tau|^6\xi_1^2\xi_2^4+(24)\tau_1|\tau|^6\xi_1\xi_2^5+4|\tau|^8\xi_2^6)b_0^3k^2\delta_2(k)b_0^2k^4\delta_2(k).$\\

To integrate these terms over the $\xi$-plane we pass to the
following coordinates

\[\xi_1= r \cos \theta - r \frac{\tau_1}{\tau_2} \sin \theta,
\,\,\,\,\,
\xi_2=\frac{r}{\tau_2}\sin \theta ,\]
where $\theta$ ranges from 0 to $2\pi$ and $r$ ranges from 0 to $\infty$.

After multiplication by the Jacobian of the change of coordinates
which is equal to $\frac{r}{\tau_2}$,
and integrating with respect to $\theta$ from 0 to $2\pi$, one gets,
up to an overall factor of
$\frac{2\pi}{\tau_2}$\\
$ $\\
$-rb_0^2k\delta_1^2(k)
-2\tau_1rb_0^2k\delta_1\delta_2(k)
-|\tau|^2rb_0^2k\delta_2^2(k)+
4r^3b_0^3k^3\delta_1^2(k)+
4r^3b_0^3k^2\delta_1(k)^2\\+
2r^3b_0^2k^2\delta_1(k)b_0\delta_1(k)+
4r^3b_0^2k\delta_1(k)b_0k\delta_1(k)+
8\tau_1r^3b_0^3k^3\delta_1\delta_2(k)
+4\tau_1r^3b_0^3k^2\delta_1(k)\delta_2(k)\\+
4\tau_1r^3b_0^3k^2\delta_2(k)\delta_1(k)+
2\tau_1r^3b_0^2k^2\delta_1(k)b_0\delta_2(k)
+2\tau_1r^3b_0^2k^2\delta_2(k)b_0\delta_1(k)+
4\tau_1r^3b_0^2k\delta_1(k)b_0k\delta_2(k)+\\
4\tau_1r^3b_0^2k\delta_2(k)b_0k\delta_1(k)+
4|\tau|^2r^3b_0^3k^3\delta_2^2(k)+
4|\tau|^2r^3b_0^3k^2\delta_2(k)^2+
4|\tau|^2r^3b_0^2k\delta_2(k)b_0k\delta_2(k)+\\
2|\tau|^2r^3b_0^2k^2\delta_2(k)b_0\delta_2(k)
-4r^5b_0^4k^5\delta_1^2(k)-
4r^5b_0^4k^4\delta_1(k)^2-
6r^5b_0^3k^4\delta_1(k)b_0\delta_1(k)-\\
12r^5b_0^3k^3\delta_1(k)b_0k\delta_1(k)-
6r^5b_0^3k^2\delta_1(k)b_0k^2\delta_1(k)-
2r^5b_0^2k^2\delta_1(k)b_0^2k^2\delta_1(k)-\\
2r^5b_0^2k\delta_1(k)b_0^2k^3\delta_1(k)-
8\tau_1r^5b_0^4k^5\delta_1\delta_2(k)-
4\tau_1r^5b_0^4k^4\delta_1(k)\delta_2(k)-
4\tau_1r^5b_0^4k^4\delta_2(k)\delta_1(k)-\\
6\tau_1r^5b_0^3k^4\delta_1(k)b_0\delta_2(k)-
6\tau_1r^5b_0^3k^4\delta_2(k)b_0\delta_1(k)-
12\tau_1r^5b_0^3k^3\delta_1(k)b_0k\delta_2(k)-\\
12\tau_1r^5b_0^3k^3\delta_2(k)b_0k\delta_1(k)-
6\tau_1r^5b_0^3k^2\delta_1(k)b_0k^2\delta_2(k)-
6\tau_1r^5b_0^3k^2\delta_2(k)b_0k^2\delta_1(k)-\\
4\tau_1r^5b_0^2k^2\delta_1(k)b_0^2k^2\delta_2(k)-
2\tau_1r^5b_0^2k\delta_1(k)b_0^2k^3\delta_2(k)-
2\tau_1r^5b_0^2k\delta_2(k)b_0^2k^3\delta_1(k)-\\
4|\tau|^2r^5b_0^4k^5\delta_2^2(k)-
4|\tau|^2r^5b_0^4k^4\delta_2(k)^2-
6|\tau|^2r^5b_0^3k^4\delta_2(k)b_0\delta_2(k)-
12|\tau|^2r^5b_0^3k^3\delta_2(k)b_0k\delta_2(k)-\\
6|\tau|^2r^5b_0^3k^2\delta_2(k)b_0k^2\delta_2(k)-
2|\tau|^2r^5b_0^2k^2\delta_2(k)b_0^2k^2\delta_2(k)-
2|\tau|^2r^5b_0^2k\delta_2(k)b_0^2k^3\delta_2(k)+\\
4r^7b_0^4k^6\delta_1(k)b_0\delta_1(k)+
8r^7b_0^4k^5\delta_1(k)b_0k\delta_1(k)+
4r^7b_0^4k^4\delta_1(k)b_0k^2\delta_1(k)+
2r^7b_0^3k^4\delta_1(k)b_0^2k^2\delta_1(k)+\\
4r^7b_0^3k^3\delta_1(k)b_0^2k^3\delta_1(k)+
2r^7b_0^3k^2\delta_1(k)b_0^2k^4\delta_1(k)+
4\tau_1r^7b_0^4k^6\delta_1(k)b_0\delta_2(k)+\\
4\tau_1r^7b_0^4k^6\delta_2(k)b_0\delta_1(k)+
8\tau_1r^7b_0^4k^5\delta_1(k)b_0k\delta_2(k)+
8\tau_1r^7b_0^4k^5\delta_2(k)b_0k\delta_1(k)+\\
4\tau_1r^7b_0^4k^4\delta_1(k)b_0k^2\delta_2(k)+
4\tau_1r^7b_0^4k^4\delta_2(k)b_0k^2\delta_1(k)+
2\tau_1r^7b_0^3k^4\delta_1(k)b_0^2k^2\delta_2(k)+\\
2\tau_1r^7b_0^3k^4\delta_2(k)b_0^2k^2\delta_1(k)+
4\tau_1r^7b_0^3k^3\delta_1(k)b_0^2k^3\delta_2(k)+
4\tau_1r^7b_0^3k^3\delta_2(k)b_0^2k^3\delta_1(k)+\\
2\tau_1r^7b_0^3k^2\delta_1(k)b_0^2k^4\delta_2(k)+
2\tau_1r^7b_0^3k^2\delta_2(k)b_0^2k^4\delta_1(k)+
4|\tau|^2r^7b_0^4k^6\delta_2(k)b_0\delta_2(k)+\\
8|\tau|^2r^7b_0^4k^5\delta_2(k)b_0k\delta_2(k)+
4|\tau|^2r^7b_0^4k^4\delta_2(k)b_0k^2\delta_2(k)+
2|\tau|^2r^7b_0^3k^4\delta_2(k)b_0^2k^2\delta_2(k)+\\
4|\tau|^2r^7b_0^3k^3\delta_2(k)b_0^2k^3\delta_2(k)+
2|\tau|^2r^7b_0^3k^2\delta_2(k)b_0^2k^4\delta_2(k),
$\\
$ $\\
where
\[b_0=(r^2k^2+1)^{-1}.\]

\subsection{Terms with all $b_0$ on the left.}

The integral $\int_{0}^{\infty}\bullet \,\, dr$ of these terms, namely\\
$ $\\
$-rb_0^2k\delta_1^2(k)-
2\tau_1rb_0^2k\delta_1\delta_2(k)-
|\tau|^2rb_0^2k\delta_2^2(k)+
4r^3b_0^3k^3\delta_1^2(k)+
4r^3b_0^3k^2\delta_1(k)^2+\\
8\tau_1r^3b_0^3k^3\delta_1\delta_2(k)+
4\tau_1r^3b_0^3k^2\delta_1(k)\delta_2(k)+
4\tau_1r^3b_0^3k^2\delta_2(k)\delta_1(k)+
4|\tau|^2r^3b_0^3k^3\delta_2^2(k)+\\
4|\tau|^2r^3b_0^3k^2\delta_2(k)^2-
4r^5b_0^4k^5\delta_1^2(k)-
4r^5b_0^4k^4\delta_1(k)^2-
8\tau_1r^5b_0^4k^5\delta_1\delta_2(k)-\\
4\tau_1r^5b_0^4k^4\delta_1(k)\delta_2(k)-
4\tau_1r^5b_0^4k^4\delta_2(k)\delta_1(k)-
4|\tau|^2r^5b_0^4k^5\delta_2^2(k)-
4|\tau|^2r^5b_0^4k^4\delta_2(k)^2
$\\
$ $\\
gives (up to an overall factor of $\frac{2\pi}{\tau_2}$)
\begin{eqnarray} \label{allb_0onleft}
&&-\frac{1}{6}k^{-1}\delta_1^2(k)-\frac{|\tau|^2}{6}k^{-1}\delta_2^2(k)-
\frac{\tau_1}{3}k^{-1}\delta_1\delta_2(k)+\frac{1}{3}k^{-2}\delta_1(k)^2
+\nonumber \\
&&\frac{|\tau|^2}{3}k^{-2}\delta_2(k)^2\newline
+\frac{\tau_1}{3}k^{-2}\delta_1(k)\delta_2(k)
+\frac{\tau_1}{3}k^{-2}\delta_2(k)\delta_1(k).
\end{eqnarray}

\subsection{Terms with $b_0$ in the middle.} \label{b_0inthemiddle}
They are the following terms\\
$ $\\
$2r^3b_0^2k^2\delta_1(k)b_0\delta_1(k)+
4r^3b_0^2k\delta_1(k)b_0k\delta_1(k)+
2\tau_1r^3b_0^2k^2\delta_1(k)b_0\delta_2(k)+\\
2\tau_1r^3b_0^2k^2\delta_2(k)b_0\delta_1(k)+
4\tau_1r^3b_0^2k\delta_1(k)b_0k\delta_2(k)+
4\tau_1r^3b_0^2k\delta_2(k)b_0k\delta_1(k)+\\
4|\tau|^2r^3b_0^2k\delta_2(k)b_0k\delta_2(k)+
2|\tau|^2r^3b_0^2k^2\delta_2(k)b_0\delta_2(k)-
6r^5b_0^3k^4\delta_1(k)b_0\delta_1(k)-\\
12r^5b_0^3k^3\delta_1(k)b_0k\delta_1(k)
-6r^5b_0^3k^2\delta_1(k)b_0k^2\delta_1(k)-
6\tau_1r^5b_0^3k^4\delta_1(k)b_0\delta_2(k)-\\
6\tau_1r^5b_0^3k^4\delta_2(k)b_0\delta_1(k)-
12\tau_1r^5b_0^3k^3\delta_1(k)b_0k\delta_2(k)-
12\tau_1r^5b_0^3k^3\delta_2(k)b_0k\delta_1(k)-\\
6\tau_1r^5b_0^3k^2\delta_1(k)b_0k^2\delta_2(k)
-6\tau_1r^5b_0^3k^2\delta_2(k)b_0k^2\delta_1(k)-
6|\tau|^2r^5b_0^3k^4\delta_2(k)b_0\delta_2(k)-\\
12|\tau|^2r^5b_0^3k^3\delta_2(k)b_0k\delta_2(k)-
6|\tau|^2r^5b_0^3k^2\delta_2(k)b_0k^2\delta_2(k)
+4r^7b_0^4k^6\delta_1(k)b_0\delta_1(k)+\\
8r^7b_0^4k^5\delta_1(k)b_0k\delta_1(k)+
4r^7b_0^4k^4\delta_1(k)b_0k^2\delta_1(k)+
4\tau_1r^7b_0^4k^6\delta_1(k)b_0\delta_2(k)+\\
4\tau_1r^7b_0^4k^6\delta_2(k)b_0\delta_1(k)+
8\tau_1r^7b_0^4k^5\delta_1(k)b_0k\delta_2(k)
+8\tau_1r^7b_0^4k^5\delta_2(k)b_0k\delta_1(k)+\\
4\tau_1r^7b_0^4k^4\delta_1(k)b_0k^2\delta_2(k)
+4\tau_1r^7b_0^4k^4\delta_2(k)b_0k^2\delta_1(k)+
4|\tau|^2r^7b_0^4k^6\delta_2(k)b_0\delta_2(k)+\\
8|\tau|^2r^7b_0^4k^5\delta_2(k)b_0k\delta_2(k)+
4|\tau|^2r^7b_0^4k^4\delta_2(k)b_0k^2\delta_2(k).
$

\subsection{The rest of the terms (with $b_0^2$ in the middle).}
These terms, namely\\
$ $\\
$-2r^5b_0^2k^2\delta_1(k)b_0^2k^2\delta_1(k)-
2r^5b_0^2k\delta_1(k)b_0^2k^3\delta_1(k)-
4\tau_1r^5b_0^2k^2\delta_1(k)b_0^2k^2\delta_2(k)-\\
2\tau_1r^5b_0^2k\delta_1(k)b_0^2k^3\delta_2(k)-
2\tau_1r^5b_0^2k\delta_2(k)b_0^2k^3\delta_1(k)-
2|\tau|^2r^5b_0^2k^2\delta_2(k)b_0^2k^2\delta_2(k)-\\
2|\tau|^2r^5b_0^2k\delta_2(k)b_0^2k^3\delta_2(k)+
2r^7b_0^3k^4\delta_1(k)b_0^2k^2\delta_1(k)+
4r^7b_0^3k^3\delta_1(k)b_0^2k^3\delta_1(k)+\\
2r^7b_0^3k^2\delta_1(k)b_0^2k^4\delta_1(k)+
2\tau_1r^7b_0^3k^4\delta_1(k)b_0^2k^2\delta_2(k)+
2\tau_1r^7b_0^3k^4\delta_2(k)b_0^2k^2\delta_1(k)+\\
4\tau_1r^7b_0^3k^3\delta_1(k)b_0^2k^3\delta_2(k)+
4\tau_1r^7b_0^3k^3\delta_2(k)b_0^2k^3\delta_1(k)+
2\tau_1r^7b_0^3k^2\delta_1(k)b_0^2k^4\delta_2(k)+\\
2\tau_1r^7b_0^3k^2\delta_2(k)b_0^2k^4\delta_1(k)
+2|\tau|^2r^7b_0^3k^4\delta_2(k)b_0^2k^2\delta_2(k)+
4|\tau|^2r^7b_0^3k^3\delta_2(k)b_0^2k^3\delta_2(k)+\\
2|\tau|^2r^7b_0^3k^2\delta_2(k)b_0^2k^4\delta_2(k)
$\\
$ $\\
can be replaced by terms with $b_0$ in the middle using
\[\partial_r(b_0)=-2rk^2b_0^2,\]
and integration by parts. For example
\[\int_{0}^{\infty}r^4b_0^2k^2\delta_1(k)rb_0^2k^2\delta_1(k)=
\frac{1}{2}\int_{0}^{\infty}\partial_r(r^4b_0^2)k^2\delta_1(k)b_0\delta_1(k). \]
After this replacement, and combining the result with the terms listed in Subsection
\ref{b_0inthemiddle},
we obtain the following terms (up to an overall factor of
$\frac{2\pi}{\tau_2}$)\\
$ $\\
$T=-2|\tau|^2r^7b_0^4k^4\delta_2(k)b_0k^2\delta_2(k)
-2\tau_1r^7b_0^4k^4\delta_2(k)b_0k^2\delta_1(k)
-2\tau_1r^7b_0^4k^4\delta_1(k)b_0k^2\delta_2(k)-\\
2r^7b_0^4k^4\delta_1(k)b_0k^2\delta_1(k)
-4|\tau|^2r^7b_0^4k^5\delta_2(k)b_0k\delta_2(k)
-4\tau_1r^7b_0^4k^5\delta_2(k)b_0k\delta_1(k)-\\
4\tau_1r^7b_0^4k^5\delta_1(k)b_0k\delta_2(k)
-4r^7b_0^4k^5\delta_1(k)b_0k\delta_1(k)
-2|\tau|^2r^7b_0^4k^6\delta_2(k)b_0\delta_2(k)-\\
2\tau_1r^7b_0^4k^6\delta_2(k)b_0\delta_1(k)
-2\tau_1r^7b_0^4k^6\delta_1(k)b_0\delta_2(k)
-2r^7b_0^4k^6\delta_1(k)b_0\delta_1(k)+\\
4|\tau|^2r^5b_0^3k^3\delta_2(k)b_0k\delta_2(k)
+4\tau_1r^5b_0^3k^3\delta_2(k)b_0k\delta_1(k)
+4\tau_1r^5b_0^3k^3\delta_1(k)b_0k\delta_2(k)+\\
4r^5b_0^3k^3\delta_1(k)b_0k\delta_1(k)
+4|\tau|^2r^5b_0^3k^4\delta_2(k)b_0\delta_2(k)
+4r^5b_0^3k^4\delta_1(k)b_0\delta_1(k)+\\
8\tau_1r^5b_0^3k^4\delta_1(k)b_0\delta_2(k)
-2|\tau|^2r^3b_0^2k^2\delta_2(k)b_0\delta_2(k)
+2\tau_1r^3b_0^2k^2\delta_2(k)b_0\delta_1(k)-\\
6\tau_1r^3b_0^2k^2\delta_1(k)b_0\delta_2(k)
-2r^3b_0^2k^2\delta_1(k)b_0\delta_1(k).$

\subsection{Integrating the terms with $b_0$ in the middle.} \label{b_0inmiddle}
The computation of $\int_{0}^{\infty} \bullet \, dr$ of these terms is
achieved by the following lemma of Connes and Tretkoff proved
in \cite{contre}.

\newtheorem{integration}[main1]{Lemma}
\begin{integration} \label{integration}
For any $\rho \in A_{\theta}^{\infty}$ and every non-negative
integer $m$, one has
\[\int_0^{\infty}\frac{k^{2m+2}u^m}{(k^2u+1)^{m+1}}\rho
\frac{1}{(k^2u+1)}du=\mathcal{D}_m(\rho),\]
where $\mathcal{D}_m=\mathcal{L}_m(\Delta)$, $\Delta$ is
the modular automorphism introduced in
Section \ref{preliminaries}, and $\mathcal{L}_m$ is the
modified logarithm:
\begin{eqnarray}
\mathcal{L}_m(u)
&=&\int_0^{\infty}\frac{x^m}{(x+1)^{m+1}}\frac{1}{(xu+1)}dx \nonumber \\
&=&(-1)^m(u-1)^{-(m+1)}\big( \log u - \sum_{j=1}^{m}(-1)^{j+1}
\frac{(u-1)^j}{j} \big). \nonumber
\end{eqnarray}
\end{integration}

Now we write $T$ as a sum of three terms $T=T_1+T_2+T_3$ in a way that the
above lemma will be applied to each $T_j$ by choosing $m=j$ for $j=1,2,3.$ We let\\
\begin{eqnarray}
T_1&=&-2|\tau|^2r^3b_0^2k^2\delta_2(k)b_0\delta_2(k)
+2\tau_1r^3b_0^2k^2\delta_2(k)b_0\delta_1(k)- \nonumber \\
&&6\tau_1r^3b_0^2k^2\delta_1(k)b_0\delta_2(k)
-2r^3b_0^2k^2\delta_1(k)b_0\delta_1(k). \nonumber
\end{eqnarray}
With the substitution $u=r^2$, we have $rdr=\frac{1}{2}du$, and using
Lemma \ref{integration}, $\int_{0}^{\infty}\bullet \, dr$ of these
terms (up to an overall factor of $\frac{2\pi}{\tau_2}$) is computed
respectively as follows:
\begin{eqnarray} \label{D_1}
&& -|\tau|^2 k^{-2} \int_{0}^{\infty} \frac{k^4 u}{(k^2u+1)^2}
\delta_2(k) \frac{1}{k^2u+1}\delta_2(k) \, du
+\tau_1 k^{-2} \int_{0}^{\infty} \frac{k^4 u}{(k^2u+1)^2}
\delta_2(k) \frac{1}{k^2u+1}\delta_1(k) \, du  \nonumber \\
&& - 3 \tau_1 k^{-2} \int_{0}^{\infty} \frac{k^4 u}{(k^2u+1)^2}
\delta_1(k) \frac{1}{k^2u+1} \delta_2(k) \, du
-k^{-2} \int_{0}^{\infty} \frac{k^4 u}{(k^2u+1)^2}
\delta_1(k) \frac{1}{k^2u+1}\delta_1(k) \,du  = \nonumber \\
&&  -|\tau|^2 k^{-2} \mathcal{D}_1(\delta_2(k)) \delta_2(k)
+\tau_1 k^{-2} \mathcal{D}_1(\delta_2(k)) \delta_1(k)-
 3 \tau_1 k^{-2} \mathcal{D}_1(\delta_1(k)) \delta_2(k) - k^{-2} \mathcal{D}_1
(\delta_1(k)) \delta_1(k). \nonumber \\
&&
\end{eqnarray}

Then, we let

\begin{eqnarray}
T_2&=&
4|\tau|^2r^5b_0^3k^3\delta_2(k)b_0k\delta_2(k)
+4\tau_1r^5b_0^3k^3\delta_2(k)b_0k\delta_1(k)
+4\tau_1r^5b_0^3k^3\delta_1(k)b_0k\delta_2(k)+ \nonumber \\
&&4r^5b_0^3k^3\delta_1(k)b_0k\delta_1(k)
+4|\tau|^2r^5b_0^3k^4\delta_2(k)b_0\delta_2(k)
+4r^5b_0^3k^4\delta_1(k)b_0\delta_1(k)+ \nonumber \\
&& 8\tau_1r^5b_0^3k^4\delta_1(k)b_0\delta_2(k) , \nonumber
\end{eqnarray}
and compute $\int_{0}^{\infty}\bullet \, dr$ of these
terms respectively (up to an overall factor of $\frac{2\pi}
{\tau_2}$) as follows.
Using $k\delta_j(k)k=k^2 \Delta^{1/2}(\delta_j(k))$, we have
\begin{eqnarray} \label{D_2-1}
 2|\tau|^2  \int_{0}^{\infty} \frac{k^3u^2}{(k^2u+1)^3}  \delta_2(k)
\frac{k}{k^2u+1}  \delta_2(k) \,du &=&
 2 |\tau|^2 k^{-2} \int_{0}^{\infty} \frac{k^6u^2}{(k^2u+1)^3} \Delta^{1/2}
(\delta_2(k)) \frac{1}{k^2u+1} \delta_2(k) du \nonumber \\
&=&  2 |\tau|^2 k^{-2} \mathcal{D}_2 \Delta^{1/2}
(\delta_2(k)) \delta_2(k);
\end{eqnarray}
similarly one has
\begin{eqnarray}\label{D_2-2}
 2\tau_1  \int_{0}^{\infty} \frac{k^3u^2}{(k^2u+1)^3}  \delta_2(k)
\frac{1}{k^2u+1} k \delta_1(k) \,du
= 2 \tau_1 k^{-2} \mathcal{D}_2 \Delta^{1/2} (\delta_2(k)) \delta_1(k),
\end{eqnarray}

\begin{eqnarray}\label{D_2-3}
 2\tau_1  \int_{0}^{\infty} \frac{k^3u^2}{(k^2u+1)^3}  \delta_1(k)
\frac{1}{k^2u+1} k \delta_2(k) \,du
= 2 \tau_1 k^{-2} \mathcal{D}_2 \Delta^{1/2} (\delta_1(k)) \delta_2(k),
\end{eqnarray}

\begin{eqnarray}\label{D_2-4}
 2  \int_{0}^{\infty} \frac{k^3u^2}{(k^2u+1)^3}  \delta_1(k)
\frac{1}{k^2u+1} k \delta_1(k) \,du
= 2  k^{-2} \mathcal{D}_2 \Delta^{1/2} (\delta_1(k)) \delta_1(k).
\end{eqnarray}

Also, we have
\begin{eqnarray}\label{D_2-5}
 2 |\tau|^2  k^{-2} \int_{0}^{\infty} \frac{k^6u^2}{(k^2u+1)^3}
\delta_2(k) \frac{1}{k^2u+1}  \delta_2(k) \,du
= 2  |\tau|^2 k^{-2} \mathcal{D}_2  (\delta_2(k)) \delta_2(k),
\end{eqnarray}

\begin{eqnarray}\label{D_2-6}
 2   k^{-2} \int_{0}^{\infty} \frac{k^6u^2}{(k^2u+1)^3}  \delta_1(k)
\frac{1}{k^2u+1}  \delta_1(k) \,du
= 2  k^{-2} \mathcal{D}_2  (\delta_1(k)) \delta_1(k).
\end{eqnarray}
In order to obtain symmetry in the result, we treat the last term of
$T_2$ more carefully. Since
\[\partial_r(r^4b_0^2)= 4r^3b_0^2-4r^5k^2b_0^3,\]
we can write:
\begin{eqnarray} \label{D_2-7-0}
&&2\tau_1\int_{0}^{\infty}4r^5b_0^3k^4\delta_1(k)b_0\delta_2(k) \,dr =
 2\tau_1\int_{0}^{\infty} (4r^3b_0^2-\partial_r(r^4b_0^2))k^2
\delta_1(k)b_0\delta_2(k) \,dr  \nonumber \\
&& \quad = 2\tau_1\int_{0}^{\infty} 4r^3b_0^2k^2
\delta_1(k)b_0\delta_2(k) \,dr
-2\tau_1\int_{0}^{\infty} \partial_r(r^4b_0^2)k^2
\delta_1(k)b_0\delta_2(k) \,dr \nonumber \\
&& \quad = 4 \tau_1 k^{-2} \int_{0}^{\infty}
\frac{k^4 u}{(k^2u+1)^2} \delta_1(k)
\frac{1}{(k^2u+1)^2} \delta_2(k) \, du + 2\tau_1
\int_{0}^{\infty} r^4b_0^2 k^2
\delta_1(k)\partial_r(b_0) \delta_2(k) \,dr \nonumber \\
&& \quad = 4 \tau_1 k^{-2} \mathcal{D}_1(\delta_1(k))
- 4 \tau_1 \int_{0}^{\infty} r^5b_0^2 k^2
\delta_1(k) k^2 b_0^2 \delta_2(k) \,dr.
\end{eqnarray}
Since we are working under the trace, we can first split the
latter integral to two parts and
then compute it as follows:
\begin{eqnarray}
&&- 4 \tau_1 \int_{0}^{\infty} r^5b_0^2 k^2 \delta_1(k) k^2 b_0^2
\delta_2(k) \,dr \nonumber \\
&& \quad = - 2 \tau_1 \int_{0}^{\infty} r^5b_0^2 k^2 \delta_1(k)
k^2 b_0^2 \delta_2(k) \,dr
- 2 \tau_1 \int_{0}^{\infty} r^5b_0^2 k^2 \delta_2(k) k^2 b_0^2
\delta_1(k) \,dr \nonumber \\
&& \quad = -\tau_1  \int_{0}^{\infty} \partial_r (r^4b_0^2) k^2
\delta_1(k) b_0 \delta_2(k) \,dr
 -\tau_1  \int_{0}^{\infty} \partial_r (r^4b_0^2) k^2 \delta_2(k) b_0
\delta_1(k) \,dr \nonumber
\end{eqnarray}
\begin{eqnarray}
 &=& -\tau_1  \int_{0}^{\infty} (4r^3b_0^2 - 4r^5k^2b_0^3) k^2
\delta_1(k) b_0 \delta_2(k) \,dr -
 \tau_1  \int_{0}^{\infty} (4r^3b_0^2 - 4r^5k^2b_0^3) k^2
\delta_2(k) b_0 \delta_1(k) \,dr \nonumber \\
&=& -4\tau_1  \int_{0}^{\infty} r^3b_0^2  k^2 \delta_1(k)
b_0 \delta_2(k) \,dr
+4\tau_1  \int_{0}^{\infty} r^5k^2b_0^3 k^2 \delta_1(k) b_0 \delta_2(k)
\,dr -\nonumber \\
&&   4\tau_1  \int_{0}^{\infty} r^3b_0^2  k^2 \delta_2(k) b_0
\delta_1(k) \,dr
+4\tau_1  \int_{0}^{\infty} r^5k^2b_0^3 k^2 \delta_2(k) b_0 \delta_1(k)
\,dr \nonumber
\end{eqnarray}
\begin{eqnarray}
&=& -2 \tau_1 k^{-2} \int_{0}^{\infty} \frac{k^4u}{(k^2u+1)^2} \delta_1(k)
\frac{1}{k^2u+1} \delta_2(k) \, du + \nonumber \\
&&  2 \tau_1 k^{-2} \int_{0}^{\infty} \frac{k^6u^2}{(k^2u+1)^3} \delta_1(k)
\frac{1}{k^2u+1} \delta_2(k) \, du - \nonumber \\
&&  2 \tau_1 k^{-2} \int_{0}^{\infty} \frac{k^4u}{(k^2u+1)^2} \delta_2(k)
\frac{1}{k^2u+1} \delta_1(k) \, du + \nonumber \\
&&  2 \tau_1 k^{-2} \int_{0}^{\infty} \frac{k^6u^2}{(k^2u+1)^3} \delta_2(k)
\frac{1}{k^2u+1} \delta_1(k) \, du \nonumber
\end{eqnarray}
\begin{eqnarray}
&=& -2\tau_1 k^{-2} \mathcal{D}_1(\delta_1(k)) \delta_2(k) +
2\tau_1 k^{-2} \mathcal{D}_2(\delta_1(k)) \delta_2(k) - \nonumber \\
&&   2\tau_1 k^{-2} \mathcal{D}_1(\delta_2(k)) \delta_1(k) +
2\tau_1 k^{-2} \mathcal{D}_2(\delta_2(k)) \delta_1(k).  \nonumber
\end{eqnarray}
Therefore from \eqref{D_2-7-0} we have
\begin{eqnarray}\label{D_2-7}
8\tau_1 &\int_{0}^{\infty}& r^5b_0^3k^4\delta_1(k)b_0\delta_2(k) \,dr  \nonumber \\
&=&  2\tau_1 k^{-2} \mathcal{D}_1(\delta_1(k)) \delta_2(k) +
2\tau_1 k^{-2} \mathcal{D}_2(\delta_1(k)) \delta_2(k) -
   \nonumber \\
&& 2\tau_1 k^{-2} \mathcal{D}_1(\delta_2(k)) \delta_1(k) +
2\tau_1 k^{-2} \mathcal{D}_2(\delta_2(k)) \delta_1(k).
\end{eqnarray}

Now, we let

\begin{eqnarray}
T_3&=&-2|\tau|^2r^7b_0^4k^4\delta_2(k)b_0k^2\delta_2(k)
-2\tau_1r^7b_0^4k^4\delta_2(k)b_0k^2\delta_1(k)
-2\tau_1r^7b_0^4k^4\delta_1(k)b_0k^2\delta_2(k)- \nonumber \\
&&2r^7b_0^4k^4\delta_1(k)b_0k^2\delta_1(k)
-4|\tau|^2r^7b_0^4k^5\delta_2(k)b_0k\delta_2(k)
-4\tau_1r^7b_0^4k^5\delta_2(k)b_0k\delta_1(k)- \nonumber \\
&&4\tau_1r^7b_0^4k^5\delta_1(k)b_0k\delta_2(k)
-4r^7b_0^4k^5\delta_1(k)b_0k\delta_1(k)
-2|\tau|^2r^7b_0^4k^6\delta_2(k)b_0\delta_2(k)- \nonumber \\
&&2\tau_1r^7b_0^4k^6\delta_2(k)b_0\delta_1(k)
-2\tau_1r^7b_0^4k^6\delta_1(k)b_0\delta_2(k)
-2r^7b_0^4k^6\delta_1(k)b_0\delta_1(k) \nonumber
\end{eqnarray}
and compute $\int_{0}^{\infty}\bullet \, dr$ of these
terms respectively (up to an overall factor of $\frac{2\pi}
{\tau_2}$) as follows. Using $\delta_j(k)k^2=k^2\Delta(\delta_j(k))$,
we have
\begin{eqnarray}\label{D_3-1}
 -|\tau|^2 \int_{0}^{\infty} \frac{k^4u^3}{(k^2u+1)^4}\delta_2(k)\frac{k^2}{k^2u+1}
\delta_2(k)du&=&
-|\tau|^2 k^{-2} \int_{0}^{\infty} \frac{k^8u^3}{(k^2u+1)^4}\Delta(\delta_2(k))
\frac{1}{k^2u+1}\delta_2(k)\,du \nonumber \\
&=&-|\tau|^2 k^{-2} \mathcal{D}_3\Delta(\delta_2(k))\delta_2(k);
\end{eqnarray}
similarly one has
\begin{eqnarray}\label{D_3-2}
-\tau_1 \int_{0}^{\infty} \frac{k^4u^3}{(k^2u+1)^4} (\delta_2(k))\frac{1}
{k^2u+1}k^2 \delta_1(k) \, du=
-\tau_1 k^{-2} \mathcal{D}_3\Delta(\delta_2(k)) \delta_1(k),
\end{eqnarray}
\begin{eqnarray}\label{D_3-3}
-\tau_1 \int_{0}^{\infty} \frac{k^4u^3}{(k^2u+1)^4} (\delta_1(k))\frac{1}
{k^2u+1}k^2 \delta_2(k) \, du=
-\tau_1 k^{-2} \mathcal{D}_3\Delta(\delta_1(k)) \delta_2(k),
\end{eqnarray}
\begin{eqnarray}\label{D_3-4}
- \int_{0}^{\infty} \frac{k^4u^3}{(k^2u+1)^4} (\delta_1(k))\frac{1}
{k^2u+1}k^2 \delta_1(k) \, du=
- k^{-2} \mathcal{D}_3\Delta(\delta_1(k)) \delta_1(k).
\end{eqnarray}
Also, using $k\delta_j(k)k=k^2 \Delta^{1/2}(\delta_j(k))$, we have
\begin{eqnarray}\label{D_3-5}
-4|\tau|^2&\int_{0}^{\infty}&\frac{k^5u^3}{(k^2u+1)^4}\delta_2(k)\frac{k}
{k^2u+1}\delta_2(k) \, du   \nonumber \\
&=& -2|\tau|^2 k^{-2} \int_{0}^{\infty}\frac{k^6u^3}{(k^2u+1)^4}\Delta^{1/2}
(\delta_2(k))
\frac{1}{k^2u+1}\delta_2(k) \, du \nonumber \\
&=&-2|\tau|^2 k^{-2} \mathcal{D}_3\Delta^{1/2}(\delta_2(k))\delta_2(k);
\end{eqnarray}
similarly one has
\begin{eqnarray}\label{D_3-6}
-2\tau_1 \int_{0}^{\infty}\frac{k^5u^3}{(k^2u+1)^4}\delta_2(k)
\frac{1}{k^2u+1}k\delta_1(k) \,du=
-2\tau_1 k^{-2} \mathcal{D}_3\Delta^{1/2}(\delta_2(k))\delta_1(k),
\end{eqnarray}

\begin{eqnarray}\label{D_3-7}
-2\tau_1 \int_{0}^{\infty}\frac{k^5u^3}{(k^2u+1)^4}\delta_1(k)
\frac{1}{k^2u+1}k\delta_2(k) \,du=
-2\tau_1 k^{-2} \mathcal{D}_3\Delta^{1/2}(\delta_1(k))\delta_2(k),
\end{eqnarray}

\begin{eqnarray}\label{D_3-8}
-2 \int_{0}^{\infty}\frac{k^5u^3}{(k^2u+1)^4}\delta_1(k)
\frac{1}{k^2u+1}k\delta_1(k) \,du=
-2 k^{-2} \mathcal{D}_3\Delta^{1/2}(\delta_1(k))\delta_1(k).
\end{eqnarray}
For the rest of the terms in $T_3$, we have
\begin{eqnarray}\label{D_3-9}
-|\tau|^2 \int_{0}^{\infty}\frac{k^6u^3}{(k^2u+1)^4}\delta_2(k)
\frac{1}{k^2u+1}k\delta_2(k) \,du=
-|\tau|^2 k^{-2} \mathcal{D}_3(\delta_2(k))\delta_2(k),
\end{eqnarray}

\begin{eqnarray}\label{D_3-10}
-\tau_1 \int_{0}^{\infty}\frac{k^6u^3}{(k^2u+1)^4}\delta_2(k)
\frac{1}{k^2u+1}k\delta_1(k) \,du=
-\tau_1 k^{-2} \mathcal{D}_3(\delta_2(k))\delta_1(k),
\end{eqnarray}

\begin{eqnarray}\label{D_3-11}
-\tau_1 \int_{0}^{\infty}\frac{k^6u^3}{(k^2u+1)^4}\delta_1(k)
\frac{1}{k^2u+1}k\delta_2(k) \,du=
-\tau_1 k^{-2} \mathcal{D}_3(\delta_1(k))\delta_2(k),
\end{eqnarray}

\begin{eqnarray}\label{D_3-12}
- \int_{0}^{\infty}\frac{k^6u^3}{(k^2u+1)^4}\delta_1(k)
\frac{1}{k^2u+1}k\delta_1(k) \,du=
- k^{-2} \mathcal{D}_3(\delta_1(k))\delta_1(k).
\end{eqnarray}

\section{Independence of $\zeta(0)$ from the Weyl Factor}

The calculations of Section \ref{computationofzeta(0)}
allow us to express the value at the origin of the zeta function of
the operator $\triangle'\sim k\triangle k$ by a complex formula:
\begin{eqnarray} \label{premainlemma}
\zeta(0)+1 &=& \frac{-2 \pi}{\tau_2}\tau_0 \Big ( -\frac{1}{6}k^{-1}\delta_1^2(k)-\frac{|\tau|^2}
{6}k^{-1}\delta_2^2(k)-
\frac{\tau_1}{3}k^{-1}\delta_1\delta_2(k)+\frac{1}{3}k^{-2}\delta_1(k)^2
+\nonumber \\
&&\frac{|\tau|^2}{3}k^{-2}\delta_2(k)^2
 +\frac{\tau_1}{3}k^{-2}\delta_1(k)\delta_2(k)
+\frac{\tau_1}{3}k^{-2}\delta_2(k)\delta_1(k)
-|\tau|^2 k^{-2} \mathcal{D}_1(\delta_2(k)) \delta_2(k) + \nonumber \\
&& \tau_1 k^{-2} \mathcal{D}_1(\delta_2(k)) \delta_1(k)-
3 \tau_1 k^{-2} \mathcal{D}_1(\delta_1(k)) \delta_2(k) - k^{-2} \mathcal{D}_1
(\delta_1(k)) \delta_1(k)  + \nonumber \\
&&2 |\tau|^2 k^{-2} \mathcal{D}_2 \Delta^{1/2} (\delta_2(k)) \delta_2(k)
+ 2 \tau_1 k^{-2} \mathcal{D}_2 \Delta^{1/2} (\delta_2(k)) \delta_1(k)+ \nonumber \\
&&2 \tau_1 k^{-2} \mathcal{D}_2 \Delta^{1/2} (\delta_1(k)) \delta_2(k)
+ 2  k^{-2} \mathcal{D}_2 \Delta^{1/2} (\delta_1(k)) \delta_1(k) + \nonumber \\
&&2  |\tau|^2 k^{-2} \mathcal{D}_2  (\delta_2(k)) \delta_2(k) +
2  k^{-2} \mathcal{D}_2  (\delta_1(k)) \delta_1(k) +
2\tau_1 k^{-2} \mathcal{D}_1(\delta_1(k)) \delta_2(k) + \nonumber \\
&&2\tau_1 k^{-2} \mathcal{D}_2(\delta_1(k)) \delta_2(k) -
2\tau_1 k^{-2} \mathcal{D}_1(\delta_2(k)) \delta_1(k) +
2\tau_1 k^{-2} \mathcal{D}_2(\delta_2(k)) \delta_1(k) \nonumber \\
&& -|\tau|^2 k^{-2} \mathcal{D}_3\Delta(\delta_2(k))\delta_2(k)
-\tau_1 k^{-2} \mathcal{D}_3\Delta(\delta_2(k)) \delta_1(k)
-\tau_1 k^{-2} \mathcal{D}_3\Delta(\delta_1(k)) \delta_2(k) - \nonumber \\
&&k^{-2} \mathcal{D}_3\Delta(\delta_1(k)) \delta_1(k)
-2|\tau|^2 k^{-2} \mathcal{D}_3\Delta^{1/2}(\delta_2(k))\delta_2(k)-
2\tau_1 k^{-2} \mathcal{D}_3\Delta^{1/2}(\delta_2(k))\delta_1(k)- \nonumber \\
&&2\tau_1 k^{-2} \mathcal{D}_3\Delta^{1/2}(\delta_1(k))\delta_2(k)
-2 k^{-2} \mathcal{D}_3\Delta^{1/2}(\delta_1(k))\delta_1(k)
-|\tau|^2 k^{-2} \mathcal{D}_3(\delta_2(k))\delta_2(k) - \nonumber \\
&&\tau_1 k^{-2} \mathcal{D}_3(\delta_2(k))\delta_1(k)
-\tau_1 k^{-2} \mathcal{D}_3(\delta_1(k))\delta_2(k)
- k^{-2} \mathcal{D}_3(\delta_1(k))\delta_1(k)
\Big ).
\end{eqnarray}
Note that we have considered an overall change of sign
since we set $\lambda=-1$ in our calculations in Section \ref{computationofzeta(0)}.\\

In the following lemma, similar to the work of Connes and Tretkoff \cite{contre},
we modify \eqref{premainlemma} to prepare the ground for the proof of the main result
of this paper which relies on using functional calculus on the modular automorphism.

\newtheorem{mainlemma}[main1]{Lemma}
\begin{mainlemma} \label{mainlemma}
For an irrational number $\theta$, let $k$ be an invertible positive element of
$A_{\theta}^{\infty}$. Then the value
$\zeta(0)$ of the zeta function $\zeta$ of the operator
$\triangle' \sim k \triangle k$ is given
by
\begin{eqnarray} \label{zeta(0)+1}
\zeta(0)+1=
&& \frac{2\pi}{\tau_2} \varphi (f({\Delta})(\delta_1(k))\delta_1(k))
+ \frac{2\pi|\tau|^2}{\tau_2} \varphi (f(\Delta)(\delta_2(k))\delta_2(k))
                         \nonumber \\
&& + \frac{2 \pi \tau_1}{\tau_2} \varphi(f(\Delta)(\delta_1(k)) \delta_2(k) )
+ \frac{2 \pi \tau_1}{\tau_2} \varphi(f(\Delta)(\delta_2(k)) \delta_1(k) ),
\end{eqnarray}
where $\varphi(x)=\tau_0(xk^{-2})$, $\tau_0$ is the unique trace on
$A_{\theta}$, $\Delta$ is the modular automorphism, and
\[f(u)=\frac{1}{6}u^{-1/2}-\frac{1}{3}+\mathcal{L}_1(u)-2(1+u^{1/2})
\mathcal{L}_2(u)+(1+u^{1/2})^2\mathcal{L}_3(u),\nonumber\]
where for each positive integer $m$, $\mathcal{L}_m$ is the modified logarithm:
\[\mathcal{L}_m(u)=(-1)^m(u-1)^{-(m+1)}\big( \log u -
\sum_{j=1}^{m}(-1)^{j+1}\frac{(u-1)^j}{j} \big).\nonumber \]
\begin{proof}
Under the trace, one can replace the first three terms in
\eqref{premainlemma}, namely
$-\frac{1}{6}k^{-1}\delta_1^2(k)-\frac{|\tau|^2}{6}k^{-1}\delta_2^2(k)-
\frac{\tau_1}{3}k^{-1}\delta_1\delta_2(k)$
by
\begin{eqnarray}
&&-\frac{1}{6}k^{-2}\Delta^{-1/2}(\delta_1(k))\delta_1(k)-
\frac{|\tau|^2}{6}k^{-2}\Delta^{-1/2}(\delta_2(k))\delta_2(k)- \nonumber \\
&&\frac{\tau_1}{6}k^{-2}\Delta^{-1/2}(\delta_1(k))\delta_2(k)-
\frac{\tau_1}{6}k^{-2}\Delta^{-1/2}(\delta_2(k))\delta_1(k), \nonumber
\end{eqnarray}
because for $i,j=1,2$ one has
\begin{eqnarray}
\tau_0(k^{-1}\delta_i\delta_j(k))&=&-\tau_0(\delta_i(k^{-1})\delta_j(k))=
\tau_0(k^{-1} \delta_i(k) k^{-1} \delta_j(k))=\tau_0(k^{-2} k \delta_i(k) k^{-1}
\delta_j(k)) \nonumber \\
&=& \tau_0(k^{-2} \Delta^{-1/2} (\delta_i(k))  \delta_j(k)). \nonumber
\end{eqnarray}
After this replacement, the terms of \eqref{premainlemma}
 are easily compared with the expressions in the
statement of this lemma.
\end{proof}
\end{mainlemma}

In the following theorem, we show that the value at the origin of the zeta function
of the operator $\triangle' \sim k \triangle k$ is independent of the choice of the
invertible positive element $k\in A_{\theta}^{\infty}$ considered as the Weyl
factor by which we changed the metric within its conformal class. Thus, for an irrational
number $\theta$, we have
proved the analogue of the Gauss-Bonnet theorem for the noncommutative two
torus  $A_{\theta}$ endowed with the conformal structure associated
with a complex number $\tau \in \mathbb{C} \setminus \mathbb{R}$.

\newtheorem{maintheorem}[main1]{Theorem}
\begin{maintheorem} \label{maintheorem}
For an irrational number $\theta$, let $k$ be an invertible positive element of
$A_{\theta}^{\infty}$. Then the value
$\zeta(0)$ of the zeta function $\zeta$ of the operator
$\triangle' \sim k \triangle k$ is independent of $k$.
\begin{proof}
The function $f(u)$ in the statement of Lemma \ref{mainlemma} is of the form
$h(\log u)$ where $h$ is the following entire function \cite{contre}:
\begin{equation}
h(x)=-\frac{e^{-x/2}(-1+3e^{x/2}+3e^x+6e^{3x/2}x-3e^{2x}-3e^{5x/2}+e^{3x})}
{6(-1+e^{x/2})^4(1+e^{x/2})^2}.
\end{equation}

Now, using the trace property of $\tau_0$, the fact that left multiplication by
$k^{-1}$ commutes
with any function of $\Delta$, and the following identities ($cf.$ \cite{contre})
\begin{equation} k^{-1} \delta_j (k) = 2 \frac{\Delta^{1/2}-1}{\log \Delta}
(\delta_j(\log k)), \,\,\, j=1,2,
\end{equation}
\begin{equation}  \delta_j (k)k^{-1} = -2 \frac{\Delta^{-1/2}-1}{\log \Delta}
(\delta_j(\log k)), \,\,\, j=1,2,
\end{equation}
\begin{equation} \label{logDeltaproperty}
\tau_0(a F(\log \Delta)(b)) = \tau_0 (F(- \log \Delta)(a)b), \,\,\, \forall
a,b \in A_{\theta}^{\infty},
\end{equation}
where $F$ is an entire function, for $i,j=1,2$, we have
\begin{eqnarray} \label{mainformula}
\varphi(f(\Delta)(\delta_i(k)) \delta_j(k)) &=& \tau_0 (f(\Delta)
(\delta_i(k)) \delta_j(k)k^{-2}) \nonumber \\
&=&\tau_0 (f(\Delta)(k^{-1}\delta_i(k)) \delta_j(k)k^{-1}) \nonumber \\
&=& \tau_0 (h(\log \Delta)(k^{-1}\delta_i(k)) \delta_j(k)k^{-1}) \nonumber \\
&=& \tau_0 \Big ( h\big(\log \Delta \big) \big(2 \frac{\Delta^{1/2}-1}
{\log \Delta} (\delta_i(\log k)) \big ) \big ( -2 \frac{\Delta^{-1/2}-1}
{\log \Delta} (\delta_j(\log k))\big )\Big ) \nonumber \\
&=& \tau_0 (K(\log \Delta) (\delta_i(\log k)) \delta_j(\log k)),
\end{eqnarray}
where the function $K$ (using \eqref{logDeltaproperty}) is given by
\begin{equation}
K(x)=-2\frac{e^{x/2}-1}{-x} h(x) 2\frac{e^{x/2}-1}{x}=4x^{-2}(e^{x/2}-1)^2h(x).
\end{equation}
In fact \cite{contre}:
\begin{equation}
K(x)= -\frac{x- \text{sh}\left[\frac{x}{2}\right]- \text{sh}[x]+\frac
13\text{sh}\left[\frac{3 x}{2}\right]}{
x^2\text{sh}\left[\frac{x}{2}\right]^2}.
\end{equation}
In particular, $K$ is an odd function. Thus, using the trace property of $\tau_0$
and \eqref{logDeltaproperty}, from \eqref{mainformula} one has
\begin{eqnarray}
\varphi(f(\Delta)(\delta_i(k)) \delta_j(k))&=& \tau_0 (K(\log \Delta) (\delta_i(\log k))
\delta_j(\log k)) \nonumber \\
&=& \tau_0 ( \delta_j(\log k) K(\log \Delta) (\delta_i(\log k))) \nonumber \\
&=& \tau_0 ( K(-\log \Delta)(\delta_j(\log k))  \delta_i(\log k)) \nonumber \\
&=& - \tau_0 ( K(\log \Delta)(\delta_j(\log k))  \delta_i(\log k)) \nonumber \\
&=& - \varphi(f(\Delta)(\delta_j(k)) \delta_i(k)).
\end{eqnarray}
Therefore $
\varphi(f(\Delta)(\delta_j(k)) \delta_j(k)) = 0$ for $j=1,2$ and
$\varphi(f(\Delta)(\delta_1(k)) \delta_2(k)) = -\varphi(f(\Delta)(\delta_2(k)) \delta_1(k))$.
Now, the desired independence follows from \eqref{zeta(0)+1}:\begin{eqnarray}
\zeta(0)+1&=&
 \frac{2\pi}{\tau_2} \varphi (f({\Delta})(\delta_1(k))\delta_1(k))
+ \frac{2\pi|\tau|^2}{\tau_2} \varphi (f(\Delta)(\delta_2(k))\delta_2(k))
                         \nonumber \\
&& + \frac{2 \pi \tau_1}{\tau_2} \varphi(f(\Delta)(\delta_1(k)) \delta_2(k) )
+ \frac{2 \pi \tau_1}{\tau_2} \varphi(f(\Delta)(\delta_2(k)) \delta_1(k) ) \nonumber \\
&=& 0.
\end{eqnarray}
\end{proof}
\end{maintheorem}

\end{document}